\documentclass{ucthesis}
\usepackage{amssymb,amsmath,amsthm,amsbsy}
\usepackage[all]{xy}

\newtheorem{thm}{Theorem}[section]
\newtheorem{prop}[thm]{Proposition} 
\newtheorem{lemma}[thm]{Lemma}
\newtheorem{cor}[thm]{Corollary}

\theoremstyle{definition} 
\newtheorem{dfn}[thm]{Definition}

\theoremstyle{remark}
\newtheorem{rmk}[thm]{Remark}
\newtheorem{ex}[thm]{Example}

\newcommand{\la}[4]{
\xymatrix{#1 \ar[r] \ar@<2pt>[d] \ar@<-2pt>[d] & #2 \ar@<2pt>[d] \ar@<-2pt>[d] \\
#3 \ar[r] & #4}}
\newcommand{\dvb}[4]{
\xymatrix{#1 \ar[r] \ar[d] & #2 \ar[d] \\
#3 \ar[r] & #4}}

\newcommand{\pdiff}[2]{\frac{\partial #1}{\partial #2}}

\newcommand{\arrows}{\rightrightarrows}
\newcommand{\defequal}{\stackrel{\mbox {\tiny {def}}}{=}}

\newcommand{\reals}{{\mathbb R}}

\newcommand{\integers}{{\mathbb Z}}
\newcommand{\sgn}{{\rm sgn~}}
\newcommand{\Hom}{{\rm Hom}}

\newcommand{\bitimes}[2]{\,_{#1}\!\times_{#2}}

\newcommand{\calO}{\mathcal{O}}
\newcommand{\lie}{\mathcal{L}}
\newcommand{\livf}{\mathcal{X}_{LI}}
\newcommand{\vect}{\mathcal{X}}
\DeclareMathOperator{\sym}{S}
\newcommand{\bigS}{\mbox{{\large{$\sym$}}}}
\renewcommand{\bigwedge}{\mbox{{\Large{$\wedge$}}}}

\DeclareMathOperator{\im}{im}
\DeclareMathOperator{\ev}{ev}
\DeclareMathOperator{\rank}{rank}
\DeclareMathOperator{\ad}{ad}
\DeclareMathOperator{\Ad}{Ad}
\DeclareMathOperator{\Ob}{Ob}
\DeclareMathOperator{\Mor}{Mor}
\DeclareMathOperator{\Lie}{Lie}

\begin{document}
\title{Supergroupoids, double structures, and equivariant cohomology}
\author{Rajan Amit Mehta}
\degreesemester{Spring} \degreeyear{2006} \degree{Doctor of
Philosophy} \chair{Professor Alan Weinstein}
\othermembers{Professor Robion Kirby \\
  Professor Kam-Biu Luk }
\numberofmembers{3} \prevdegrees{B.A. (University of Pennsylvania) 1999} \field{Mathematics} \campus{Berkeley}
    \maketitle

     \begin{abstract}
$Q$-groupoids and $Q$-algebroids are, respectively, supergroupoids and superalgebroids that are equipped with compatible homological vector fields.  These new objects are closely related to the double structures of Mackenzie; in particular, we show that $Q$-groupoids are intermediary objects between Mackenzie's $\mathcal{LA}$-groupoids and double complexes, which include as a special case the simplicial model of equivariant cohomology.  There is also a double complex associated to a $Q$-algebroid, which in the above special case is the BRST model of equivariant cohomology.  Other special cases include models for the Drinfel'd double of a Lie bialgebra and Ginzburg's equivariant Poisson cohomology.  Finally, a supergroupoid version of the van Est map is used to give a homomorphism from the double complex of a $Q$-groupoid to that of a $Q$-algebroid.
    \abstractsignature
    \end{abstract}
    \begin{frontmatter}

\begin{dedication}\null\vfil{\large\begin{center}{\em \`a ma lionne Vanitha}\end{center}}\vfil\null\end{dedication}

\tableofcontents 

\begin{acknowledgements} 
Above all, I wish to thank my advisor Alan Weinstein for several years of guidance, encouragement, and support.  I am constantly impressed by his vision and intuition for mathematics as well as his ability to communicate ideas.  I can only hope that some of those skills have rubbed off onto me.

I would like to thank Alfonso Gracia-Saz for many useful comments and helpful discussions, as well as Camille Laurent-Gengoux, Eli Lebow, David Spivak, Xiang Tang, Marco Zambon, and Chenchang Zhu for interesting discussions on topics relating to this work.  Thanks to Bonnie Huggins, James Kelley, Walter Kim, Sasha Peterka, and Aaditya Rangan for helping me to navigate the labryinths of Evans Hall.  

Thanks to my parents and to my cousin Kamlesh and his family for their support.  Finally, I want to give a shout-out my wonderful friends in the Bay Area and around the world.
\end{acknowledgements}
		\end{frontmatter}

    \chapter{Introduction}
The concept of classifying space of a (topological) category $\mathcal{C}$ dates back to Segal \cite{segal}, who defined the classifying space $B\mathcal{C}$ to be the geometric realization of a simplicial space
\begin{equation}\label{eqn:nerveint}
\xymatrix{
\cdots \ar@<6pt>[r]\ar@<2pt>[r]\ar@<-2pt>[r]\ar@<-6pt>[r] & \Mor_2 \mathcal{C} \ar@<4pt>[r]\ar[r]\ar@<-4pt>[r] & \Mor \mathcal{C} \ar@<2pt>[r]\ar@<-2pt>[r] & \Ob \mathcal{C},
}
\end{equation}
where $\Mor_q \mathcal{C}$ is the space of composable $q$-tuplets of morphisms.  This construction essentially generalized Milnor's \cite{milnor:universal1} \cite{milnor:universal2} construction of the classifying space of a topological group $\Gamma$, which may be viewed as a category with a single object and for which all morphisms are invertible.

In the case where $\mathcal{C} = G$ is a topological groupoid, Buffet and Lor \cite{buffet} have shown that $BG$ classifies the homotopy classes of principal $G$-bundles\footnote{A similar result for general categories has been obtained by Moerdijk \cite{moerdijk}; also see \cite{weiss}.}.  In light of their result, the cohomology ring $H^\bullet(BG)$ may be viewed as the ring of characteristic classes for principal $G$-bundles.  In certain special cases, there are even more tangible interpretations.
\begin{itemize}
	\item If $G = X \times \Gamma$ is the action groupoid for the action of a topological group $\Gamma$ on a topological space $X$, then $H^\bullet(BG)$ is equal to the equivariant cohomology $H_\Gamma^\bullet(X)$ (see Appendix \ref{appendix}).
	\item If $G$ is an \'{e}tale groupoid representing an orbifold $X$, then $H^\bullet(BG)$ is equal to the orbifold sheaf cohomology \cite{mp:orb} \cite{mp:orb2}.
	\item If $G = \Gamma^q$ is the Haefliger groupoid \cite{bott-haefliger}, then $H^\bullet(BG)$ consists of the ``universal characteristic classes'' for codimension $q$ foliations.  One may also consider the holonomy groupoid of a foliation \cite{crainic-moerdijk}.
\end{itemize}

Let $C^{p,q} \defequal S^p(G^{(q)})$ be the space of singular cochains on the space of composable $q$-tuplets of elements of $G$.  Then there are operators $d: C^{p,q} \to C^{p+1,q}$, the singular coboundary operator, and $\delta: C^{p,q} \to C^{p, q+1}$, the coboundary operator associated to the simplicial structure (\ref{eqn:nerveint}).  It is a well-known result (see, e.g. \cite{bss}, \cite{dupont}) that the total cohomology of the double complex $(C^{p,q};d,\delta)$ computes $H^\bullet(BG)$.  This double complex was employed by Bott, Shulman, and Stasheff \cite{bott} \cite{bss} \cite{shulman} in their study of classifying spaces and foliations.  Furthermore, if $G$ is a Lie groupoid, then differential forms may be used in the place of singular cochains \cite{bss}, resulting in the \emph{de Rham double complex} of $G$ \cite{ltx}.

The main objective of this thesis is to investigate the relationship between double complexes such as the de Rham double complex and double structures such as those introduced by Mackenzie in \cite{mac:dblie1} and \cite{mac:dblie2}.  These are double structures in the sense of Ehresmann and include the following examples:
\begin{itemize}
\item double Lie groupoids, or groupoid objects in the category of Lie groupoids,
\item double Lie algebroids, or Lie algebroid objects in the category of Lie algebroids\footnote{Because the usual definition of a Lie algebroid is not completely expressed in terms of morphisms (in particular, the bracket is not a bundle map), it is a nontrivial issue to define double Lie algebroids.  Mackenzie has given a definition in \cite{mac:dbl2} (also see \cite{mac:dbl}).  There is an equivalent homological definition that is the subject of work in progress with A. Gracia-Saz \cite{gm}.}, 
\item $\mathcal{LA}$-groupoids, or groupoid objects in the category of Lie algebroids.
\end{itemize}
We collectively refer to these three types of objects as \emph{Mackenzie doubles}.  In \cite{mac:dblie1} and \cite{mac:dblie2}, Mackenzie has dealt extensively with the Lie functor and its application to these structures.  In particular, there are two Lie functors (a horizontal one and a vertical one) from the category of double Lie groupoids to the category of $\mathcal{LA}$-groupoids, and there is a Lie functor from the category of $\mathcal{LA}$-groupoids to the category of double Lie algebroids.

In this work, we introduce two new double structures that are closely related to the Mackenzie doubles.  These structures, called \emph{$Q$-groupoids} and \emph{$Q$-algebroids}, are groupoid and Lie algebroid structures, respectively, in the category of $Q$-manifolds\footnote{The term \emph{$Q$-manifold} is due to Schwarz \cite{schwarz}; a $Q$-manifold is a supermanifold equipped with a homological vector field.  See \S\ref{sub:qman}.}.  More concretely, a $Q$-groupoid (resp. $Q$-algebroid) is a supergroupoid (resp. superalgebroid) equipped with a homological vector field satisfying certain compatibility conditions.

There is a well-known functor \cite{aksz} \cite{vaintrob}, which we denote $[-1]$, from the category of Lie algebroids to the category of $Q$-manifolds, that takes a Lie algebroid $A \to M$ to the $Q$-manifold $([-1]A, d_A)$, where $[-1]A$ is the supermanifold with structure sheaf $\Gamma\left(\bigwedge A^*\right)$ and $d_A$ is the differential associated to the algebroid structure of $A$.  The $[-1]$ functor may be applied to the Mackenzie doubles in the following ways.  There are two $[-1]$ functors from the category of double Lie algebroids to the category of $Q$-algebroids, and there is a $[-1]$ functor from the category of $\mathcal{LA}$-groupoids to the category of $Q$-groupoids.  Finally, there is a $[-1]$ functor from the category of $Q$-algebroids to the category of \emph{double $Q$-manifolds}, or supermanifolds equipped with two compatible homological vector fields.  The diagram in Table \ref{table:qdoubles} describes the various categories of double structures and the functors between them.

\begin{table}
\caption{Functors between the categories of Mackenzie doubles, $Q$-groupoids, and $Q$-algebroids}
\label{table:qdoubles}
\begin{equation*}
\xymatrix{
\mbox{Double Lie groupoids} \ar^{\Lie_H}[r] \ar^{\Lie_V}[d] & \mbox{$\mathcal{LA}$-groupoids} \ar^{[-1]}[r] \ar^{\Lie}[d] & \mbox{$Q$-groupoids} \ar^{\Lie}[d] \\
\mbox{$\mathcal{LA}$-groupoids} \ar^{\Lie}[r] \ar^{[-1]}[d] & \mbox{Double Lie algebroids} \ar^{[-1]_H}[r] \ar^{[-1]_V}[d] & \mbox{$Q$-algebroids} \ar^{[-1]}[d] \\
\mbox{$Q$-groupoids} \ar^{\Lie}[r] & \mbox{$Q$-algebroids} \ar^{[-1]}[r] & \mbox{Double $Q$-manifolds}
}
\end{equation*}
\end{table}

Our key observation is that, if $\mathcal{G}$ is a $Q$-groupoid with homological vector field $\psi$, then the Eilenberg-Maclane complex $\left(C^\infty(\mathcal{G}^{(q)}), \delta\right)$ naturally extends to a double complex $\left(C^\infty_p (\mathcal{G}^{(q)}), \delta, \psi\right)$, where $C^\infty_p (\mathcal{G}^{(q)})$ is the space of degree $p$ functions on $\mathcal{G}^{(q)}$.  In other words, there is a (contravariant) functor from the category of $Q$-groupoids to the category of double complexes.  If this functor is composed with the $[-1]$ functor, we obtain a functor from the category of $\mathcal{LA}$-groupoids to the category of double complexes.  

Two examples \cite{cdw}  \cite{mac-xu:bialg} of $\mathcal{LA}$-groupoids are the \emph{tangent prolongation} $TG \arrows TM$ of a groupoid $G \arrows M$ and the \emph{cotangent prolongation} $T^*G \arrows A(G)^*$ of a Poisson groupoid, where $A(G) \to M$ is the Lie algebroid of $G$.  The associated $Q$-groupoids are, respectively, $([-1]TG, d)$, where $d$ is the de Rham differential of $G$, and $([-1]T^*G, d_\pi)$, where $d_\pi$ is the Lichnerowicz \cite{lich} differential for Poisson cohomology.  The double complex that arises from $[-1]TG$ is equal to the de Rham double complex of $G$.  In the case of $[-1]T^*G$, one obtains a double complex that computes the Poisson-invariant groupoid cohomology of $G$.

By applying the Lie functor to a $Q$-groupoid, we may obtain a $Q$-algebroid;  for example, the $Q$-algebroid of $[-1]TG$ is $[-1]TA$, where the homological vector field is the de Rham differential of $A$ (see \S\ref{sub:lieonetg}).  The $Q$-algebroid of $[-1]T^*G$ is $[-1]T^*A$, where the homological vector field is the Lichnerowicz differential for the induced Poisson structure on $A$.  By applying the $[-1]$ functor to a $Q$-algebroid, we obtain a \emph{double $Q$-manifold}, or a supermanifold equipped with two commuting homological vector fields.  This produces another double complex, which in a sense contains the infinitesimal data of the double complex that we have obtained at the groupoid level.  

The relationship between the two double complexes is illustrated (and in fact, our approach is inspired by) the case of an action groupoid $G = M \times \Gamma$, where $\Gamma$ is a Lie group\footnote{Since we will be using infinitesimal models, we will assume that $\Gamma$ is connected and compact.} acting on $M$.  In this case, the Lie algebroid of $G$ is $A = M \times \mathfrak{g}$, where $\mathfrak{g}$ is the Lie algebra of $\Gamma$.  The double complex that arises from $[-1]TG$ is equal to the simplicial model of equivariant cohomology \cite{ltx}, and the double complex that arises from $[-1]TA$ is equal to the BRST model \cite{kalkman}, which also is used to compute equivariant cohomology.  However, the BRST model is only a model for $M \times E\Gamma$, and in order to obtain the equivariant cohomology one must restrict to a basic subcomplex.  The basic subcomplex of the BRST model is known to be equivalent to the Cartan model \cite{cartan}.

We consider the double complex of $[-1]TA$, where $A$ is an arbitrary Lie algebroid, as a ``generalized BRST model''.  We extend to this case an isomorphism of Mathai-Quillen \cite{mq} and Kalkman \cite{kalkman} that relates the BRST model to the Weil model\footnote{Our description of the isomorphism is coordinate-free and therefore may provide insight even into the cases considered by Mathai-Quillen and Kalkman.}.  As in our motivating example, a basic subcomplex is required in order to obtain interesting cohomology.  We propose a definition of a basic subcomplex, requiring a choice of a connection on $A$.  In the case of an action algebroid, there is a canonical connection and the basic subcomplex agrees with that of the BRST-Cartan model.  We expect that, if $G$ is a source-connected and source-compact groupoid, then the basic cohomology arising from the $Q$-algebroid $[-1]TA$ is independent of the choice of connection and isomorphic to $H^\bullet (BG)$; however, we do not prove this.  We also expect that these ``algebroid characteristic classes'' should be closely related to those of Crainic and Fernandes \cite{crainic:vanest} \cite{cf} \cite{fernandes}.  We plan to explore these topics in future work.

The structure of the thesis is as follows.  In Chapter \ref{chapter:supergeom}, we provide an introduction to the theory of $\integers$-graded supermanifolds.  Because our main interest is in double complexes and cohomology theories, we prefer to work with $\integers$-gradings instead of the usual $\integers_2$-gradings.  Although the idea of $\integers$-graded supermanifolds is not new, the existing literature does not possess a general introduction to the theory, so we have provided one in order to have a logical foundation for the later chapters.  Superalgebroids and $Q$-algebroids are introduced in \S\ref{section:superalg}.

In Chapter \ref{chapter:supergpds}, we define supergroupoids and describe the Lie functor from the category of supergroupoids to the category of superalgebroids.

The heart of the thesis is Chapter \ref{chapter:double}, in which we describe the double complexes arising from $\mathcal{LA}$-groupoids, $Q$-groupoids, and $Q$-algebroids.  In this chapter we focus on the relatively familiar cases of $[-1]TG$ and $[-1]TA$, particularly in relation to the various models of equivariant cohomology.

In Chapter \ref{ch:examples}, we describe the following other examples of $Q$-groupoids and $Q$-algebroids:
\begin{itemize}
\item If $G$ is a Poisson-Lie group, then the $Q$-groupoid $[-1]T^*G \arrows [-1]\mathfrak{g}^*$ is a special case of the $Q$-groupoid arising from a Poisson groupoid.  The associated $Q$-algebroid is $[-1]T^*\mathfrak{g} = [-1]\mathfrak{g}^* \oplus \mathfrak{g}$, and the induced double complex is equal to the Chevalley-Eilenberg complex for the Drinfel'd double.  This homological point of view toward the double was studied by Lecomte, Roger \cite{lr} and Kosmann-Schwarzbach \cite{ks1}.  More generally, if $G$ is a Poisson groupoid, then the double complex of the $Q$-algebroid $[-1]T^* A$ may be identified with Roytenberg's  \cite{royt} ``commuting Hamiltonians'' description of the Lie bialgebroid $(A, A^*)$.
\item If a Lie algebroid $A$ is equipped with a compatible action of $TG$, where $G$ is a Lie group, then one can form an $\mathcal{LA}$-groupoid of the form 
\begin{equation}\la{A \times TG}{M \times G}{A}{M.}\end{equation}
We identify the double complex of the associated $Q$-algebroid with a complex introduced by Ginzburg \cite{ginzburg} in the context of equivariant Poisson cohomology.  Thus we may view the double complex of the $Q$-groupoid $[-1]A \times [-1]TG \arrows [-1]A$ as a global, as opposed to infinitesimal, model for ``equivariant algebroid cohomology''.
\item If $M$ is a Poisson manifold equipped with the Poisson action of a Poisson-Lie group $G$, then one can form an $\mathcal{LA}$-groupoid that in a sense integrates the ``matched pair'' of Lie algebroids introduced by Lu \cite{lu} in the study of Poisson homogeneous spaces.
\end{itemize}

In Chapter \ref{ch:vanest}, we describe an extension of the van Est map \cite{crainic:vanest} \cite{vanest} \cite{wx} that maps the double complex of a $Q$-groupoid to the double complex of its $Q$-algebroid.  The results of this chapter are from \cite{mw}.

Throughout the thesis, we will often omit the prefix ``super-''.  We will generally use calligraphic letters (e.g. $\mathcal{M}$, $\mathcal{G}$) to denote objects that are assumed to be in the category of supermanifolds.  For emphasis, we will occasionally use the term ``ordinary'' to mean ``not super''.  Unless otherwise stated, everything that is not super is assumed to be in the smooth category.

When writing formulas in coordinates, we consistently use the Einstein summation convention.  Since we usually deal with graded objects, signs involving gradings often appear in formulae, e.g. (\ref{eqn:modhom}); such a formula is, of course, only valid for homogeneous elements.
    
    \chapter{$\integers$-graded supergeometry}
    \label{chapter:supergeom}
Supermanifolds (with a $\integers_2$-grading) were introduced in the physics literature \cite{ss} to provide a formalism to describe supersymmetric field theories.  Informally, a $\integers_2$-graded supermanifold is a manifold with additional ``odd'' (or, in the physicists' terminology, \emph{fermionic}) coordinate functions that anticommute with each other.  The formal definition \cite{kostant} is sheaf-theoretic, defining a supermanifold as a sheaf of $\integers_2$-graded algebras that locally looks like the algebra of functions on a \emph{superspace}.  Some early mathematical texts on the subject were written by Berezin \cite{berezin}, Kostant \cite{kostant}, and Leites \cite{leites}; more recent treatments include \cite{bdm}, \cite{manin}, \cite{tuynman},  and \cite{varadarajan}.

Clearly, $\integers_2$-graded manifolds have been studied in much depth; in contrast, very little has been written about $\integers$-graded supermanifolds.  Kostant \cite{kostant} and Tuynman \cite{tuynman} do mention the possibility of using a $\integers$-grading as well as more general gradings, but the theory is developed only for the $\integers_2$ case.

Our choice to work with $\integers$-graded supermanifolds is somewhat inspired by the graded manifolds and $N$-manifolds that have appeared in the work of Kontsevich \cite{kont:deformation}, Roytenberg \cite{royt:graded}, Severa \cite{severa}, and Voronov \cite{voronov}, primarily in relation to Poisson geometry, Lie algebroids, and Courant algebroids.   Although there are certain similarities, our $\integers$-graded supermanifolds differ from theirs, all of which require an additional $\integers$-grading (the \emph{weight}) on the structure sheaf of a $\integers_2$-graded supermanifold.  Such structures include some peculiarities; for example, if $x$ is an even coordinate of weight $1$, then $e^x$ is a function that does not finitely decompose into functions of homogeneous weight.  It is concerns such as these that lead us to consider supermanifolds that are inherently $\integers$-graded.  We point out, however, that for Voronov's graded manifolds, weight is allowed to be independent of parity, making them somewhat more versatile than ours. 

It should be emphasized that, although $\integers_2$-graded supermanifolds and $\integers$-graded supermanifolds bear many formal similarities, they are not strictly comparable.  While Batchelor's theorem \cite{bat} \cite{berezin} \cite{gawedzki} asserts that any $\integers_2$-graded supermanifold admits a compatible $\integers$-grading, the $\integers$-grading is not canonical.  Conversely, it is not possible in general to obtain a $\integers_2$-graded supermanifold by applying a ``mod $2$'' functor to a $\integers$-graded supermanifold, since $\integers$-graded supermanifolds may have coordinates that are of even, nonzero degrees, in which functions are polynomial, as opposed to smooth.  The target of the ``mod $2$'' functor is actually the category of Konechny and Schwarz's \emph{partially formal supermanifolds} \cite{kon-sch}.

Since, to our knowledge, the foundations of $\integers$-graded supermanifold theory have not previously been published, there is a sense in which all of the material in this chapter is new (though it perhaps has been ``in the air'').  However, the exposition, particularly in \S\ref{section:superm}, \S\ref{section:gla}, and \S\ref{section:calc}, leans heavily on the parallel exposition in \cite{bdm}, and most of the results are obvious extrapolations of the analogous results in the $\integers_2$ case.  The more interesting results, then, are those that depend on the fact that $-1 \neq 1$, such as the observation that tangent bundles and cotangent bundles have different dimension (Example \ref{ex:cotangent}).

As a final note, this chapter is not meant to be an exhaustive treatise on $\integers$-graded supermanifold theory, but rather the minimal set of results that are necessary for the theory of supergroupoids which follows.  Many elements that one might desire from a full supermanifold theory, such as a $\integers$-graded Berezinian and integration theory or a result analogous to Batchelor's theorem, are not addressed here.  Also, the ``functor of points'' approach to supermanifolds is not utilized here in favor of more explicit constructions, although some material, for example the Cartan calculus \cite{ks:gorms}, arises more naturally in that setting.

\section{Supermanifolds}\label{section:superm}

In this section, $\integers$-graded supermanifolds are defined, and some basic examples are given.  Supermanifolds are locally modelled after superdomains or superspaces.  Whereas in the $\integers_2$ case the dimension of a supermanifold is described by a pair of nonnegative integers $(p,q)$, the dimension of a $\integers$-graded supermanifold is described by a sequence of nonnegative integers $\{p_i\}$.  All graded objects are assumed to be $\integers$-graded unless otherwise indicated.

\subsection{Superdomains}

Let $\{ p_i \}_{i \in \integers}$ be a nonnegative integer-valued sequence.  Denote by $\calO^{ \{p_i\} }$ the sheaf of graded, graded-commutative algebras on $\reals^{p_0}$ defined by
\begin{equation}
\calO^{ \{p_i\} } (U) = C^\infty(U)\left[\bigcup_{i \neq 0} \{\xi_i^1, \dots, \xi_i^{p_i}\} \right]
\end{equation}
for any open set $U \subseteq \reals^{p_0}$, where $\xi_i^k$ is of degree $-i$.  In particular, $\xi_i^k \xi_j^\ell = (-1)^{ij}\xi_j^\ell\xi_i^k$.  

\begin{dfn} The \emph{coordinate superspace} $\reals^{ \{p_i\} }$ is the pair $\left( \reals^{p_0} , \calO^{ \{p_i\} }\right)$.
\end{dfn}

\begin{rmk} The basic premise of supergeometry is that we treat $\reals^{ \{p_i\} }$ as if it were a space whose sheaf of ``smooth functions'' is $\calO^{ \{p_i\} }$.  Following this idea, we write $C^\infty\left(\reals^{ \{p_i\} }\right) \defequal \calO^{ \{p_i\} }$.
\end{rmk}

There is a natural surjection of sheaves $\ev: C^\infty\left(\reals^{ \{p_i\} }\right) \to C^\infty\left(\reals^{p_0}\right)$, called the \emph{evaluation map}, where the kernel is the ideal generated by all elements of nonzero degree.

\begin{dfn}
A \emph{superdomain} $\mathcal{U}$ of dimension $\{ p_i \}$ is a pair $\left(U, C^\infty(\mathcal{U})\right)$, where $U$ is an open subset of $\reals^{p_0}$ and $C^\infty(\mathcal{U}) \defequal \calO^{ \{p_i\} }|_U$.
A morphism of superdomains $\mu: \mathcal{U} \to \mathcal{V}$ consists of a smooth map $\mu_0: U \to V$ and a morphism of sheaves of graded algebras $\mu^*: C^\infty(\mathcal{V}) \to C^\infty(\mathcal{U})$ over $\mu_0$, such that $\ev \circ \mu^* = \mu_0^* \circ \ev$.
\end{dfn}

\begin{rmk} To describe a superdomain, it is only necessary to specify the dimension $\{ p_i \}$ and an open set $U \in \reals^{p_0}$.
\end{rmk}

\subsection{Supermanifolds}

\begin{dfn}\label{dfn:supermfld}
A \emph{supermanifold} $\mathcal{M}$ of dimension $\{p_i\}$ is a pair $(M, C^\infty(\mathcal{M}))$, where $M$ (the \emph{support}) is a topological space and $C^\infty(\mathcal{M})$ is a sheaf on $M$ of graded algebras (the \emph{function sheaf}) that is locally isomorphic to a superdomain of dimension $\{p_i\}$.  A morphism of supermanifolds $\mu: \mathcal{M} \to \mathcal{N}$ consists of a smooth map $\mu_0: M \to N$ and a morphism of sheaves of graded algebras $\mu^*: C^\infty(\mathcal{N}) \to C^\infty(\mathcal{M})$ over $\mu_0$, such that $\ev \circ \mu^* = \mu_0^* \circ \ev$.
\end{dfn}

\begin{rmk}
It is immediate from Definition \ref{dfn:supermfld} that if $\mathcal{M} = (M, C^\infty(\mathcal{M}))$ is a dimension $\{p_i\}$ supermanifold, then the topological space $M$ automatically has the structure of a $p_0$-dimensional manifold.  The evaluation map describes an embedding of $M$ into $\mathcal{M}$ (see Definition \ref{dfn:embed}).
\end{rmk}

\begin{ex} \label{ex:vbshift}Let $\pi: E \to M$ be a vector bundle.  The supermanifold $[j]E$ has support $M$ and function sheaf $\bigwedge \Gamma(E^*) = \Gamma\left(\bigwedge E^*\right)$ if $j$ is odd or $\bigS \Gamma (E^*) = \Gamma\left(\bigS E^*\right)$ if $j \neq 0$ is even; in either case the grading is the natural $\integers$-grading multiplied by $-j$.  The dimension of $[j]E$ is $\{p_i\}$, where $p_0 = \dim M$, $p_j = \rank E$, and $p_i = 0$ for all other values of $i$.  This example is discussed in more detail in \S\ref{sub:vb}.
\end{ex}

\begin{rmk} The notation in Example \ref{ex:vbshift} differs from that in the existing literature, e.g. \cite{kont:deformation}, \cite{royt:graded}, \cite{severa}, where the notation $E[-j]$ has been used for the supermanifold which we denote by $[j]E$.  There are two separate distinctions at work here.  The first is that, following a suggestion of Weinstein, we have placed the ``degree shift'' operator on the left in order to emphasize the fact that it is a functor (see Definition \ref{dfn:shift} below)\footnote{This convention also agrees with the fact that $[j]E$ is the $\integers$-graded analogue of the $\integers_2$-graded supermanifold $\Pi E$.}.  The second distinction is that, in the spirit of supergeometry, we have interpreted the operation $[j]$ to be a geometric, as opposed to an algebraic, operation.  In other words, the fibres of $[j]E$ are of degree $j$, whereas in the previous literature $E[j]$ has been characterized by the property that the linear functions are of degree $j$.  Because the degree of a vector space is opposite (in sign) to the degree of its linear functions (see Remarks \ref{rmk:dual} and \ref{rmk:dual2} below), our degree shift operator differs by a sign from the degree shift operator in the existing literature.
\end{rmk}

\begin{dfn} \label{dfn:embed} A morphism $\mu: \mathcal{M} \to \mathcal{N}$ of supermanifolds is an \emph{embedding} if $\mu_0$ is an embedding and $\mu^*$ is a surjection of sheaves.
\end{dfn}

\begin{rmk} \label{rmk:subideal} If $\mu: \mathcal{M} \to \mathcal{N}$ is an embedding, then $\ker \mu^* \subseteq C^\infty(\mathcal{N})$ is a sheaf of ideals that completely determines the supermanifold structure of $\mathcal{M}$ and the embedding.  In this case, we denote the sheaf of ideals as $I(\mathcal{M}) \defequal \ker \mu^*$.  Not every sheaf of ideals corresponds to an embedded submanifold.
\end{rmk}

\begin{dfn} A morphism $\mu: \mathcal{M} \to \mathcal{N}$ of supermanifolds is a \emph{surjection} if $\mu_0$ is an surjection and $\mu^*$ is an injection of sheaves.
\end{dfn}

\begin{dfn} Let $\mathcal{M} = (M, C^\infty(\mathcal{M}))$ be a supermanifold.  An \emph{open submanifold} of $\mathcal{M}$ is of the form $\mathcal{N} = (N, C^\infty(\mathcal{M})|_N)$, where $N$ is an open submanifold of $M$.
\end{dfn}

\begin{rmk} Because of the local structure of supermanifolds, there exists on any supermanifold $\mathcal{M}$ of dimension $\{p_i\}$ an open cover $\{\mathcal{U}_\alpha\}$ such that each $\mathcal{U}_\alpha$ is isomorphic as a supermanifold to a superdomain $\mathcal{V}_\alpha \subseteq \reals^{\{p_i\}}$.  If we choose for every $\alpha$ an isomorphism $\psi_\alpha: \mathcal{U}_\alpha \to \mathcal{V}_\alpha$, then the supermanifold structure is completely determined by the transition maps $\varphi_{\alpha \beta} = \psi_\beta \circ \psi_\alpha^{-1} : \psi_\alpha(\mathcal{U}_\alpha \cap \mathcal{U}_\beta) \to \psi_\beta(\mathcal{U}_\alpha \cap \mathcal{U}_\beta)$.  The transition maps form a \v{C}ech $1$-cocycle with values in the automorphism sheaf of $\reals^{\{p_i\}}$.
\end{rmk}

\section{Vector bundles}\label{sub:vb}

This section deals with vector bundles in the category of supermanifolds.  There are two typical approaches to vector bundles: the topological approach is to consider a vector bundle as being built out of local trivializations and glueing maps, whereas the algebraic approach is to consider a vector bundle as a locally free module.  The main goal of this section is to describe the two points of view and prove an analogue of the Serre-Swan theorem, which states that the two points of view are equivalent.

After a brief introduction to graded linear algebra, vector bundles are defined in \S\ref{section:bundles} from the topological point of view.  The $\Gamma$ functor from vector bundles to their modules of sections is defined, and the equivalence of the two points of view is shown in Theorem \ref{thm:serreswan}.  Finally, in \S\ref{section:dvb}, double vector bundles are discussed.  Double vector bundles play a significant role in later sections, e.g. \S\ref{sec:oneta}.

\subsection{Graded linear algebra}\label{section:gla}

Let $\{V_i\}$, $i \in \integers$, be a collection of vector spaces.  Then $V = \bigoplus V_i$ is a \emph{graded vector space}.  We use the notation $|v|$ to denote the grading of a homogeneous element; in other words, if $v \in V_i$, then $|v| = i$.  

The grading on the tensor product of two graded vector spaces is as follows:
\begin{equation}
(V \otimes W)_i = \bigoplus_{j+k = i}(V_j \otimes W_k).
\end{equation}

\begin{rmk} The class of quasi-finite-dimensional graded vector spaces (those such that $\dim V_i < \infty$ for each $i$) is not closed under the tensor product.  For example, if $\dim V_i = 1$ for all $i$, then $(V \otimes V)_i$ is infinite-dimensional.  For this reason, one should either restrict to finite-dimensional graded vector spaces (where $\sum \dim V_i < \infty$) or allow infinite-dimensional graded vector spaces (with no restrictions on $\dim V_i$)\footnote{Another possibility is to consider graded vector spaces whose grading is bounded below (or above).  The class of all such vector spaces is closed under the tensor product.}.  In what follows, we will assume that all graded vector spaces are finite-dimensional.
\end{rmk}

The key aspect of linear superalgebra is the sign rule, which is the convention that the ``canonical'' isomorphism between $V \otimes W$ and $W \otimes V$ is the one that sends $v \otimes w$ to $(-1)^{|v||w|}w \otimes v$ for homogeneous $v \in V, w \in W$.  In particular, the graded symmetric algebra $\bigS(V)$ is such that elements of odd degree anticommute with each other.

\begin{dfn} Let $V$ and $W$ be graded vector spaces.  Then $\Hom(V,W)$ is a graded vector space, where the degree $i$ subspace $\Hom_i(V,W)$ is equal to $\bigoplus_{j \in \integers} \Hom(V_j, W_{i+j})$.
\end{dfn}

\begin{dfn}The \emph{dual} $V^*$ of a graded vector space is equal to $\Hom(V, \reals)$.
\end{dfn}

\begin{rmk} \label{rmk:dual}As graded vector spaces, $V$ and $V^*$ are in general non-isomorphic, since they have ``opposite'' gradings.  For example, if $V$ is purely of degree $j$, then $V^*$ is purely of degree $-j$.  
\end{rmk}

\begin{rmk} \label{rmk:dual2}
The correspondence between vector spaces and linear manifolds does not immediately extend to the graded case.  However, there is a functorial relationship between graded vector spaces and linear supermanifolds that gives an equivalence of categories.  We will now briefly describe this relationship.

Suppose that $V$ is a graded vector space of dimension $\{k_i\}$.  Then $\reals^{\{k_i\}}$ may be (noncanonically) identified with $V$ in the sense that a homogeneous basis of $V^*$ can be identified with coordinates on $\reals^{\{k_i\}}$.  The group of invertible degree $0$ elements of $\mathrm{End}(V)$ can then be identified with the group of linear coordinate transformations of $\reals^{\{k_i\}}$.  One can use this data to canonically associate to $V$ a ``linear'' supermanifold $\mathcal{V}$, as follows.  

The supermanifold $\mathcal{V}$ has support $V_0$ and has a collection of isomorphisms $\phi_a: \mathcal{V} \to \reals^{\{k_i\}}$, where $a$ is any basis of $V$.  Given any two bases $a$ and $b$, the transformation map $\phi_b \circ \phi_a^{-1}: \reals^{\{k_i\}} \to \reals^{\{k_i\}}$ is the coordinate transformation induced by the change of basis from $b$ to $a$.  Since all of the coordinate transformations obtained in this way are linear, $\mathcal{V}$ possesses a well-defined subspace $C^\infty_{lin}(\mathcal{V}) \subseteq C^\infty(\mathcal{V})$ of linear functions.  The correspondence between $V$ and $\mathcal{V}$ is an equivalence of categories between the category of graded vector spaces and the category of linear supermanifolds.
\end{rmk}

\begin{dfn} \label{dfn:shift} Let $j \in \integers$.  The \emph{degree-shifting functor} $[j]$ acts on the category of graded vector spaces as follows:
\begin{equation}
[j]V = [j]\reals \otimes V,
\end{equation}
where $[j]\reals$ is the one-dimensional vector space of degree $j$.
\end{dfn}

\begin{rmk} The degree-shifting functor may be described more explicitly by the property $([j]V)_i = V_{i-j}$.  However, Definition \ref{dfn:shift} is useful for doing computations in light of the sign rule.  In particular, the graded vector space $V$ may be written as $\bigoplus [i]V_{(i)}$, where each $V_{(i)}$ is a degree $0$ vector space of the same dimension as $V_i$.  A general element of $V$ is then a finite sum $v = \epsilon^i v_i$, $v_i \in V_{(i)}$, where $\epsilon^i$ is the standard basis vector of $[i]\reals$.  This point of view allows linear algebraic formulae to be naturally extended to graded linear algebra via the so-called ``even rules''\footnote{For example, a Lie superbracket should satisfy the skew-symmetry identity $[\epsilon^i v_i, \epsilon^j w_j] = \epsilon^i\epsilon^j[ v_i,  w_j] = -\epsilon^i\epsilon^j[ w_j,  v_i] = (-1)^{1+ij} \epsilon^j\epsilon^i[ w_j,  v_i] = (-1)^{1+ij}[\epsilon^j w_j, \epsilon^i v_i]$.  Thus the rule for Lie superalgebras has been derived from the rule for (degree $0$) Lie algebras.} \cite{bdm}.
\end{rmk}

The following identities are immediate from the above definitions.

\begin{prop} \label{prop:jarule}Let $V$ and $W$ be graded vector spaces.  Then
\begin{enumerate}
\item $[j][k]V = [j+k]V$
\item $[j](V \oplus W) = [j]V \oplus [j]W$,
\item $[j](V \otimes W) = [j]V \otimes W = V \otimes [j]W$,
\item $([j]V)^* = [-j](V^*)$,
\item $\Hom_i (V,W) = \Hom_0 (V, [-i]W) = \Hom_0 ([i]V, W)$.
\end{enumerate}
\end{prop}

\begin{rmk} \label{rmk:jarule} One should be careful when using the identities of Proposition \ref{prop:jarule}, since it is possible to draw diagrams that commute only up to a sign.  For example, there is a natural isomorphism between $[j][k]V$ and $[k][j]V$ that sends $\epsilon^j \epsilon^k v$ to $(-1)^{jk}\epsilon^k \epsilon^j v$.  On the other hand, the identifications $[j][k]V = [j+k]V = [k][j]V$ describe an isomorphism that sends $\epsilon^j \epsilon^k v$ to $\epsilon^k \epsilon^j v$ .  Thus, when identifying $[j][k]V$ and $[k][j]V$, it is necessary to specify which isomorphism is being used.
\end{rmk}

\subsection{Bundles and sections}\label{section:bundles}

Let $\mathcal{M}$ be a supermanifold with support $M$, and let $U$ be an open subset of $M$.  Denote by $\mathcal{M}|_U$ the supermanifold with support $U$ whose function sheaf is the restriction of $C^\infty(\mathcal{M})$ to $U$.

\begin{dfn} \label{dfn:bundle-top}A \emph{vector bundle} of rank $\{k_i\}$ over $\mathcal{M}$ is a supermanifold $\mathcal{E}$ and a surjection $\pi: \mathcal{E} \to \mathcal{M}$ equipped with an atlas of local trivializations $\mathcal{E}|_{\pi_0^{-1}(U)} \cong \mathcal{M}|_U \times \reals^{\{k_i\}}$ such that the transition function between any two local trivializations is linear in the fibre coordinates.
\end{dfn}

\begin{rmk}Any operation on graded vector spaces has an associated operation on vector bundles.  In particular, the degree-shifting functor $[j]$ has the effect of lowering the degrees of the fibre coordinates by $j$, producing the vector bundle $[j]\pi: [j]\mathcal{E} \to \mathcal{M}$ (see Remark \ref{rmk:dual}).
\end{rmk}

\begin{ex} Let $E \to M$ be a graded vector bundle, i.e. a sequence of vector bundles $E_i$ ($i \neq 0$) such that $E = \bigoplus_i E_i$.  The vector bundle $\bigoplus_i [i] E_i \to M$ is denoted by $[\integers] E$.  The rank of $[\integers] E$ is $\{k_i\}$, where $k_i = \rank E_i$.
\end{ex}

\begin{dfn} Let $\mathcal{E}$ and $\mathcal{E'}$ be vector bundles over $\mathcal{M}$.  A degree $j$ \emph{vector bundle morphism} from $\mathcal{E}$ to $\mathcal{E'}$ over $\mathcal{M}$ is a map $\tau: \mathcal{E} \to [-j]\mathcal{E}'$ which is linear in the fibre coordinates and such that $[-j]\pi' \circ \tau = \pi$. \end{dfn}

\begin{dfn}Let $\pi: \mathcal{E} \to \mathcal{M}$ be a vector bundle.  A \emph{homogeneous section} of degree $j$ is a map $X: \mathcal{M} \to [-j]\mathcal{E}$ such that $[-j]\pi \circ X = id_\mathcal{M}$.  The space of degree $j$ sections is denoted $\Gamma_j (\mathcal{E})$.  The \emph{space of sections} is $\Gamma(\mathcal{E}) \defequal \bigoplus_j \Gamma_j (\mathcal{E})$
\end{dfn}

\begin{prop} \label{prop:shiftsecs}$\Gamma([j]\mathcal{E}) \cong [j]\Gamma(\mathcal{E})$.\end{prop}
\begin{proof} Since $[-k]\mathcal{E} = [- (j+k) + j] \mathcal{E} = [-(j+k)][j]\mathcal{E}$, it is immediate from the definition that degree $k$ sections of $\mathcal{E}$ are in one-to-one correspondence\footnote{This bijection is, however, only canonical ``up to sign'' since one must choose a convention for making the identification of $[-j]\mathcal{E}$ and $[-(j+k)][k]\mathcal{E}$ (see Remark \ref{rmk:jarule}).} with degree $(j+k)$ sections of $[j]\mathcal{E}$.
\end{proof}

A degree $j$ section $X$ of $\mathcal{E}$ is completely determined by its action (via pullback) on $C^\infty_{lin}([-j]\mathcal{E})$, the graded space of functions on $[-j]\mathcal{E}$ that are linear in the fibre coordinates. 
In a local trivialization, let $\{\xi^a\}$ be a set of fibre coordinates on $\mathcal{E}$, and let $\{\hat{\xi}^a\}$ be the corresponding set of fibre coordinates on $[-j]\mathcal{E}$.  Then a set of functions $\{f^a\} \in C^\infty(\mathcal{M})$ such that $|f^a| = |\hat{\xi}^a| = |\xi^a| + j$ defines locally a section $X \in \Gamma_j (\mathcal{E})$ with the property $X^*(\hat{\xi}^a) = f^a$.  In particular, the \emph{frame of sections} $\{X_a\}$ dual to the coordinates $\{\xi^a\}$ is such that $|X_a| = -|\xi^a|$ and $X_a^*(\overline{\xi}^b) = \delta_a^b$, where the $\overline{\xi}^b$ are the appropriately grading-shifted fibre coordinates.

It is clear from the local picture that $C^\infty_{lin}(\mathcal{E})$ may be identified with $\Gamma(\mathcal{E}^*)$.  Our next objective will be to describe a pairing of $\Gamma(\mathcal{E})$ and $C^\infty_{lin}(\mathcal{E})$ that gives an isomorphism of $C^\infty_{lin}(\mathcal{E})$ and the module dual $\Hom_{C^\infty(\mathcal{M})} (\Gamma(\mathcal{E}), C^\infty(\mathcal{M}))$.

Locally, any linear function $\alpha \in C^\infty_{lin}(\mathcal{E})$ can be written in the form $\alpha = \xi^a g_a$, where $g_a \in C^\infty(\mathcal{M})$.  The map $\xi^a g_a \mapsto \hat{\xi}^a g_a$ gives a right $C^\infty(\mathcal{M})$-module isomorphism between $C^\infty_{lin}(\mathcal{E})$ and $C^\infty_{lin}([-j]\mathcal{E})$.  Similarly, the map $g_a \xi^a \mapsto g_a \hat{\xi}^a$ gives a left $C^\infty(\mathcal{M})$-module isomorphism.  Because both maps are of degree $j$, it is not possible to respect the left and right module structures simultaneously.

If $X \in \Gamma_j (\mathcal{E})$, then we can use the right module isomorphism\footnote{If we were to use the left module isomorphism instead, then the discussion which follows could be carried out for the opposite pairing $\langle \cdot , \cdot \rangle:C^\infty_{lin}(\mathcal{E})  \otimes \Gamma(\mathcal{E}) \to C^\infty(\mathcal{M})$.}
 of $C^\infty_{lin}(\mathcal{E})$ and $C^\infty_{lin}([-j]\mathcal{E})$ to obtain a degree $j$ map $X_0^*: C^\infty_{lin}(\mathcal{E}) \to C^\infty(\mathcal{M})$.  It is straightforward to check that $X_0^*$ is a right $C^\infty(\mathcal{M})$-module homomorphism, and that it satisfies
\begin{equation} \label{eqn:modhom}X_0^* (f\alpha) = (-1)^{|f||X|} f X_0^*(\alpha) \end{equation}
for any $f \in C^\infty(\mathcal{M})$ and $\alpha \in C^\infty_{lin}(\mathcal{E})$.

\begin{dfn} The pairing $\langle \cdot , \cdot \rangle: \Gamma(\mathcal{E}) \otimes C^\infty_{lin}(\mathcal{E}) \to C^\infty(\mathcal{M})$ is defined by 
\begin{equation}\label{eqn:pairing}\langle X, \alpha \rangle \defequal X_0^*(\alpha).\end{equation}
\end{dfn}

\begin{prop} \label{prop:pairing} The pairing $\langle \cdot , \cdot \rangle$ is a right $C^\infty(\mathcal{M})$-module homomorphism and is nondegenerate.
\end{prop}

\begin{proof}
It may immediately be verified from (\ref{eqn:pairing}) that $\langle X, \alpha f\rangle = \langle X, \alpha \rangle f$, showing that the pairing is a right module homomorphism.

Let $\{\xi^a\}$ be fibre coordinates in a local trivialization.  The dual frame of sections $\{X_a\}$ clearly satisfies the property
\begin{equation}\langle X_a, \xi^b \rangle = \delta_a^b,\end{equation}
thus giving the nondegeneracy.
\end{proof}

\begin{rmk}The pairing (\ref{eqn:pairing}) gives $\Gamma(\mathcal{E})$ a $C^\infty(\mathcal{M})$-module structure defined by the equation
\begin{equation}\label{eqn:module}\langle Xf, \alpha \rangle = \langle X, f \alpha \rangle. \end{equation}
With this module structure, the pairing is a bimodule homomorphism.  It follows from the nondegeneracy that the pairing gives a $C^\infty(\mathcal{M})$-module isomorphism between $C^\infty_{lin}(\mathcal{E})$ and $\Hom_{C^\infty(\mathcal{M})} (\Gamma(\mathcal{E}), C^\infty(\mathcal{M}))$.
\end{rmk}

\begin{rmk} $\Gamma(\mathcal{E})$ can also be be considered as a sheaf on the support $M$, associating to any open set $U \subseteq M$ the space of sections of the bundle $\mathcal{E}|_{\pi_0^{-1}(U)} \to \mathcal{M}|U$.  Then $\Gamma(\mathcal{E})$ is a $C^\infty(\mathcal{M})$-module in the sheaf-theoretic sense and is locally free.  If $\mathcal{E}$ is of rank $\{k_i\}$, then $\Gamma(\mathcal{E})$ is locally isomorphic to $C^\infty(\mathcal{M}) \otimes \reals^{\{k_i\}}$.
\end{rmk}

\begin{prop} Let $\tau: \mathcal{E} \to [-j]\mathcal{E}'$ be a degree $j$ bundle morphism over $\mathcal{M}$.  Then there is an induced degree $j$ morphism of $C^\infty(\mathcal{M})$-modules $\Gamma(\tau): \Gamma(\mathcal{E}) \to \Gamma(\mathcal{E}')$.\end{prop}
\begin{proof}
The map $\tau$ gives rise to bundle morphisms $[-k]\tau: [-k]\mathcal{E} \to [-(k+j)]\mathcal{E}'$ for all $k$.  Given a degree $k$ section $X \in \Gamma(\mathcal{E})$, composition with $[-k]\tau$ gives a section $\Gamma(\tau)(X)\in \Gamma_{k+j}(\mathcal{E}')$.  By considering the module structure in terms of pullbacks, it is clear that $\Gamma(\tau)$ respects the module structure.
\end{proof}

Thus $\Gamma$ is a functor from the category of vector bundles over $\mathcal{M}$ to the category of locally free $C^\infty(\mathcal{M})$-modules.  The following is a version of the Serre-Swan theorem for vector bundles in the category of supermanifolds.

\begin{thm}\label{thm:serreswan}
The functor $\Gamma$ is an equivalence of categories.
\end{thm}

\begin{proof}
First we will show that $\Gamma$ is a full functor.  Let $\{X_a\}$ and $\{Y_b\}$ be frames of sections on bundles $\mathcal{E}$ and $\mathcal{E'}$, respectively, and let $\{\xi^a\}$ and $\{\eta^b\}$ be the respective sets of dual coordinates.  A morphism of $C^\infty(\mathcal{M})$-modules $\varphi: \Gamma(\mathcal{E}) \to \Gamma(\mathcal{E}')$ is locally of the form $\varphi: X_a \mapsto \varphi_a^b Y_b$, where $\varphi_a^b \in C^\infty(\mathcal{U})$.  There is a corresponding bundle map $\tau$, defined by the requirement that
\begin{equation}
\langle X_a, \tau^* (\eta^b) \rangle = \langle \varphi (X_a), \eta^b\rangle,
\end{equation}
which is satisfied by setting $\tau^* (\eta^b) = \xi^c \varphi^b_c$.  Thus any morphism of $C^\infty(\mathcal{M})$-modules from $\Gamma(\mathcal{E})$ to $\Gamma(\mathcal{E}')$ arises from a bundle map $\mathcal{E} \to \mathcal{E}'$.

It remains to show that $\Gamma$ is faithful and surjective.  There is an obvious bijection between the class of free $C^\infty(\mathcal{M})$-modules and the class of trivial vector bundles over $\mathcal{M}$.  Since we can identify the local structures as well as the tranformation maps, we can globally identify locally free modules and vector bundles.
\end{proof}

\begin{prop}\label{prop:dual}
\begin{enumerate}
\item $\Gamma(\mathcal{E} \otimes \mathcal{E}')$ is naturally isomorphic to $\Gamma(\mathcal{E}) \otimes_{C^\infty(\mathcal{M})} \Gamma(\mathcal{E}')$.
\item $\Gamma(\mathcal{E}^*)$ is naturally isomorphic to $\Hom_{C^\infty(\mathcal{M})}\left(\Gamma(\mathcal{E}), C^\infty(\mathcal{M})\right)$.
\end{enumerate}
\end{prop}
\begin{proof}
Let $\{X_i\}$ and $\{Y_j\}$ be frames of sections on $\mathcal{E}$ and $\mathcal{E'}$, respectively.  Using the linear algebraic fact that, for vector spaces $V$ and $W$, $(V \otimes W)^*$ is naturally isomorphic to $V^* \otimes W^*$, dual fibre coordinates $\{\xi^i\}$ on $\mathcal{E}$ and $\{\eta^i\}$ on $\mathcal{E'}$ give rise to fibre coordinates $\{\xi^i \otimes \eta^j\}$ on $\mathcal{E} \otimes \mathcal{E'}$.  Identifying $\{X_i \otimes Y_j\}$ with the dual frame of sections gives the  first identity.

The second statement follows from Proposition \ref{prop:pairing} (also see Remark \ref{rmk:dual}).
\end{proof}

\begin{ex}[Tangent bundle]  \label{ex:tangent}
For a superdomain $\mathcal{U}$ of dimension $\{p_i\}$, the tangent bundle $T\mathcal{U}$ is the trivial bundle $\mathcal{U}\times\reals^{\{p_i\}}$.  Let $\{x^i\}$ be coordinates on $\mathcal{U}$, and denote the fibre coordinates by $\{\dot{x}^i\}$ (so that $|\dot{x}^i| = |x^i|$).  Naturally associated to any morphism $\mu: \mathcal{U} \to \mathcal{V}$ is a bundle map $T\mu$, defined as follows.  Let $\{y^i, \dot{y}^i\}$ be coordinates on $T\mathcal{V}$.  Then
\begin{equation} \label{eqn:tangent}
(T\mu)^* (\dot{y}^i) = \dot{x}^j \pdiff{}{x^j} [\mu^*(y^i)] .
\end{equation}
These properties completely define the tangent bundle $T\mathcal{M}$ of any supermanifold $\mathcal{M}$.  The operation is functorial.  
\end{ex}

\begin{ex}[Cotangent bundle] \label{ex:cotangent}
The cotangent bundle $T^*\mathcal{M}$ is dual to the tangent bundle.  It is worth emphasizing that, if $\mathcal{M}$ is of dimension $\{p_i\}$ then $T\mathcal{M}$ is of dimension $\{2p_i\}$, whereas $T^*\mathcal{M}$ is of dimension $\{p_i + p_{-i}\}$.
This is a distinction that does not appear in the study of $\integers_2$-graded supermanifolds.
\end{ex}

\subsection{Double vector bundles}\label{section:dvb}
Double vector bundles were introduced by Pradines \cite{pradines} (see also \cite{mac:dblie1}).  In this section, we collect some results on double vector bundles, particularly with respect to the degree-shifting functor, that will be valuable in later sections.

\begin{dfn}A \emph{double vector bundle} $(\mathcal{D},\mathcal{A};\mathcal{B},\mathcal{M})$ is a vector bundle in the category of vector bundles, namely a square
\begin{equation}
\xymatrix{\mathcal{D} \ar[r] \ar[d] & \mathcal{A} \ar[d] \\ \mathcal{B} \ar[r] & \mathcal{M}},
\end{equation}
where $\mathcal{D}$ is a vector bundle over $\mathcal{A}$, $\mathcal{B}$ is a vector bundle over $\mathcal{M}$, $\mathcal{D}$ is a vector bundle over $\mathcal{B}$, and the structure maps (the projection, addition, and scalar multiplication maps) for $\mathcal{D} \to \mathcal{B}$ are bundle morphisms.
\end{dfn}

\begin{rmk}
The concept of double vector bundles is symmetric in the sense that 
\begin{equation}
\xymatrix{\mathcal{D} \ar[r] \ar[d] & \mathcal{A} \ar[d] \\ \mathcal{B} \ar[r] & \mathcal{M}}
\end{equation}
is a double vector bundle whenever 
\begin{equation}
\xymatrix{\mathcal{D} \ar[r] \ar[d] & \mathcal{B} \ar[d] \\ \mathcal{A} \ar[r] & \mathcal{M}}
\end{equation}
is a double vector bundle.  As a result, the degree-shifting functor $[j]$ may be applied to a double vector bundle in two different ways.  For example, applying $[j]$ to the rows yields 
\begin{equation}
\xymatrix{[j]_{\mathcal{A}}\mathcal{D} \ar[r] \ar[d] & \mathcal{A} \ar[d] \\ [j] \mathcal{B} \ar[r] & \mathcal{M}},
\end{equation}
where the subscript on the functor $[j]$ indicates to which vector bundle structure $[j]$ is being applied.
\end{rmk}

\begin{prop} \label{prop:dvb1}The horizontal and vertical degree-shifting functors commute, in the sense that $[k]_{[j]\mathcal{B}} [j]_\mathcal{A} \mathcal{D} = [j]_{[k]\mathcal{A}} [k]_\mathcal{B} \mathcal{D}$.
\end{prop}
\begin{proof} Locally, there exist coordinates on $\mathcal{D}$ of the form $\{x^i, a^i, b^i,c^i\}$, where $\{x^i\}$ is a set of coordinates on $\mathcal{M}$, $\{a^i\}$ and $\{b^i\}$ are fibre coordinates on $\mathcal{A}$ and $\mathcal{B}$, respectively, and $\{c^i\}$ are coordinates on the \emph{intersections} of the fibres of $\mathcal{D}$ over $\mathcal{A}$ and $\mathcal{B}$ (called the \emph{core} by Mackenzie \cite{mac:vdbl}).  Both $[k]_{[j]\mathcal{B}} [j]_\mathcal{A}$ and $[j]_{[k]\mathcal{A}} [k]_\mathcal{B}$ shift the grading of $\{a^i\}$, $\{b^i\}$, and $\{c^i\}$ by $-k$, $-j$, and $-(j+k)$, respectively.
\end{proof}

\begin{rmk} In terms of the coordinates $\{x^i, a^i, b^i, c^i\}$ on $\mathcal{D}$ as above, and introducing auxilliary coefficients $\epsilon^j$ and $\epsilon^k$ to represent the degree-shifts, there is another natural isomorphism $[k]_{[j]\mathcal{B}} [j]_\mathcal{A} \mathcal{D} \to [j]_{[k]\mathcal{A}} [k]_\mathcal{B} \mathcal{D}$ is as follows:
\begin{equation}
(x^i, \epsilon^k a^i, \epsilon^j b^i, \epsilon^k \epsilon^j c^i) \mapsto (x^i, \epsilon^k a^i, \epsilon^j b^i, (-1)^{jk} \epsilon^j \epsilon^k c^i).
\end{equation}
In particular, if $j$ and $k$ are both odd, then a minus sign may be introduced when identifying the coordinates on the core of the double vector bundle (cf. Remark \ref{rmk:jarule}).
\end{rmk}

Let $(\mathcal{D},\mathcal{A};\mathcal{B},\mathcal{M})$ be a double vector bundle.  There are two special types of sections of $\mathcal{D}$ as a vector bundle over $\mathcal{A}$.

\begin{dfn} A degree $k$ section $X \in \Gamma(\mathcal{D},\mathcal{A})$ is \emph{linear} (with respect to the vertical bundle structures) if $X: \mathcal{A} \to [-k]_\mathcal{A} \mathcal{D}$ is a degree $0$ bundle morphism
\begin{equation}
\xymatrix{[-k]_\mathcal{A} \mathcal{D} \ar[r] \ar[d] & \mathcal{A} \ar@/_/[l]_-{X} \ar[d]\\
[-k]\mathcal{B} \ar[r] & \mathcal{M} \ar@/_/[l]_-{X_0}.}
\end{equation}
\end{dfn}

\begin{dfn} A degree $k$ section $X \in \Gamma(\mathcal{D},\mathcal{A})$ is \emph{vertical} (with respect to the vertical bundle structures) if, for all $\alpha \in C^\infty([-k]_\mathcal{A} \mathcal{D})$ that are linear with respect to both vector bundle structures, $X^* \alpha \in C^\infty(\mathcal{M})$.
\end{dfn}

\begin{prop} \label{prop:vertandlin}Let $\Gamma_{lin}(\mathcal{D}, \mathcal{A})$ and $\Gamma_V (\mathcal{D}, \mathcal{A})$ be the spaces of linear and vertical sections, respectively.  Then \begin{equation*}\Gamma(\mathcal{D}, \mathcal{A}) = C^\infty(\mathcal{A}) \otimes \left(\Gamma_{lin}(\mathcal{D}, \mathcal{A}) \oplus \Gamma_V (\mathcal{D}, \mathcal{A})\right).\end{equation*}
\end{prop}

\begin{proof} In terms of coordinates $\{x^i, a^i, b^i, c^i\}$ on $\mathcal{D}$ as above, let $\{B_i, C_i\}$ be the frame of sections dual to the fibre coordinates $\{b^i, c^i\}$.  Then a section $X$ is linear if and only if it locally takes the form
\begin{equation}\label{eqn:locallin}
X = f^i(x) B_i + g^i_j (x) a^j C_i,
\end{equation}
and $X$ is vertical if and only if it locally takes the form
\begin{equation}\label{eqn:localvert}
X = f^i(x) C_i.
\end{equation}
In particular, the sections $B_i$ are linear, and the sections $C_i$ are vertical.  The result is immediate.
\end{proof}

The main value of Proposition \ref{prop:vertandlin} is that it allows us to describe the $C^\infty([j]\mathcal{A})$-module $\Gamma([j]_\mathcal{B} \mathcal{D}, [j] \mathcal{A})$ in terms of $\Gamma(\mathcal{D}, \mathcal{A})$, as follows.

\begin{prop}\label{prop:shiftsection}
\begin{enumerate}
\item $\Gamma_{lin}([j]_\mathcal{B} \mathcal{D}, [j] \mathcal{A})$ and $\Gamma_{lin}(\mathcal{D}, \mathcal{A})$ are isomorphic as $C^\infty(\mathcal{M})$-modules.
\item $\Gamma_V ([j]_\mathcal{B} \mathcal{D}, [j] \mathcal{A})$ and $[j]\Gamma_V(\mathcal{D}, \mathcal{A})$ are isomorphic as $C^\infty(\mathcal{M})$-modules.
\item $\Gamma([j]_\mathcal{B} \mathcal{D}, [j] \mathcal{A}) = C^\infty([j] \mathcal{A}) \otimes \left(\Gamma_{lin}(\mathcal{D}, \mathcal{A}) \oplus [j]\Gamma_V (\mathcal{D}, \mathcal{A})\right).$
\end{enumerate}
\end{prop}

\begin{proof}
If $X \in \Gamma_{lin}(\mathcal{D}, \mathcal{A})$ is of degree $k$, then it follows immediately from the fact that $[j]$ is a functor that there is a bundle morphism $[j]X: [j]\mathcal{A} \to [j]_{[-k]\mathcal{B}} [-k]_\mathcal{A} \mathcal{D}$.  Using the identification of Proposition \ref{prop:dvb1}, we may view $[j]X$ as a degree $k$ section of $[j]_\mathcal{B} \mathcal{D}$.  This proves the first statement.

For the second statement, we observe that coordinates on the core of $[-k]_\mathcal{A} \mathcal{D}$ have the same grading as coordinates on the core of $[-k-j]_\mathcal{A} [j]_\mathcal{B} \mathcal{D}$.  Thus, given a vertical degree $k$ section of $\mathcal{D}$, we can obtain a vertical degree $k + j$ section of $[j]_\mathcal{B} \mathcal{D}$ that acts in the same way on the core coordinates.

Lastly, the first and second statements combined with Proposition \ref{prop:vertandlin} give the third statement.
\end{proof}

\begin{ex}
If $\mathcal{E} \to \mathcal{M}$ is a vector bundle, then
\begin{equation}
\xymatrix{T\mathcal{E} \ar[r] \ar[d] & \mathcal{E} \ar[d] \\ T\mathcal{M} \ar[r] & \mathcal{M}}
\end{equation}
is a double vector bundle.  The following proposition may be verified directly in local coordinates.
\end{ex}

\begin{prop} \label{prop:dvb2}Let $\mathcal{E} \to \mathcal{M}$ be a vector bundle. Then $T([j]\mathcal{E}) = [j]_{T\mathcal{M}} T\mathcal{E}$.
\end{prop}

\section{Calculus on supermanifolds}\label{section:calc}

In this section, vector fields and differential forms on supermanifolds are described.  A vector field is defined to be a graded derivation of $C^\infty(\mathcal{M})$.  In parallel with a similar construction in the $\integers_2$ case, the differential forms on a supermanifold are interpreted as functions on the \emph{odd tangent bundle} $[-1]T\mathcal{M}$.  It is shown in Proposition \ref{prop:tangent} that vector fields on $\mathcal{M}$ may be identified with sections of the tangent bundle $T\mathcal{M}$.  Finally, a Cartan-type formula is derived for the de Rham differential; this will provide a starting point in \S\ref{section:superalg} for defining algebroid differentials.

\subsection{Vector fields}\label{sub:vfields}

\begin{dfn} A \emph{vector field} of degree $j$ on a supermanifold $\mathcal{M}$ is a degree $j$ derivation $\phi$ of $C^\infty(\mathcal{M})$, i.e. a linear operator such that, for any homogeneous functions $f, g \in C^\infty(\mathcal{M})$,
\begin{equation} |\phi f| = j + |f| \end{equation}
and 
\begin{equation} \phi(fg) = \phi(f)g + (-1)^{j|f|}f \phi(g). \end{equation}
The space of vector fields on $M$ is denoted $\vect(\mathcal{M})$.
\end{dfn}

\begin{rmk}
Let $\{x^i\}$ be local coordinates on $\mathcal{M}$.  Given a set of functions $\{f^i \in C^\infty(\mathcal{M})\}$, there is (locally) a unique vector field $\psi$ such that $\psi(x^i) = f^i$.  Locally, such a vector field is denoted $\psi = f^i \pdiff{}{x^i}$.  Conversely, any vector field $\phi \in \vect(\mathcal{M})$ is locally determined by the set of functions $\phi(x^i)$.  In particular, if $\mathcal{E} \to \mathcal{M}$ is a vector bundle, then a vector field on $\mathcal{E}$ is completely determined by its action on $C^\infty_{lin}(\mathcal{E})$.
\end{rmk}

\begin{dfn}The \emph{Lie bracket} of two derivations is 
\begin{equation}[\phi,\psi] = \phi\psi - (-1)^{|\phi||\psi|}\psi \phi.\end{equation}
\end{dfn}

\begin{rmk}
If $|\phi| = p$ and $|\psi|= q$, then $[\phi,\psi]$ is a derivation of degree $p+q$.  The resulting bracket gives $\vect(\mathcal{M})$ the structure of a Lie superalgebra; namely, for any homogeneous derivations $\phi$, $\psi$, $\gamma$ and function $f$:
\begin{enumerate}
\item $[\phi,\psi] = (-1)^{1+ |\phi||\psi|}[\psi,\phi]$ (antisymmetry),
\item $[\phi,f\psi]= \phi(f)\psi + (-1)^{|\phi||f|}f[\phi,\psi]$ (Leibniz rule),
\item $(-1)^{|\phi||\gamma|}[\phi,[\psi,\gamma]] + (-1)^{|\psi||\phi|}[\psi,[\gamma,\phi]] + (-1)^{|\gamma||\psi|}[\gamma,[\phi,\psi]] = 0$ (Jacobi identity).
\end{enumerate}
\end{rmk}

\subsection{The category of $Q$-manifolds}\label{sub:qman}

As in the $\integers_2$-graded case, if $\phi$ is an odd degree vector field then $[\phi,\phi] = 2\phi^2$, so the identity $[\phi,\phi] = 0$ is not automatically satisfied.  A degree $1$ vector field which satisfies this identity is said to be \emph{homological}.  

\begin{dfn}[\cite{schwarz}]
A supermanifold equipped with a homological vector field is called a \emph{$Q$-manifold}.  
\end{dfn}

\begin{dfn}Let $\mu: \mathcal{M} \to \mathcal{N}$ be a morphism of supermanifolds, and let $\phi \in \vect(\mathcal{M})$ and $\psi \in \vect(\mathcal{N})$.  Then $\phi$ and $\psi$ are \emph{$\mu$-related} if, for any $f \in C^\infty(\mathcal{N})$, 
\begin{equation}\label{eqn:rel} \mu^* (\psi f) = \phi (\mu^*f).\end{equation}
\end{dfn}

\begin{rmk}Consider the category $\mathbf{VSMan}$ whose objects are pairs $(\mathcal{M}, \phi)$, where $\mathcal{M}$ is a supermanifold and $\phi$ is a vector field on $M$, and where a morphism $\mu: (\mathcal{M},\phi) \to (\mathcal{N},\psi)$ is a supermanifold morphism $\mu: \mathcal{M} \to \mathcal{N}$ such that $\phi$ and $\psi$ are $\mu$-related.  The $Q$-manifolds comprise the class of objects of a full subcategory $\mathbf{QMan}$.
\end{rmk}

\begin{prop}Let $\mu: \mathcal{M} \to \mathcal{N}$ be an embedding, and let $\psi \in \vect(\mathcal{N})$.  There exists a vector field $\phi$ on $\mathcal{M}$ that is $\mu$-related to $\psi$ if and only if, for any $f \in C^\infty(\mathcal{V})$, where $\mathcal{V}$ is an open submanifold of $\mathcal{N}$, $\mu^*f = 0$ implies that $\mu^* \psi f = 0$.  If such a $\phi$ exists, it is unique and denoted $\psi|_{\mathcal{M}}$.
\end{prop}
\begin{proof}
($\Rightarrow$) Suppose that there exists a $\phi$ that is $\mu$-related to $\psi$.  If $\mu^* f = 0$, then by (\ref{eqn:rel}), $\mu^* \psi f = \phi(0) = 0$.

($\Leftarrow$) If $\ker \mu^*$ is invariant under $\psi$, then (\ref{eqn:rel}) well-defines $\phi$ on $\im \mu^*$.  Since $\mu$ is an embedding, $\mu^*$ is (locally) a surjection, and $\phi$ is thus uniquely determined by (\ref{eqn:rel}).  Since $\psi$ is a derivation, $\phi$ is also a derivation.
\end{proof}

\begin{prop}\label{prop:surjrel}Let $\mu: \mathcal{M} \to \mathcal{N}$ be a surjection, and let $\phi \in \vect(\mathcal{M})$.  There exists a vector field $\psi$ on $\mathcal{N}$ that is $\mu$-related to $\phi$ if and only if, for any $f \in C^\infty(\mathcal{V})$, where $\mathcal{V}$ is an open submanifold of $\mathcal{N}$, $\phi(\mu^*f) \in \im \mu^*$.  If such a $\psi$ exists, it is unique and denoted $\mu_* \phi$.
\end{prop}

\begin{proof}
($\Rightarrow$) Suppose that there exists a $\psi$ that is $\mu$-related to $\phi$.  By (\ref{eqn:rel}), $\phi(\mu^* f)$ is in $\im \mu^*$.

($\Leftarrow$) If $\phi(\mu^*f) \in \im \mu^*$, then, since $\mu^*$ is an injection, $\psi$ is uniquely determined by (\ref{eqn:rel}).  Since $\phi$ is a derivation, $\psi$ is also a derivation.
\end{proof}

\begin{rmk}\label{rmk:along}
If $\mu$ is an injection and $\psi \in \vect(\mathcal{N})$ does not restrict to $\mathcal{M}$, then we may still obtain a ``vector field along $\mathcal{M}$ in $\mathcal{N}$'', which is the map $\mu^* \circ \psi : C^\infty(\mathcal{N}) \to C^\infty(\mathcal{M})$.  

If $\mu$ is a surjection and $\phi \in \vect(\mathcal{M})$ does not push forward to a vector field on $\mathcal{N}$, then we will use $\mu_* \phi$ to denote the map $\phi \circ \mu^* : C^\infty(\mathcal{N}) \to C^\infty(\mathcal{M})$.  In analogy with the previous situation, $\mu_* \phi$ may be thought of as a ``vector field along $\mathcal{M}$ in $\mathcal{N}$''.
\end{rmk}

\subsection{(Pseudo)differential forms and the Cartan calculus} \label{sub:diff}

\begin{dfn} \label{dfn:forms} The algebra of differential forms on a supermanifold $\mathcal{M}$ is 
\begin{equation} \Omega(\mathcal{M}) = \bigS\left([1]\Gamma(T^*\mathcal{M})\right) = \Gamma (\bigS[1]T^*\mathcal{M}),
\end{equation}
where $\bigS$ denotes the graded symmetric product.
\end{dfn}

\begin{rmk}
We view the functor $\bigS [1]$ as a supergeometric generalization of the exterior product.  The notation $\bigwedge \left( \Gamma(T^*\mathcal{M}) \right)$ seems more appropriate for the convention used in, e.g. \cite{bdm} and \cite{tuynman}, where differential forms commute or anticommute according to a double-grading consisting of the ``cohomological grading'' and the ``supermanifold grading''.  For example, if $x$ is a degree $0$ coordinate and $\xi$ is a degree $1$ coordinate, then, according to their convention, $dx$ (whose double grading is $(1,0)$)and $d\xi$ (whose double grading is $(1,1)$) will anticommute with each other. On the other hand, our convention uses the total grading, where $d\xi$ commutes with everything since it is of total degree $2$.  Our reason for using the total grading is that it allows us to interpret $\Omega(\mathcal{M})$ as being (up to a suitable closure) the algebra of functions on a supermanifold, as described below\footnote{Deligne \cite{bdm} gives a nice comparison of the two conventions, listing the advantages and disadvantages of each.}.
\end{rmk}

\begin{ex}[Odd tangent bundle]\label{ex:oddtm} Let $M$ be a manifold.  Then the \emph{odd tangent bundle} $[-1]TM$ has fibre coordinates $\dot{x}^i$ of degree $1$ that transform in exactly the same way as the 1-forms $dx^i$.  The function sheaf of $[-1]TM$ is isomorphic to the sheaf of differential forms on $M$.  
\end{ex}

\begin{rmk}
There is a slight difference when extending the construction of Example \ref{ex:oddtm} to the case of supermanifolds; if a supermanifold $\mathcal{M}$ has coordinates of degree $-1$, then $[-1]T\mathcal{M}$ has fibre coordinates of degree $0$ in which functions are not necessarily polynomial.  For this reason, $C^\infty([-1]T\mathcal{M})$ is called the algebra of \emph{pseudodifferential forms} \cite{bl2}.  Since a derivation is completely determined by its action on a set of coordinates, any derivation of $\Omega(\mathcal{M})$ extends uniquely to a derivation of $C^\infty([-1]T\mathcal{M})$, i.e. a vector field on $[-1]T\mathcal{M}$. 
\end{rmk}

\begin{rmk}Let $\mu: \mathcal{M} \to \mathcal{N}$ be a morphism of supermanifolds.  Because the convention in standard differential geometry is to denote by $\mu^*$ both pullback of functions and pullback of differential forms, the ``$[-1]T$'' will often be dropped in what should properly be called ``$[-1]T\mu^*$''.
\end{rmk}

\begin{dfn} The \emph{de Rham differential} $d$ is the degree $1$ vector field on $[-1]T\mathcal{M}$ defined locally by the properties
\begin{align} d(x^i) = \dot{x}^i, && d(\dot{x}^i) = 0.
\end{align}
\end{dfn}

\begin{rmk} Note that (\ref{eqn:tangent}) is equivalent to the identity $\mu^* \circ d = d \circ \mu^*$.\end{rmk}

\begin{rmk}\label{rmk:fanddf}Let $\{x^i\}$ be coordinates on $\mathcal{M}$.  Then $\{x^i, dx^i\}$ may be naturally identified with the coordinates $\{x^i, \dot{x}^i\}$ on $[-1]T\mathcal{M}$.  So any two vector fields on $[-1]T\mathcal{M}$ are equal if they agree on functions of the form $f$ and $df$, where $f \in C^\infty(\mathcal{M})$.
\end{rmk}

\begin{dfn} Let $X \in \Gamma(T\mathcal{M})$.  The \emph{contraction operator} $\iota_X$ is the degree $|X| - 1$ vector field on $[-1]T\mathcal{M}$ such that, for any $\alpha \in C^\infty_{lin}([-1]T\mathcal{M})$,
\begin{equation}\label{eqn:contraction}\iota_X(\alpha) = \langle [-1]X, \alpha \rangle. \end{equation}
\end{dfn}

\begin{dfn} Let $X \in \Gamma(T\mathcal{M})$.  The \emph{Lie derivative} with respect to $X$ is the degree $|X|$ vector field $\lie_X$ on $[-1]T\mathcal{M}$ defined by 
\begin{equation}\label{eqn:lie}\lie_X = [\iota_X, d].\end{equation}
\end{dfn}

\begin{rmk} The Lie derivative has the property that for any $f \in C^\infty(\mathcal{M})$ and $X \in \Gamma(T\mathcal{M})$,
\begin{equation} \lie_{fX} = f \cdot \lie_X + (-1)^{|f| + |X|} df \cdot \iota_X. \end{equation}
\end{rmk}

\begin{rmk} The operators $d$, $\iota_X$, and $\lie_X$ clearly leave $\Omega(\mathcal{M})$ invariant and have cohomological degree $1$, $-1$, and $0$, respectively.
\end{rmk}

\begin{prop}\label{prop:tangent}
The space of sections of the tangent bundle $T\mathcal{M}$ is naturally isomorphic to the space of vector fields on $\mathcal{M}$.
\end{prop}

\begin{proof} Let $X \in \Gamma(T\mathcal{M})$.  For any $f \in C^\infty(\mathcal{M})$, $\lie_X f = \iota_X df \in C^\infty(\mathcal{M})$.  So the map $X \mapsto \lie_X |_{C^\infty(\mathcal{M})}$ is a map from $\Gamma(T\mathcal{M})$ to $\vect(M)$.

Let $\phi \in \vect(\mathcal{M})$.  Then $\phi^V \in \vect([-1]T\mathcal{M})$, defined by the properties
\begin{align}\phi^V (x^i) = 0, && \phi^V (\dot{x}^i) = \phi(x^i), \end{align}
restricts to a homomorphism of right modules $C^\infty_{lin}([-1]T\mathcal{M}) \to C^\infty(\mathcal{M})$.  By Proposition \ref{prop:pairing}, there is a corresponding section $X \in \Gamma(T\mathcal{M})$ such that $\phi^V = \iota_X$, thus giving a map from $\vect(M)$ to $\Gamma(T\mathcal{M})$.

It is straightforward to see that the two maps are inverses of each other.
\end{proof}

\begin{rmk} Since the space of $1$-forms $\Omega^1 (\mathcal{M}) \defequal [1]\Gamma(T^*\mathcal{M})$ is equal to $C^\infty_{lin}([-1]T\mathcal{M})$, (\ref{eqn:contraction}) describes a pairing of $[-1]\vect(\mathcal{M})$ and $\Omega^1 (\mathcal{M})$.
\end{rmk}

\begin{ex}\label{ex:oddcot} Let $M$ be an ordinary manifold, and consider the \emph{odd cotangent bundle} $[-1]T^*M$.  Then $C^\infty([-1]T^*M) = \bigS([1]\vect(M))$ is the algebra of \emph{multivector fields} on $M$ and is usually denoted by $\vect^\bullet (M)$.  The Lie bracket operation on vector fields naturally extends to a graded biderivation $[\cdot, \cdot]: \vect^q(M) \otimes \vect^{q'}(M) \to \vect^{q + q' -1} (M)$, known as the \emph{Schouten bracket}. By interpreting $\vect^\bullet(M)$ as the algebra of functions on $C^\infty([-1]T^*M)$, we may view the Schouten bracket as the degree -1 Poisson bracket associated to the canonical degree $1$ symplectic structure on $[-1]T^*M$ (see \cite{ksm}, \cite{lz}).

If $\mathcal{M}$ is a supermanifold, then $C^\infty([-1]T^*\mathcal{M})$ is the algebra of \emph{pseudomultivector fields}.  Again, the Lie bracket operation on vector fields naturally extends as a biderivation to a bracket on the algebra of pseudomultivector fields, and this bracket may be viewed as the Poisson bracket associated to the canonical degree $1$ symplectic structure on $[-1]T^*\mathcal{M}$.
\end{ex}

\begin{prop}\label{prop:cartancomm}
For any vector fields $X, Y$ on $\mathcal{M}$, the following relations are satisfied:
\begin{align}
[d,d] &=0 \\
[\iota_X, \iota_Y] &= 0 \\
[\lie_X, \iota_Y] &= \iota_{[X,Y]} \\
[\lie_X, \lie_Y] &= \lie_{[X,Y]} \\
[d, \lie_X] &= 0. \label{eqn:liedcom}
\end{align}
\end{prop}

\begin{proof}
The local coordinate expressions of the operators are as follows.  
\begin{equation} d = \dot{x}^i \pdiff{}{x^i}, \end{equation}
and for a homogeneous vector field $X = f^i \pdiff{}{x^i}$,
\begin{align}
\iota_X &= f^i \pdiff{}{\dot{x}^i},\\
\lie_X &= f^i \pdiff{}{x^i} + (-1)^{|X|} df^i \pdiff{}{\dot{x}^i}.
\end{align}
The Lie brackets may be computed directly from the local expressions.
\end{proof}

From the above commutation relations, a Cartan-type formula for the de Rham differential may be derived.  
In particular, from the equations
\begin{align}
\iota_X d &= \lie_X - (-1)^{|X|} d \iota_X \\
\iota_Y \lie_X &= (-1)^{|X|(|Y| - 1)}(\lie_X \iota_Y - \iota_{[X,Y]}),
\end{align}
it follows that for any $p$-form $\omega$ and any vector fields $X_0, \dots, X_p$,
\begin{align}
\begin{split}\label{eqn:cartdiff}
& \iota_{X_p} \cdots \iota_{X_0} d \omega = \sum_{i=0}^p (-1)^{i+ \sum_{j<i}|X_j|} \iota_{X_p} \cdots \lie_{X_i} \cdots \iota_{X_0}\omega \\
&= \sum_{i=0}^p (-1)^{i+\sum_{k<i}|X_k| + |X_i|\left(\sum_{k > i} (|X_k| - 1)\right)} \lie_{X_i}\iota_{X_p} \cdots \widehat{\iota_{X_i}} \cdots \iota_{X_0} \omega \\
& \quad + \sum_{j>i} (-1)^{j+\sum_{k<i}|X_k| + |X_i|\left(\sum_{j \geq k > i}(|X_k| - 1)\right)} \iota_{X_p} \cdots \iota_{[X_i,X_j]} \cdots \widehat{\iota_{X_i}} \cdots \iota_{X_0} \omega.
\end{split}
\end{align}

\begin{prop} \label{prop:iota}Let $\mu:\mathcal{M} \to \mathcal{N}$ be a morphism of supermanifolds and let $X \in \vect(\mathcal{M})$ and $Y \in \vect(\mathcal{N})$ be $\mu$-related vector fields.  Then $\iota_X$ and $\iota_Y$ are $[-1]T\mu$-related.
\end{prop}
\begin{proof}
For any $f \in C^\infty(\mathcal{N})$, $\iota_X \mu^* df = \iota_X d \mu^*f = X \mu^* f = \mu^* Y f = \mu^* \iota_Y df$.  By Remark \ref{rmk:fanddf}, $\iota_X \mu^* = \mu^* \iota_Y$.
\end{proof}

\begin{cor} \label{cor:lie}Let $X \in \vect(\mathcal{M})$ and $Y \in \vect(\mathcal{N})$ be $\mu$-related vector fields.  Then $\lie_X$ and $\lie_Y$ are $[-1]T\mu$-related.
\end{cor}

\section{Superalgebroids}\label{section:superalg}

In this section, superalgebroids are introduced.  It is shown, as in the ordinary case, that an algebroid structure on a vector bundle $\mathcal{A} \to \mathcal{M}$ is equivalent to a differential structure on the algebra of algebroid forms, and the differential structure is interpreted as a homological vector field on the supermanifold $[-1]\mathcal{A}$.  Morphic vector fields, which correspond in the ordinary case to infinitesimal algebroid automorphisms \cite{mac-xu}, are defined.  With superalgebroids, the possibility of morphic vector fields of nonzero degree arises.  In particular, a superalgebroid equipped with a homological morphic vector field is called a \emph{$Q$-algebroid}.  A key example of a $Q$-algebroid is $[-1]T\mathcal{A}$.  The final part of the section describes the structure of $[-1]T\mathcal{A}$ explicitly.

\subsection{Definition and examples}

\begin{dfn}
Let $\mathcal{M}$ be a supermanifold and $\mathcal{A}$ a supervector bundle over $\mathcal{M}$.  A \emph{superalgebroid} structure on $\mathcal{A}$ is a degree $0$ bundle map $\rho: \mathcal{A} \to T\mathcal{M}$ (the \emph{anchor}) and a bilinear operation $[\cdot,\cdot]: \Gamma (\mathcal{A})  \times \Gamma (\mathcal{A})  \to \Gamma (\mathcal{A}) $ such that
\begin{enumerate}
\item $[\cdot,\cdot]$ gives $\Gamma(\mathcal{A})$ the structure of a Lie superalgebra,
\item $[X, fY] = \lie_{\rho(X)}f \cdot Y + (-1)^{|X||f|}f[X,Y]$ for any $X, Y \in \Gamma (\mathcal{A}) $, $f \in C^\infty(\mathcal{M})$ (Leibniz rule), 
\end{enumerate}
where, by abuse of notation, the symbol $\rho$ is also used to denote the module morphism $\Gamma(\rho): \Gamma(\mathcal{A}) \to \vect(\mathcal{M})$.
\end{dfn}

\begin{rmk} As in the ordinary case, the property $\rho([X,Y]) = [\rho(X),\rho(Y)]$, known as the \emph{anchor identity}, follows from the Leibniz rule and Jacobi identity.
\end{rmk}

\begin{ex}[Tangent bundle]
The tangent bundle $T\mathcal{M}$ is a superalgebroid with $\rho = id$ and the Lie bracket.  Similarly, any integrable subbundle of $T\mathcal{M}$ is a superalgebroid.
\end{ex}

\begin{ex}[Action algebroid]
Let $\mathfrak{g}$ be a Lie (super)algebra and let $a: \mathfrak{g} \to \vect(\mathcal{M})$ be a Lie algebra homomorphism.  Then the trivial bundle $\mathcal{M} \times \mathfrak{g} \to \mathcal{M}$ has an algebroid structure, as follows.  Since $\Gamma(\mathcal{M} \times \mathfrak{g}) = C^\infty(\mathcal{M}) \otimes \mathfrak{g}$, the map $id \otimes a: \Gamma(\mathcal{M} \times \mathfrak{g}) \to C^\infty(\mathcal{M}) \otimes \vect(\mathcal{M})$, followed by the natural map $C^\infty(\mathcal{M}) \otimes \vect(\mathcal{M}) \to \vect(\mathcal{M})$, defines the anchor map $\mathcal{M} \times \mathfrak{g} \to T\mathcal{M}$.  The bracket arises from extending the Lie bracket on $\mathfrak{g}$ by the Leibniz rule.
\end{ex}

\begin{ex}[Cotangent bundle]\label{ex:poiss}
Let $\mathcal{M}$ be a supermanifold equipped with a degree $0$ Poisson structure, i.e. a degree $0$ bivector $\pi$ (which may also be viewed as a degree $2$ element of $\vect^2(\mathcal{M})$; see Example \ref{ex:oddcot}) satisfying $[\pi, \pi] = 0$.  Then an algebroid structure can be defined on the cotangent bundle $T^*\mathcal{M}$ in essentially the same way as in the ordinary case \cite{cdw} \cite{w:sg}, where the anchor $\pi^\sharp : \Omega^1(\mathcal{M}) \to \vect^1 (\mathcal{M})$ is defined by the property 
\begin{equation}\label{eqn:poissanch}
\iota_{\pi^\sharp \omega} \eta = \iota_\pi (\omega \wedge \eta),\end{equation}
and the bracket is given by the Koszul bracket
\begin{equation}\label{eqn:poissbrack}
[\omega, \eta]_\pi = \lie_{\pi^\sharp \omega} \eta - \lie_{\pi^\sharp \eta} \omega - d(\iota_\pi(\omega \wedge\eta)). \end{equation}

More generally, a degree $k$ Poisson structure on $\mathcal{M}$ is equivalent to an algebroid structure on $[k]T^*\mathcal{M}$, as follows\footnote{An algebroid structure on $[k]T^*\mathcal{M}$ may be thought of as a ``degree $k$ algebroid structure'' on $T^*M$.  This extends the notion of a degree $k$ Lie algebra.}.  Let the ``bivector field'' $\pi$ be a degree $2-k$ element of $\bigS^2([1-k]\vect(\mathcal{M}))$.  A degree $k$ Poisson bracket on $C^\infty(\mathcal{M})$ may be obtained from $\pi$ as the derived bracket
\begin{equation} \label{eqn:derived} \{f,g\} = \left[ [f,\pi]_k,g \right]_k, \end{equation}
where $[\cdot,\cdot]_k$ is the degree $k-1$ Schouten bracket on $C^\infty([k-1]T^*\mathcal{M})$.  As in the degree $0$ case, it may be shown that the Jacobi identity for the bracket (\ref{eqn:derived}) is equivalent to the condition $[\pi,\pi]_k = 0$.

We may define an algebroid structure on $[k]T^*\mathcal{M}$ by appropriately interpreting the formulae (\ref{eqn:poissanch}) and (\ref{eqn:poissbrack}).  It is simpler, however, to describe the algebroid structure via the $Q$-manifold approach of \S\ref{sub:algq}.  This is done in Example \ref{ex:oddpoiss}.  
\end{ex}

\subsection{Algebroids and $Q$-manifolds}\label{sub:algq}

As in the ordinary case, an algebroid structure on $\mathcal{A}$ induces a differential $d_{\mathcal{A}} : \bigS^\bullet [1]\Gamma (\mathcal{A}^*) \to \bigS^{\bullet + 1} [1]\Gamma (\mathcal{A}^*)$, defined as follows:
\begin{align}\label{eqn-diff}
\begin{split}
& \iota_{X_p} \cdots \iota_{X_0} d_{\mathcal{A}} \omega = \\
&\quad \sum_{i=0}^p (-1)^{i+\sum_{k<i}|X_k| + |X_i|\left(\sum_{k > i} (|X_k| - 1)\right)} \lie_{\rho(X_i)} \iota_{X_p} \cdots \hat{\iota_{X_i}} \cdots \iota_{X_0} \omega \\
& \quad + \sum_{j>i} (-1)^{j+\sum_{k<i}|X_k| + |X_i|\left(\sum_{j \geq k > i}(|X_k| - 1)\right)} \iota_{X_p} \cdots \iota_{[X_i,X_j]} \cdots \hat{\iota_{X_i}} \cdots \iota_{X_0} \omega
\end{split}
\end{align}
for $\omega \in \bigS^p [1]\Gamma (\mathcal{A}^*)$ and any $X_0,\dots X_p \in \Gamma (\mathcal{A}) $.  Equation (\ref{eqn-diff}) is a generalization of equation (\ref{eqn:cartdiff}) for the de Rham differential.

\begin{prop}\label{prop-diff}
$d_{\mathcal{A}}$ satisfies the following properties:
\begin{enumerate}
\item $d_{\mathcal{A}}$ is a degree $1$ derivation of $\bigS [1]\Gamma (\mathcal{A}^*)$
\item $d_{\mathcal{A}}^2=0$
\end{enumerate}
\end{prop}

\begin{proof} 
Let $\{x^i\}$ be local coordinates on $\mathcal{M}$, and let $\{X_\alpha\}$ be a frame of sections of $\mathcal{A}$ where $|X_\alpha| = p_\alpha$.  Then the anchor may be described locally in terms of functions $\rho_\alpha^i$ such that $\rho (X_\alpha) = \rho_\alpha^i \pdiff{}{x^i}$.  Similarly, the bracket may be described in terms of the \emph{structure functions} $c_{\alpha \beta}^\gamma$, where $[X_\alpha, X_\beta] = c_{\alpha \beta}^\gamma X_\gamma$.  

If $\{\lambda^\alpha\}$ are the fibre coordinates dual to $\{X_\alpha\}$, the local description of $d_{\mathcal{A}}$ is
\begin{equation} \label{eqn:localgd}
d_{\mathcal{A}} = \lambda^\alpha \rho_\alpha^i \pdiff{}{x^i} - (-1)^{p_\alpha (p_\beta - 1)}\frac{1}{2} \lambda^\alpha \lambda^\beta c_{\alpha \beta}^\gamma \pdiff{}{\lambda^\gamma}.
\end{equation}
It may be verified by a direct computation that $d_{\mathcal{A}}$, defined locally in this manner, satisfies (\ref{eqn-diff}).

Since $d_{\mathcal{A}}$ is a derivation, it is only necessary to check that the statement holds on $C^\infty(\mathcal{M})$ and $\Gamma(\mathcal{A}^*)$.  Let $f \in C^\infty(\mathcal{M})$.  Then 
\begin{equation}
\begin{split}
\iota_Y \iota_X d_{\mathcal{A}}^2 f &= (-1)^{|X|(|Y|-1)} \left[\rho(X) \iota_Y d_{\mathcal{A}} f - (-1)^{|X||Y|} \rho(Y)\iota_X d_{\mathcal{A}} f - \iota_{[X,Y]} d_{\mathcal{A}} f \right] \\
&= (-1)^{|X|(|Y|-1)} \left[\rho(X) \rho(Y) f - (-1)^{|X||Y|} \rho(Y)\rho(X) f - \rho([X,Y] f \right] \\
&= (-1)^{|X|(|Y|-1)} \left([\rho(X), \rho(Y)] f - \rho([X,Y]) f \right),
\end{split}
\end{equation}
which vanishes by the anchor identity.  A similar (though much more tedious) calculation shows that, for $\omega \in \Gamma (\mathcal{A}^*)$, $\iota_Z \iota_Y \iota_X d_{\mathcal{A}}^2 \omega = 0$ as a result of the anchor and Jacobi identities.
\end{proof}

The most important feature of $d_{\mathcal{A}}$ is that it completely encodes the algebroid structure of $\mathcal{A}$ in the following sense.

\begin{lemma}\label{lemma-deriv}
Let $\mathcal{A} \to \mathcal{M}$ be a supervector bundle.  Given a degree $1$ derivation $d_{\mathcal{A}}: \bigS^\bullet [1]\Gamma (\mathcal{A}^*) \to \bigS^{\bullet + 1} [1]\Gamma (\mathcal{A}^*)$, there is a unique bundle map $\rho: \mathcal{A} \to T\mathcal{M}$ and a unique graded skew-symmetric bilinear map $[\cdot,\cdot]: \Gamma \mathcal{A} \otimes \Gamma \mathcal{A} \to \Gamma \mathcal{A} $, satifying the Leibniz rule, such that (\ref{eqn-diff}) holds.  The Jacobi identity is satisfied if and only if $d_{\mathcal{A}}^2 = 0$.
\end{lemma}

\begin{proof}
Applying (\ref{eqn-diff}) to $f \in C^\infty(\mathcal{M})$ and $\omega \in \Gamma (\mathcal{A}^*)$, the following equations are obtained which uniquely define $\rho$ and $[\cdot,\cdot]$:
\begin{align}\label{eqn:anchor}
\rho(X)f &= \iota_X d_{\mathcal{A}} f = [\iota_X, d_{\mathcal{A}}] f, \\
\begin{split}\label{eqn:bracket}
\iota_{[X,Y]} \omega &=  \rho(X) \iota_Y \omega - (-1)^{|X||Y|} \rho(Y)\iota_X \omega - (-1)^{|X|(|Y|-1)}\iota_Y \iota_X d_{\mathcal{A}} \omega \\
&= \left[ [\iota_X, d_{\mathcal{A}}], \iota_Y \right] \omega.
\end{split}
\end{align}
The last part of the proposition follows from the proof of Proposition \ref{prop-diff}.
\end{proof}

\begin{rmk}Generalizing the role that $[-1]T\mathcal{M}$ plays in the theory of differential forms, $C^\infty([-1]\mathcal{A})$ may be called the algebra of \emph{algebroid pseudoforms}.  The differential $d_{\mathcal{A}}$ extends uniquely to a homological vector field on $[-1]\mathcal{A}$.  The Lie derivative with respect to a section $X$ may be defined as in (\ref{eqn:lie}), and the commutation identities of Proposition \ref{prop:cartancomm} hold.
\end{rmk}

Morphisms of superalgebroids are defined in terms of the differentials; specifically, if $\mathcal{A} \to \mathcal{M}$ and $\mathcal{A}' \to \mathcal{N}$ are superalgebroids, then a (degree $0$) bundle map $\tau: \mathcal{A} \to \mathcal{A}'$ is an algebroid morphism if $d_{\mathcal{A}}$ and $d_{\mathcal{A}'}$ are $[-1]\tau$-related, or, in other words, if $[-1]\tau$ is a morphism of $Q$-manifolds.

We have thus generalized to the case of superalgebroids an important result of \cite{aksz} and \cite{vaintrob}:
\begin{thm}\label{thm:minusone}
The functor $[-1]$ is a functor from the category of superalgebroids to the category of $Q$-manifolds.
\end{thm}

\begin{rmk}\label{rmk:anti}The $Q$-manifolds in the image of the $[-1]$ functor take a specific form; they are vector bundles such that their homological vector fields are of degree $1$ in the fibre coordinates.  Voronov \cite{voronov} calls such objects \emph{anti-algebroids}.  One can recover an algebroid by applying the $[1]$ functor to an anti-algebroid.
\end{rmk}

\begin{ex}[Odd cotangent bundles] \label{ex:oddpoiss}  Let $\pi$ be a degree $2-k$ element of $\bigS^2([1-k]\vect(\mathcal{M}))$ satisfying $[\pi,\pi]_k = 0$.  As in Example \ref{ex:poiss}, $\pi$ determines a degree $k$ Poisson bracket $\{\cdot, \cdot\}$ on $C^\infty(\mathcal{M})$.  On the other hand, the operator $d_\pi \defequal [\pi, \cdot]_k$ is a homological vector field on $[k-1]T^*\mathcal{M}$, giving $[k-1]T^*\mathcal{M}$ the structure of an anti-algebroid.  Thus there is an algebroid structure on $[k]T^*\mathcal{M}$, where the anchor $\pi^\sharp$ and the bracket $[\cdot,\cdot]_\pi$ are defined by (\ref{eqn:anchor}) and (\ref{eqn:bracket}).  In particular, if we denote by $d_n$ the ``degree $n$ de Rham operator'', then $\pi^\sharp$ and $[\cdot,\cdot]_\pi$ may be shown to satisfy the properties
\begin{align}
\pi^\sharp (d_n f)(h) &= \{f,h\}, \\
[d_n f, d_n g]_\pi &= d_n \{f,g\},
\end{align}
for any $f$, $g$, $h \in C^\infty(\mathcal{M})$.
\end{ex}

\subsection{Morphic vector fields}

In this section, we define morphic vector fields on superalgebroids and describe their basic properties.  Morphic vector fields correspond, in the ordinary case, to infinitesimal algebroid automorphisms \cite{mac-xu}.  

Let $\mathcal{E} \stackrel{\pi}{\to} \mathcal{M}$ be a vector bundle.  

\begin{dfn}\label{dfn:linvf}A vector field $\Xi$ on $\mathcal{E}$ is \emph{linear} if, for every $\alpha \in C^\infty_{lin}(\mathcal{E})$, $\Xi(\alpha) \in C^\infty_{lin}(\mathcal{E})$.  The space of linear vector fields is denoted $\vect_{lin}(\mathcal{E})$.
\end{dfn}

\begin{rmk}It follows from the derivation property that, if $\Xi$ is a linear vector field on $\mathcal{E}$, then $C^\infty(\mathcal{M})$ is invariant under the action of $\Xi$ as a derivation.  By Proposition \ref{prop:surjrel}, there exists a unique vector field $\phi$ on $\mathcal{M}$ that is $\pi$-related to $\Xi$.  We call $\phi$ the \emph{base vector field}.
\end{rmk} 

\begin{prop}\label{prop:liniso} For all $j$, $\vect_{lin}(\mathcal{E})$ and $\vect_{lin}([-j]\mathcal{E})$ are naturally isomorphic to each other as $C^\infty(\mathcal{M})$-modules.
\end{prop}

\begin{proof}
Let $\Xi \in \vect_{lin}(\mathcal{E})$.  Using the left $C^\infty(\mathcal{M})$-module isomorphism $C^\infty_{lin}(\mathcal{E}) \to C^\infty_{lin}([-j]\mathcal{E})$, $\alpha \mapsto \hat{\alpha}$, we may define a linear operator $\Xi_j$ on $C^\infty_{lin}([-j]\mathcal{E})$ by the property
\begin{equation}\label{eqn:hatmap}
\Xi_j (\hat{\alpha}) = \widehat{\Xi(\alpha)}.
\end{equation}
Using (\ref{eqn:hatmap}), the fact that $\widehat{f \alpha} = f\hat{\alpha}$, and the fact that $\Xi$ is a linear vector field, it is a simple calculation to see that
\begin{equation}\label{eqn:opdev}
\Xi_j (f \hat{\alpha}) = \phi (f) \hat{\alpha} + (-1)^{|\Xi||f|} f \Xi_j (\hat{\alpha}).
\end{equation}
Equation (\ref{eqn:opdev}) implies that $\Xi_j$ may be extended as a derivation to a linear vector field on $[-j]\mathcal{E}$ with the same degree and the same base vector field as $\Xi$.

It is immediate from the definition that $(f\Xi)_j = f \Xi_j$.  Since the map $\Xi \mapsto \Xi_j$ is a degree $0$ map, the fact that it is a left module isomorphism implies that it is also a right module isomorphism.
\end{proof}

\begin{rmk} In local coordinates, if we write 
\[\Xi = \phi^i (x) \pdiff{}{x^i} + \varphi^a_b(x) \xi^b \pdiff{}{\xi^a},\]
then the local description of $\Xi_j$ is simply
\[\Xi_j = \phi^i (x) \pdiff{}{x^i} + \varphi^a_b(x) \hat{\xi}^b \pdiff{}{\hat{\xi}^a}.\]
\end{rmk}

\begin{prop}\label{prop:hathom}
The map $\Xi \mapsto \Xi_j$ is a Lie algebra homomorphism.
\end{prop}

\begin{proof}
Since $\Xi_j$ is defined by the property (\ref{eqn:hatmap}), we have that
\begin{equation}
\begin{split}
[\Xi_j, \Xi'_j] (\hat{\alpha}) &= \Xi_j \Xi_j^\prime (\hat{\alpha}) - (-1)^{|\Xi||\Xi'|}\Xi_j^\prime \Xi_j (\hat{\alpha}) \\
&= \Xi_j \left(\widehat{\Xi^\prime (\alpha)}\right) - (-1)^{|\Xi||\Xi'|}\Xi^\prime_j \left(\widehat{\Xi (\alpha)}\right) \\
&= \widehat{\Xi \Xi^\prime (\alpha)} - (-1)^{|\Xi||\Xi'|}\widehat{\Xi^\prime \Xi (\alpha)} \\
&= \widehat{[\Xi, \Xi'] (\alpha)},
\end{split}
\end{equation}
from which it follows that $[\Xi_j, \Xi'_j] = [\Xi,\Xi']_j$.
\end{proof}

\begin{dfn} Let $\mathcal{A} \to \mathcal{M}$ be an algebroid.  A linear vector field $\Xi$ on $\mathcal{A}$ is called \emph{morphic} if $[d_\mathcal{A}, \Xi_1] = 0$.
\end{dfn}

\begin{rmk}When there is no possibility of ambiguity, $\Xi_1 \in \vect([-1]\mathcal{A})$ will also be called morphic if $\Xi$ is morphic.
\end{rmk}

\begin{rmk} For any $X \in \Gamma(\mathcal{A})$, the vector field $\lie_X \in \vect([-1]\mathcal{A})$ is morphic (see (\ref{eqn:liedcom})).  When $\mathcal{A} = T\mathcal{M}$, every morphic vector field is of the form $\lie_X$ for some $X \in \vect(\mathcal{M})$.  However, an algebroid may in general have ``outer symmetries'' that do not arise from sections.  For example, consider an algebroid whose bracket is always $0$.  For such an algebroid, $d_{\mathcal{A}}=0$, so any linear vector field is morphic.  However, in this case $\lie_X=0$ for all $X \in \Gamma(\mathcal{A})$.
\end{rmk}

\begin{dfn} Let $\Xi$ be a morphic vector field on an algebroid $\mathcal{A}$.  Then the operator $D_\Xi$ on $\Gamma(\mathcal{A})$ is defined by the equation
\begin{equation}\label{eqn:morph}
\iota_{D_\Xi X} = [\Xi_1, \iota_X].
\end{equation}
\end{dfn}

\begin{prop} Let $\mathcal{A} \to \mathcal{M}$ be an algebroid.  Then
\begin{enumerate}
\item The space of morphic vector fields on $\mathcal{A}$ is closed under the Lie bracket.
\item If $\Xi$ is a morphic vector field and $X \in \Gamma(\mathcal{A})$ then $[\Xi_1, \lie_X] = \lie_{D_\Xi X}$.
\end{enumerate}
\end{prop}

\begin{proof}
By the Jacobi identity, it is clear that if $\Xi$ and $\Xi'$ are linear vector fields such that $[\Xi_1, d_{\mathcal{A}}] = [\Xi^\prime_1, d_{\mathcal{A}}] = 0$, then $\left[ [\Xi_1, \Xi^\prime_1], d_{\mathcal{A}}\right] = 0$.  Using Proposition \ref{prop:hathom}, we conclude that $[\Xi, \Xi']$ is a morphic vector field.

For the second part of the proposition, the Jacobi identity and the fact that $[\Xi_1, d_{\mathcal{A}}] = 0$ imply that 
\[ \left[\Xi_1, [\iota_X, d_{\mathcal{A}}]\right] = \left[[\Xi_1, \iota_X], d_{\mathcal{A}},\right] \]
or, simply, $[\Xi_1, \lie_X] = \lie_{D_\Xi X}$.
\end{proof}

\begin{prop}\label{prop:morphop}Let $\Xi$ be a morphic vector field with base vector field $\phi$.  Then 
\begin{enumerate}
\item For any $X \in \Gamma(\mathcal{A})$, $\rho(D_\Xi X) = [\phi, \rho(X)]$,
\item $D_\Xi$ is a derivation with respect to the Lie bracket.
\end{enumerate}
\end{prop}
\begin{proof}
If $f \in C^\infty(\mathcal{M})$, then by (\ref{eqn:morph}), 
\begin{equation}\begin{split}
\rho(D_\Xi X) f &= \iota_{D_\Xi X} d_{\mathcal{A}} f = [\Xi_1, \iota_X]d_{\mathcal{A}} f  \\
&= \phi \rho(X) f - (-1)^{|\Xi|(|X|-1)} \iota_X \Xi_1 d_{\mathcal{A}} f \\
&= \phi \rho(X) f - (-1)^{|\Xi||X|} \rho(X) \phi f \\
&= [\phi, \rho(X)].
\end{split}\end{equation}
This proves the first statement.

The second statement follows from the Jacobi identity, as follows.
\begin{equation}\begin{split}
\iota_{D_\Xi [X,Y]} &= \left[\Xi_1, [\lie_X, \iota_Y]\right]\\
 &= \left[ [\Xi_1, \lie_X], \iota_Y \right] + (-1)^{|\Xi||X|} \left[\lie_X, [\Xi_1, \iota_Y] \right] \\
 &= \iota_{[D_\Xi X, Y]} + (-1)^{|\Xi||X|} \iota_{[X, D_\Xi Y]}.\qedhere
\end{split}\end{equation}
\end{proof}

\begin{rmk}From (\ref{eqn:morph}), the following formula may be derived:
\begin{equation}\label{eqn:liemorph}
\begin{split}
\iota_{X_q} \cdots \iota_{X_1} \Xi_1 \omega =& 
(-1)^{|\Xi|\sum_{k=1}^q (|X_k| - 1)} \phi(\iota_{X_q} \cdots \iota_{X_1} \omega)\\
 &- \sum_{i=1}^q (-1)^{|\Xi|\sum_{k=1}^i (|X_k| - 1)}\iota_{X_q} \cdots \iota_{X_{i+1}} \iota_{D_\Xi X_i} \iota_{X_{i-1}} \cdots \iota_{X_1} \omega
 \end{split}
 \end{equation}
\end{rmk}

\begin{lemma}\label{lem:morphderiv}
Let $D$ be a linear operator on $\Gamma (\mathcal{A})$ and let $\phi \in \vect(\mathcal{M})$ such that
\begin{enumerate} 
\item $D (fX) = \phi(f)X + (-1)^{|f||D|}f D(X)$ for all $f \in C^\infty(\mathcal{M})$ and $X \in \Gamma(\mathcal{A})$, and
\item $D$ is a derivation of the Lie bracket.
\end{enumerate}
Then, for any $X \in \Gamma(\mathcal{A})$, $\rho(D(X)) = [\phi, \rho(X)]$.
\end{lemma}
\begin{proof}Using the Leibniz rule, we have that for any $X, Y \in \Gamma(\mathcal{A})$ and $f \in C^\infty(\mathcal{M})$,
\begin{equation}\begin{split}\label{morph1}
D \left( [X, fY] \right) =& D\left(\rho(X)(f) \cdot Y + (-1)^{|X||f|} f[X,Y] \right) \\
=& \phi \circ \rho(X)(f) \cdot Y + (-1)^{|D|(|X| + |f|)} \rho(X)(f) \cdot D(Y) + (-1)^{|X||f|} \phi(f) [X,Y] \\
&+ (-1)^{(|D| + |X|)|f|} f[D(X), Y] + (-1)^{(|D| + |X|)|f| + |D||X|} f[X,D(Y)].
\end{split}\end{equation}
On the other hand, 
\begin{equation}\begin{split}\label{morph2}
D \left( [X, fY] \right) =& [D(X),fY] + (-1)^{|D||X|} [X,D(fY)] \\
=& [D(X),fY] + (-1)^{|D||X|}[X, \phi(f)\cdot Y + (-1)^{|D||f|} f D(Y)] \\
=& \rho(D(X))(f) \cdot Y + (-1)^{(|D| + |X|)|f|} f[D(X),Y] + (-1)^{|D||X|} \rho(X) \circ \phi(f) \cdot Y \\
&+ (-1)^{|X||f|} \phi(f)[X,Y] + (-1)^{|D|(|X| + |f|)}\rho(X)(f) \cdot D(Y) \\
&+ (-1)^{|D|(|X|+|f|) + |X||f|} f[X,D(Y)].
\end{split}\end{equation}
After equating the results of (\ref{morph1}) and (\ref{morph2}) and cancelling terms, we obtain the equation
\begin{equation}
\phi \circ \rho(X)(f) \cdot Y = \rho(D(X))(f) \cdot Y + (-1)^{|D||X|} \rho(X) \circ \phi(f) \cdot Y,
\end{equation}
which may be written more simply as
\begin{equation}\label{morph3}
[\phi, \rho(X)] (f) \cdot Y = \rho(D(X))(f) \cdot Y.
\end{equation} 
Since (\ref{morph3}) holds for all $f$ and $Y$, we conclude that $[\phi, \rho(X)] = \rho(D(X))$.
\end{proof}

\begin{prop}\label{prop:morphderiv}
Let $D$ be a linear operator on $\Gamma (\mathcal{A})$ satisfying the hypotheses of Lemma \ref{lem:morphderiv}.  Then there exists a morphic vector field $\Xi \in \vect(\mathcal{A})$ such that $D = D_\Xi$.
\end{prop}
\begin{proof} Given such a $D$, let $\Xi_1$ be the degree $|D|$ linear operator on $C^\infty_{lin}([-1]\mathcal{A})$ defined by the property
\begin{equation} \iota_X \Xi_1(\alpha) = \phi(\iota_X \alpha) - (-1)^{|D|(|X|-1)} \iota_{D(X)} \alpha,
\end{equation}
for all $X \in \Gamma(\mathcal{A})$ and $\alpha \in C^\infty_{lin}([-1]\mathcal{A})$.  Since, for any $f \in C^\infty(\mathcal{M})$, 
\begin{equation}
\iota_X \Xi_1(f\alpha) = \iota_X \left( \phi(f) \alpha - (-1)^{|D||f|} f \Xi_1(\alpha) \right),
\end{equation}
we may extend $\Xi_1$ uniquely to a degree $|D|$ linear vector field on $[-1]\mathcal{A}$ with base vector field $\phi$.

To show that $\Xi$ is morphic, it is sufficient to check that $[d_{\mathcal{A}}, \Xi_1] (\alpha) = 0$ for all $\alpha \in C^\infty_{lin}([-1]\mathcal{A})$.  A (somewhat long) calculation using (\ref{eqn-diff}) and (\ref{eqn:liemorph}) reveals that, for any $X, Y \in \Gamma(\mathcal{A})$,
\begin{equation}\begin{split}
\iota_Y \iota_X [d_{\mathcal{A}}, \Xi_1] \alpha = &\pm \left( \rho(D(X)) - \phi \circ \rho(X) + (-1)^{|D||X|}\rho(X) \circ \phi \right) \iota_Y \alpha \\
& \pm  \left( \rho(D(Y)) - \phi \circ \rho(Y) + (-1)^{|D||Y|}\rho(Y) \circ \phi \right) \iota_X \alpha \\
& \pm \left( \iota_{D[X,Y]} - \iota_{[D(X),Y]} - (-1)^{|X||D|} \iota_{[X, D(Y)]} \right) \alpha,
\end{split}
\end{equation}
where the signs for each line are omitted.  The first two lines vanish by Lemma \ref{lem:morphderiv}, and the final line vanishes because $D$ is by assumption a derivation of the Lie bracket.  It follows that $[d_{\mathcal{A}}, \Xi_1] = 0$.
\end{proof}

\begin{dfn} A Lie algebroid equipped with a homological morphic vector field $\Xi$ is called a \emph{Q-algebroid}.
\end{dfn}

\subsection{The $Q$-algebroid structure of $[-1]T\mathcal{A}$}\label{sec:oneta}

Let $\mathcal{A} \to \mathcal{M}$ be an algebroid.  Then
\begin{equation}
\xymatrix{T\mathcal{A} \ar[r] \ar[d] & \mathcal{A} \ar[d] \\ T\mathcal{M} \ar[r] & \mathcal{M}}
\end{equation}
is a double vector bundle.  Applying the $[-1]$ functor to the rows results in the vector bundle $[-1]_\mathcal{A} T\mathcal{A} \to [-1]T\mathcal{M}$.  It will now be shown that this vector bundle naturally has the structure of a $Q$-algebroid.

It is well-known in the ordinary case (see e.g.\ \cite{mac-xu}) that an algebroid $A \to M$ induces an algebroid structure on the \emph{prolongation} $TA \to TM$.  This fact continues to hold in the super case, but it is more natural to describe the algebroid structure of $[-1]T\mathcal{A} \to [-1]T\mathcal{M}$, from which the algebroid structure of $T\mathcal{A}$ may be obtained by applying the $[1]$ functor (see \S\ref{sec:vba}).

The construction relies on the fact that the sections of $[-1]T\mathcal{A}$ are spanned by two types of ``lifts'' of sections of $\mathcal{A}$.  These lifts are analogous to the vertical and complete lifts of Yano and Ishihara \cite{yi}.

\begin{dfn} Let $X \in \Gamma_j(\mathcal{A})$.  The (odd) \emph{complete lift} $X^C$ of $X$ is a degree $j$ section of $[-1]T\mathcal{A}$ defined by
\begin{equation} X^C = [-1]TX: [-1]T\mathcal{M} \to [-1]T\left([-j]\mathcal{A}\right),
\end{equation}
where the target $[-1]_{[-j]\mathcal{A}} T\left([-j]_\mathcal{M}\mathcal{A}\right)$ is identified with $[-j]_{[-1]T\mathcal{M}} \left([-1]_\mathcal{A} T\mathcal{A}\right)$ (see Propositions \ref{prop:dvb1} and \ref{prop:dvb2}). 
\end{dfn}

\begin{rmk} In terms of functions, $X^C$ is uniquely determined by the properties
\begin{align} X^{C*} (\alpha) = X^* (\alpha), && X^{C*} (d\alpha) = d X^* (\alpha), \end{align}
for $\alpha \in C^\infty(\mathcal{A})$.
\end{rmk}

\begin{dfn} Let $X \in \Gamma(\mathcal{A})$.  The (odd) \emph{vertical lift} $X^V$ of $X$ is a section of $[-1]T\mathcal{A}$ defined by the properties
\begin{align} X^{V*} (\alpha) = 0, && X^{V*} (d\alpha) = (-1)^{|X| + |\alpha| + 1} X^* (\alpha), \end{align}
for $\alpha \in C^\infty_{lin} (\mathcal{A})$.
\end{dfn}

\begin{prop} \label{prop:lifts}Let $f \in C^\infty(\mathcal{M})$, $X \in \Gamma(\mathcal{A})$.  Then
\begin{align}
(fX)^C &= f \cdot X^C + (-1)^{|X| + |f|} df \cdot X^V, \\
(fX)^V &= f \cdot X^V.
\end{align}
\end{prop}

\begin{proof}
The identities follow directly from the above definitions.
\end{proof}

\begin{prop} $\Gamma([-1]T\mathcal{A})$ is spanned by the complete and vertical lifts of $\Gamma(\mathcal{A})$.  More precisely, $\Gamma([-1]T\mathcal{A})$ is equal to $C^\infty([-1]T\mathcal{M}) \otimes \left( \left\{X^C\right\} \oplus \left\{X^V\right\}\right)$, modulo the relations of Proposition \ref{prop:lifts}.
\end{prop}

\begin{proof}Let $\{X_\alpha\}$ be a local frame of sections of $\mathcal{A}$ dual to fibre coordinates $\{\lambda^\alpha\}$.  Then it is clear from the definitions that $\{X^C_\alpha, -X^V_\alpha\}$ is a local frame of sections of $[-1]T\mathcal{A}$ dual to the fibre coordinates $\{\lambda^\alpha, \dot{\lambda}^\alpha\}$.
\end{proof}

The algebroid structure of $[-1]T\mathcal{A}$ is as follows.  The anchor $\widetilde{\rho} : [-1]T\mathcal{A} \to T([-1]T\mathcal{M})$ is defined by letting
\begin{align} \label{eqn:otanch}
\widetilde{\rho}(X^C) = \lie_{\rho(X)}, && \widetilde{\rho}(X^V) = \iota_{\rho(X)}, 
\end{align}
and then extending by $C^\infty([-1]T\mathcal{M})$-linearity.  Because the relations of Proposition \ref{prop:lifts} are identical to the relations satisfied by Lie derivative and contraction operators, $\widetilde{\rho}$ is well-defined.  The bracket is defined by letting
\begin{equation}
\begin{split} \label{eqn:otbrack}
\left[X^C, Y^C\right] &= \left[X,Y\right]^C, \\
\left[X^C, Y^V\right] &= \left[X,Y\right]^V, \\
\left[X^V, Y^V\right] &= 0,
\end{split}
\end{equation}
and extending by the Leibniz rule.  

One could now verify directly that the bracket (\ref{eqn:otbrack}) satisfies the Jacobi identity.  However, we will instead describe the vector field $d_{[-1]T\mathcal{A}} \in \vect\left([-1]_{[-1]T\mathcal{M}}\left([-1]T\mathcal{A}\right)\right)$ defined by (\ref{eqn-diff}) and see that $\left(d_{[-1]T\mathcal{A}}\right)^2 = 0$.

By Propositions \ref{prop:dvb1} and \ref{prop:dvb2}, $[-1]_{[-1]T\mathcal{M}}\left([-1]T\mathcal{A}\right)$ may be identified with $[-1]_{[-1]\mathcal{A}} T \left([-1]\mathcal{A}\right)$.  The algebroid structure on $\mathcal{A}$ induces the differential $d_\mathcal{A} \in \vect([-1]\mathcal{A})$ whose Lie derivative operator is $\lie_{d_\mathcal{A}} \in \vect\left([-1]_{[-1]\mathcal{A}}T ([-1]\mathcal{A})\right)$.

\begin{thm}\label{thm:lieda}$d_{[-1]T\mathcal{A}} = \lie_{d_\mathcal{A}}$.
\end{thm}
\begin{proof}
From the local coordinate description (\ref{eqn:localgd}) of $d_\mathcal{A}$, the Lie derivative may be computed to be 
\begin{equation} 
\begin{split}
\lie_{d_\mathcal{A}} = & \lambda^\alpha \rho_\alpha^i \pdiff{}{x^i} 
- (-1)^{p_\alpha (p_\beta - 1)}\frac{1}{2} \lambda^\alpha \lambda^\beta c_{\alpha \beta}^\gamma \pdiff{}{\lambda^\gamma} 
+ (-1)^{p_\alpha} \lambda^\alpha d\rho^i_\alpha \pdiff{}{\dot{x}^i} 
- (-1)^{p_\alpha p_\beta} \lambda^\alpha \dot{\lambda}^\beta c_{\alpha \beta}^{\gamma} \pdiff{}{\dot{\lambda}^\gamma} \\
& - \dot{\lambda}^\alpha \rho_\alpha^i \pdiff{}{\dot{x}^i} 
+ (-1)^{(p_\alpha +1) p_\beta} \frac{1}{2} \lambda^\alpha \lambda^\beta dc_{\alpha \beta}^\gamma \pdiff{}{\dot{\lambda}^\gamma}.
\end{split}
\end{equation}
If $\{X_\alpha\}$ are the sections dual to the fibre coordinates $\{\lambda^\alpha\}$ on $\mathcal{A}$, then $\{X^C_\alpha, - X^V_\alpha\}$ are the sections dual to $\{\lambda^\alpha,\dot{\lambda}^\alpha\}$ on $[-1]T\mathcal{A}$.  Using equations (\ref{eqn:anchor}) and (\ref{eqn:bracket}), it is straightforward to check that (\ref{eqn:otanch}) and (\ref{eqn:otbrack}) are satisfied.
\end{proof}

\begin{cor}
\begin{enumerate}
\item $\left(d_{[-1]T\mathcal{A}}\right)^2 = 0$ and, equivalently, the anchor and Jacobi identities hold on $[-1]T\mathcal{A}$.
\item The de Rham differential $d \in \vect\left([-1]T([-1]\mathcal{A})\right)$ satisfies the equation $\left[d_{[-1]T\mathcal{A}}, d \right] = 0$.
\end{enumerate}
\end{cor}

\begin{cor}
$[-1]T\mathcal{A}$ is a $Q$-algebroid with the algebroid structure described above and the morphic vector field $d$.
\end{cor}
    
    \chapter{Supergroupoids}
    \label{chapter:supergpds}
Just as a Lie groupoid is a groupoid object in the category of manifolds, a supergroupoid is a groupoid object in the category of supermanifolds.  The goal of this chapter is to define supergroupoids and describe the Lie functor from the category of supergroupoids to the category of superalgebroids.

The approach taken here is based on the simplicial structure associated to a groupoid, and all the definitions and proofs rely only on properties of this simplicial structure.  Consequently, although the emphasis is on $\integers$-graded supergroupoids, the results hold equally well for $\integers_2$-graded supergroupoids.  

The simplicial point of view is also suitable for other generalizations.  As is pointed out in Remark \ref{rmk:liecat}, the definition of the Lie functor makes sense on a Lie (super)category.  This result corresponds to the fact that, in the ordinary case, inverses exist near the units of a Lie category, and the submanifold of invertible elements is a Lie groupoid.  We also expect that this approach will be a good starting point for describing the Lie functor from $2$-groupoids to $2$-algebroids \cite{baez-crans}.

When considering the Lie functor, it is natural to ask when it is possible to go in the other direction, specifically, when a superalgebroid can be integrated to a supergroupoid.  This question is not addressed here and would make an interesting subject of future study.

\section{Definitions}

The definition of a supergroupoid uses in an essential way the concept of fibre product.  In this section, fibre products are defined and some basic lemmas are given regarding vector fields on fibre products.  Then supergroupoids are defined, and the associated simplicial structures are described.

\subsection{Fibre product}

\begin{dfn} Let $\mathcal{Z} \subseteq \mathcal{M}$ be an embedded submanifold, and let $\mu: \mathcal{N} \rightarrow \mathcal{M}$ be a morphism of supermanifolds.  Then the \emph{preimage} $\mu^{-1} (\mathcal{Z})$ is the pair $(\mu_0^{-1}(Z), \mathcal{O}_{\mu^{-1}(\mathcal{Z})})$, where $\mathcal{O}_{\mu^{-1}(\mathcal{Z})} \defequal C^\infty(\mathcal{N})/\langle\mu^*I(\mathcal{Z})\rangle)$ (see Remark \ref{rmk:subideal}).
\end{dfn}

\begin{rmk} As in the ordinary case, the preimage of a submanifold is not always a submanifold.  In this case, it is simply a ringed space.  
\end{rmk}

\begin{dfn} \label{dfn:submersion} Let $\mathcal{M}$ and $\mathcal{P}$ be supermanifolds of dimension $\{p_i\}$ and $\{p_i + q_i\}$, respectively.  A map $\pi: \mathcal{P} \to \mathcal{M}$ is a \emph{submersion} if it is locally isomorphic to the canonical projection $\reals^{\{p_i\}} \times \reals^{\{q_i\}} \to \reals^{\{p_i\}}$.
\end{dfn}

\begin{rmk} Manin \cite{manin} gives the following submersion theorem for supermanifolds: A morphism $\pi: \mathcal{P} \to \mathcal{M}$ is a submersion if and only if $\pi_*: \vect(\mathcal{P}) \to \Gamma(\pi^*(T\mathcal{M}))$ is a surjection.
\end{rmk}

Let $\mathcal{P}$ and $\mathcal{P}^\prime$ be supermanifolds equipped with maps $\pi$ and $\pi^\prime$, respectively, to a supermanifold $\mathcal{M}$.  

\begin{dfn} \label{dfn:fibreprod} The \emph{fibre product} $\mathcal{P} \bitimes{\pi}{\pi^\prime} \mathcal{P}^\prime$ is a supermanifold equipped with maps $p_1: \mathcal{P} \bitimes{\pi}{\pi^\prime} \mathcal{P}^\prime \to \mathcal{P}$ and $p_2: \mathcal{P} \bitimes{\pi}{\pi^\prime} \mathcal{P}^\prime \to \mathcal{P}^\prime$ such that
\begin{equation}\label{eqn:fibrediag}
\xymatrix{
\mathcal{P} \bitimes{\pi}{\pi^\prime} \mathcal{P}^\prime \ar[r]^-{p_1} \ar[d]_{p_2} & \mathcal{P} \ar[d]^\pi \\
\mathcal{P}^\prime \ar[r]^{\pi^\prime} & \mathcal{M}}\end{equation}
commutes and satisfying the universal property that any supermanifold $\mathcal{N}$ with maps $\alpha: \mathcal{N} \to \mathcal{P}$ and $\alpha^\prime : \mathcal{N} \to \mathcal{P}^\prime$ such that $\pi \circ \alpha = \pi^\prime \circ \alpha^\prime $ uniquely factors through $\mathcal{P} \bitimes{\pi}{\pi^\prime} \mathcal{P}^\prime$.
\end{dfn}

\begin{rmk} The universal property implies that the fibre product, if it exists, is unique up to canonical isomorphism.  However, even in the ordinary situation, the fibre product does not always exist as a smooth manifold.  As in the ordinary case, there does always exist a ``topological fibre product'' that lies in the category of ringed spaces.
\end{rmk}

Denote by $\Delta_{\mathcal{M}}$ the diagonal map from $\mathcal{M}$ to $\mathcal{M} \times \mathcal{M}$.  If the domain is clear from the context, we will drop the subscript $\mathcal{M}$.  The image $\Delta(\mathcal{M})$ is a submanifold of $\mathcal{M} \times \mathcal{M}$.

\begin{prop} Let $\boldsymbol{\pi} \defequal \pi \times \pi' : \mathcal{P} \times \mathcal{P}' \to \mathcal{M} \times \mathcal{M}$.  Then $\boldsymbol{\pi}^{-1}(\Delta(\mathcal{M}))$ satisfies the conditions of Definition \ref{dfn:fibreprod} in the category of ringed spaces.
\end{prop}
\begin{proof}
Let $\mathcal{N}$ be as in Definition \ref{dfn:fibreprod}.  By the universal property of the Cartesian product, there is a unique map $\alpha \times \alpha^\prime: \mathcal{N} \to \mathcal{P} \times \mathcal{P}'$ such that $\alpha = p_1 \circ (\alpha \times \alpha')$ and $\alpha' = p_2 \circ (\alpha \times \alpha')$.  Then, since $\pi \circ \alpha = \pi^\prime \circ \alpha^\prime$, it follows that $\boldsymbol{\pi}^* I(\Delta(\mathcal{M})) \subseteq \ker (\alpha \times \alpha')^*$; therefore the sheaf morphism $(\alpha \times \alpha')^*$ passes to the quotient sheaf $\mathcal{O}_{\boldsymbol{\pi}^{-1}(\Delta(\mathcal{M}))}$.
\end{proof}

\begin{prop} \label{prop:fpsub} If either $\pi$ or $\pi'$ is a submersion, then $\boldsymbol{\pi}^{-1}(\Delta(\mathcal{M}))$ is a submanifold of $\mathcal{P} \times \mathcal{P}'$.
\end{prop}

\begin{proof} Let $\{x^i\}$ be a set of coordinates on a neighborhood $\mathcal{U} \subseteq \mathcal{M}$.  Then, on the product neighborhood $\mathcal{U} \times \mathcal{U} \subseteq \mathcal{M} \times \mathcal{M}$, we can take coordinates of the form $\{x_1^i, x_2^i\}$.  The ideal corresponding to $\boldsymbol{\pi}^{-1}(\Delta(\mathcal{M}))$ is then generated by the functions $\boldsymbol{\pi}^*( x^i_2 - x^i_1 )$.  
Let $\{z^i\}$ be a set of coordinates on a neighborhood $\mathcal{V} \subseteq \pi^{-1}(\mathcal{U})$, and, assuming that $\pi'$ is a submersion, choose coordinates of the form $\{\pi'^* x^i, y^i\}$ on a neighborhood $\mathcal{V}' \subseteq (\pi')^{-1}(\mathcal{U})$.  Then, taking $\{\overline{x}^i, y^i, z^i\}$, where $\overline{x}^i \defequal \pi'^* x^i - \boldsymbol{\pi}^* x^i_1$, to be coordinates on $\mathcal{V} \times \mathcal{V}'$, we see that, locally, $I(\boldsymbol{\pi}^{-1}(\Delta(\mathcal{M})))$ is generated by the coordinate functions $\overline{x}^i$, and therefore $\boldsymbol{\pi}^{-1}(\Delta(\mathcal{M}))$ is a submanifold of $\mathcal{P} \times \mathcal{P}'$.
\end{proof}

In what follows, it will be assumed that $\pi'$ is a surjective submersion.

\begin{lemma} \label{lemma:fpint}The fibre product satisfies the property that any function $f \in C^\infty(\mathcal{P} \bitimes{\pi}{\pi^\prime} \mathcal{P}^\prime)$ that can be locally written as both $p_1^* g$ and $p_2^*g'$, where $g \in C^\infty(\mathcal{P})$, $g' \in C^\infty(\mathcal{P}')$, can also be written as $p_1^* \pi^* h = p_2^* \pi'^* h$, where $h \in C^\infty(\mathcal{M})$.
\end{lemma}

\begin{proof}
In terms of the local description given in Proposition \ref{prop:fpsub}, it is clear that if $f = p_2^* g'$, then $f$ may be written locally in terms of the coordinates $\{\pi'^* x^i, y^i\}$, where $\{x^i\}$ are coordinates on $\mathcal{M}$.  On the other hand, if $f = p_1^* g$ then it cannot depend on the coordinates $\{y^i\}$.  Thus $f$ only depends on the coordinates on $\mathcal{M}$.
\end{proof}

\begin{lemma}\label{lemma:equal}
Any two vector fields on $\mathcal{P} \bitimes{\pi}{\pi^\prime} \mathcal{P}^\prime$ that agree on $\im p_1^*$ and $\im p_2^*$ are equal.
\end{lemma}

\begin{proof}
Since a vector field is completely determined by its action as a derivation on a set of coordinate functions, the result follows immediately from the local description of Proposition \ref{prop:fpsub}.
\end{proof}

\begin{lemma}\label{lemma:lift}
Let $X$ and $X^\prime$ be vector fields on $\mathcal{P}$ and $\mathcal{P}^\prime$, respectively.  There is a unique vector field $\widetilde{X}$ on $\mathcal{P} \bitimes{\pi}{\pi^\prime} \mathcal{P}^\prime$ that is $p_1$-related to $X$ and $p_2$-related to $X^\prime$ if and only if there exists a vector field $Y$ on $\mathcal{M}$ that is both $\pi$-related to $X$ and $\pi^\prime$-related to $X^\prime$.
\end{lemma}

\begin{proof}
$(\Rightarrow)$. Suppose that such an $\widetilde{X}$ exists.  Then, for any $f \in C^\infty(\mathcal{M})$,
\[\widetilde{X} p_2^* \pi'^* f = p_2^* X' \pi'^* f,\]
while on the other hand, 
\[\widetilde{X} p_1^* \pi^* f = p_1^* X \pi^* f.\] 
By Lemma \ref{lemma:fpint}, it follows that $\widetilde{X} p_2^* \pi'^* f$ may be written as $p_2^* \pi'^* h$ for some $h \in C^\infty(\mathcal{M})$.  By Proposition \ref{prop:surjrel}, there exists a unique vector field $Y$ on $\mathcal{M}$ that is $(\pi' \circ p_2)$-related to $\widetilde{X}$.  Since $\widetilde{X}$ is $p_2$-related to $X'$, $Y$ is $\pi'$-related to $X'$.  Since $\pi' \circ p_2 = \pi \circ p_1$ and $X$ is $p_1$-related to $X$, $Y$ is $\pi$-related to $X$.

$(\Leftarrow)$. Suppose such a $Y$ exists.  We wish to define a vector field $\widetilde{X}$ by the properties 
\begin{align}\label{eqn:lemlift1}\widetilde{X} (p_1^* g) &= p_1^* X (g), \\
\label{eqn:lemlift2}\widetilde{X} (p_2^* g') &= p_2^* X' (g'),
\end{align}
for any $g \in C^\infty(\mathcal{P})$, $g' \in C^\infty(\mathcal{P}')$.  By Lemma \ref{lemma:fpint}, if $f$ can be written as both $f = p_1^* g$ and $f = p_2^* g'$, then $f$ may also be written as $p_1^* \pi^* h = p_2^* \pi'^* h$ for some $h \in C^\infty(\mathcal{M})$.  Then properties (\ref{eqn:lemlift1}) and (\ref{eqn:lemlift2}) both reduce to $\widetilde{X} f = p_1^* \pi^* Y(h) = p_2^* \pi'^* Y(h)$, so $\widetilde{X}$ is well-defined.
\end{proof}

\begin{rmk} Given any two maps $\alpha: \mathcal{N} \to \mathcal{P}$ and $\alpha^\prime : \mathcal{N} \to \mathcal{P}^\prime$, the image of the product map $\alpha \times \alpha^\prime: \mathcal{N} \to \mathcal{P} \times \mathcal{P}'$ lies in the submanifold $\mathcal{P} \bitimes{\pi}{\pi^\prime} \mathcal{P}^\prime$ if and only if $\pi \circ \alpha = \pi^\prime \circ \alpha^\prime$.  If this is the case then the map $\alpha \times \alpha^\prime: \mathcal{N} \to \mathcal{P} \bitimes{\pi}{\pi^\prime} \mathcal{P}^\prime$ will be called \emph{well-defined}.
\end{rmk}

\subsection{Supergroupoids} 

A supergroupoid is simply a groupoid structure in the category of supermanifolds.  A groupoid structure on a supermanifold $\mathcal{G}$ includes, in particular, a pair of surjective submersions $s,t: \mathcal{G} \to \mathcal{M}$.  The fibre product $\mathcal{G} \bitimes{s}{t} \mathcal{G}$ is denoted by $\mathcal{G}^{(2)}$ and in general
\[\mathcal{G}^{(q)} \defequal \mathcal{G} \bitimes{s}{t} \cdots \bitimes{s}{t} \mathcal{G}\]
for $q > 0$, while $\mathcal{G}^{(0)} \defequal \mathcal{M}$.

\begin{dfn}\label{dfn-gpd}
A \emph{Lie supergroupoid} is a pair of supermanifolds $(\mathcal{G}, \mathcal{M})$ equipped with two surjective submersions $s,t : \mathcal{G} \to \mathcal{M}$ and maps $m: \mathcal{G}^{(2)} \to \mathcal{G}$ (multiplication), $e: \mathcal{M} \to \mathcal{G}$ (identity), and $i: \mathcal{G}\to \mathcal{G}$ (inverse) such that the following diagrams are well-defined and commute:

\begin{enumerate}
\item (associativity)
\xymatrix{
		& \mathcal{G}^{(3)} \ar[dl]_{id \times m} \ar[dr]^{m \times id} \\
\mathcal{G}^{(2)} \ar[dr]_m &		& \mathcal{G}^{(2)} \ar[dl]^m \\
		& \mathcal{G}
   }

\item (identity)
\xymatrix{
\mathcal{G} \ar[r]^-{\Delta_0^1} \ar@/_1pc/[rr]_{id} & \mathcal{G}^{(2)} \ar[r]^-m & \mathcal{G}
}
\xymatrix{
\mathcal{G} \ar[r]^-{\Delta_1^1} \ar@/_1pc/[rr]_{id} & \mathcal{G}^{(2)} \ar[r]^-m & \mathcal{G}
}

\item (inverse)
\xymatrix{
\mathcal{G} \ar[r]^-{i_R} \ar[dr]_t & \mathcal{G}^{(2)} \ar[r]^-m & \mathcal{G} \\
		& \mathcal{M} \ar[ur]_e&
}
\xymatrix{
\mathcal{G} \ar[r]^-{i_L} \ar[dr]_s & \mathcal{G}^{(2)} \ar[r]^-m & \mathcal{G} \\
		& \mathcal{M} \ar[ur]_e&
},
\end{enumerate}
where
\begin{align}
\Delta_0^1 &= (e \circ t \times id)\circ \Delta, \\
\Delta_1^1 &= (id \times e \circ s)\circ \Delta, \\
i_R &= (id \times i)\circ \Delta, \\
i_L &= (i \times id) \circ \Delta.
\end{align}
\end{dfn}

\begin{rmk}\label{rmk:liecat}Much of the material which follows does not rely on the existence of the inverse and is valid in the more general context of a Lie supercategory, where the map $i$ (and consequently condition (3)) is dropped and $s$ (or $t$) is not required to be a submersion.  
\end{rmk}

\begin{ex}[Pair groupoid]
Let $\mathcal{M}$ be a supermanifold.  Then the \emph{pair groupoid} of $\mathcal{M}$ is $\mathcal{G} = \mathcal{M} \times \mathcal{M} \arrows \mathcal{M}$, where $s$ and $t$ are projection to the left and right components, respectively.  The multiplication map $m: \mathcal{G}^{(2)} = \mathcal{M} \times \mathcal{M} \times \mathcal{M} \to \mathcal{M} \times \mathcal{M}$ is the map which ``forgets'' the middle component.  The identity map is the diagonal map $\Delta: \mathcal{M} \to \mathcal{M} \times \mathcal{M}$, and the inverse map $i: \mathcal{M} \times \mathcal{M} \to \mathcal{M} \times \mathcal{M}$ exchanges the two components.
\end{ex}

\begin{ex}[Action groupoid]
Let $\Gamma$ be a supergroup that acts (on the right) on a supermanifold $\mathcal{M}$, with action map $s:\mathcal{M} \times \Gamma \to \mathcal{M}$.  Then the \emph{action groupoid} is $\mathcal{G} = \mathcal{M} \times \Gamma \arrows \mathcal{M}$, where $t$ is projection onto $\mathcal{M}$ and $s$ is the given action map.  The multiplication map $m: \mathcal{G}^{(2)} = \mathcal{M} \times \Gamma \times \Gamma \to \mathcal{M} \times \Gamma$ is $id \times \mu$, where $\mu$ is the multiplication map on $\Gamma$.  The identity map $e: \mathcal{M} \cong \mathcal{M} \times \{pt.\} \to \mathcal{M} \times \Gamma$ is $id \times E$, where $E$ is the identity map on $\Gamma$.  The inverse map $i: \mathcal{M} \times \Gamma \to \mathcal{M} \times \Gamma$ is $(s \times (I \circ p_\Gamma)) \circ \Delta$, where $I$ is the inverse map on $\Gamma$ and $p_\Gamma: \mathcal{M} \times \Gamma \to \Gamma$ is projection onto $\Gamma$.
\end{ex}

\begin{dfn}
A \emph{groupoid in the category of vector bundles} is a pair of groupoids, $\Gamma \arrows \mathcal{E}$ and $\mathcal{G} \arrows \mathcal{M}$, such that $\Gamma$ is a vector bundle over $\mathcal{G}$, $\mathcal{E}$ is a vector bundle over $\mathcal{M}$, and all of the structure maps for $\Gamma$ are bundle maps over the corresponding structure maps for $\mathcal{G}$.
\end{dfn}

Groupoids in the category of vector bundles are called $\mathcal{VB}$-groupoids in \cite{mac:dbl}.  If 
\[\xymatrix{\Gamma \ar[r] \ar@<2pt>[d] \ar@<-2pt>[d] & \mathcal{G} \ar@<2pt>[d] \ar@<-2pt>[d]\\
\mathcal{E} \ar[r] & \mathcal{M}
}\]
is a $\mathcal{VB}$-groupoid, the degree-shifting functor $[j]$ may be applied to the left column of the above diagram.  It follows from functoriality that the result is again a groupoid.  Thus we have 

\begin{prop}\label{prop-bundles}
Let $(\Gamma \arrows \mathcal{E}, \mathcal{G} \arrows \mathcal{M})$ be a groupoid in the category of vector bundles.  Then there is an induced groupoid structure on $[j]\Gamma \arrows [j]\mathcal{E}$.
\end{prop}

\begin{ex}
If $\mathcal{G} \arrows \mathcal{M}$ is a groupoid, then its tangent bundle $T\mathcal{G}$ has a groupoid structure with base $T\mathcal{M}$, and the $[-1]$ functor may be applied, resulting in the \emph{odd tangent groupoid} $[-1]T\mathcal{G} \arrows [-1]T\mathcal{M}$.  The groupoid structure on $[-1]T\mathcal{G}$ may be easily described, since the structure maps of $[-1]T\mathcal{G}$ arise from the pullback of differential forms by the structure maps of $\mathcal{G}$.
\end{ex}

\begin{ex}
If $G \arrows M$ is an ordinary groupoid\footnote{We expect a similar result to hold for supergroupoids.},  the cotangent bundle $T^*G$ has a natural groupoid structure \cite{cdw} \cite{pradines2} over the dual $A^*$ of the Lie algebroid of $G$.  The $[-1]$ functor yields the ``odd cotangent groupoid'' $[-1]T^*G \arrows [-1]A^*$.  In fact, $[-1]T^*G$ has the structure of an \emph{odd symplectic groupoid}.
\end{ex}

\begin{dfn} \label{dfn:nerve} Let $\mathcal{G} \arrows \mathcal{M}$ be a groupoid.  The \emph{nerve} $N\mathcal{G}$ of $\mathcal{G}$ is the simplicial (super)manifold \cite{segal} such that $N\mathcal{G}_q = \mathcal{G}^{(q)}$, with face maps $\sigma_0^1 = s$, $\sigma_1^1 = t$, and
\begin{align}
\sigma_0^q &= p_2 \times id \times \cdots \times id, \\
\sigma_i^q &= \underbrace{id \times \cdots \times id}_{i - 1} \times m \times \underbrace{id \times \cdots id}_{q - i - 1}, & 0 < i < q, \\
\sigma_q^q &= id \times \cdots \times id \times p_1,
\end{align}
for $q > 1$, and degeneracy maps $\Delta_0^0 = e$ and 
\begin{align}
\Delta_i^q &= \underbrace{id \times \cdots \times id}_{i} \times \Delta_0^1 \times \underbrace{id \times \cdots id}_{q - i - 1}, & i < q, \label{eqn:degen-left}\\
\Delta_i^q &= \underbrace{id \times \cdots \times id}_{i - 1} \times \Delta_1^1 \times \underbrace{id \times \cdots id}_{q - i}, & 0 < i, \label{eqn:degen-right}
\end{align}
for $q > 0$.
\end{dfn}

\begin{rmk}The nerve of a Lie category is defined in the same manner (see \cite{segal}).
\end{rmk}

\begin{rmk}
The simplicial point of view introduces duplicate notation for several maps. For example, the multiplication map $m$ may also be written as $\sigma_1^2$.  In general, the ``standard'' notation will be used except when ``higher'' simplicial structures are involved.  At these times the reader may need to translate equations between the standard and simplicial notations.
\end{rmk}

\section{Lie Functor}

To describe the Lie functor from groupoids to algebroids, one considers left-invariant vector fields.  In the first part of this section, left-invariant vector fields on a supergroupoid $\mathcal{G} \to \mathcal{M}$ are defined, and a natural identification is established between the space of left-invariant vector fields and the space of vector fields along $\mathcal{M}$ that are tangent to the $t$-fibres.  The left-invariant vector fields form the space of sections of a superalgebroid, thus defining the Lie functor (at the level of objects).  Next, multiplicative vector fields, which are infinitesimal groupoid automorphisms, are described, and it is shown that a multiplicative vector field on a groupoid induces a morphic vector field on its algebroid.  Finally, left-invariant differential forms are described in order to give the dual picture.

All of the results of this section are known in the ordinary case.  However, the simplicial interpretation of the Lie functor, used here to extend the results to the super case, may be interesting even in ordinary contexts.

\subsection{Left-invariant vector fields}

\begin{rmk}
Consider the case where $G$ is an ordinary Lie group, and let $a$ be an element of $G$.  The diffeomorphism $\ell_a : G \to G$ given by left-multiplication by $a$ may be decomposed into the composition of $\sigma_a: G \to G \times G$, $\sigma_a(g) =(a,g)$, and the multiplication map $m$.  A vector field $X$ on $G$ such that $(m \circ \sigma_a)_* X = X$ for all $a$ is left-invariant in the usual sense.

For each $a$, $\sigma_a$ maps $G$ to a different vertical slice of $G \times G$, and considering the image of $X$ under all such maps, we obtain a vector field $\widetilde{X}$ on $G \times G$, determined completely by the conditions that $\widetilde{X} \in \ker p_{1*}$ and $\widetilde{X}$ is $p_2$-related to $X$.  The left-invariance condition is simply that $\widetilde{X}$ is also $m$-related to $X$.  This point of view leads to the following definition.
\end{rmk}

\begin{dfn}\label{dfn-left}
A vector field $X$ on a supergroupoid $\mathcal{G}$ is \emph{left-invariant} if there exists a vector field $\widetilde{X}$ on $\mathcal{G}^{(2)}$ such that 
\begin{enumerate}
\item $\widetilde{X} \in \ker p_{1*}$,
\item $\widetilde{X}$ is $p_2$-related to $X$, and
\item $\widetilde{X}$ is $m$-related to $X$.
\end{enumerate}
The space of left-invariant vector fields is denoted by $\livf (\mathcal{G})$.
\end{dfn}

\begin{rmk}
By Lemma \ref{lemma:lift}, the first two conditions above imply that $X$ is $t$-related to $0$, or in other words, that $X$ is tangent to the $t$-fibres.

\end{rmk}
\begin{rmk}Right-invariant vector fields are defined similarly, with the roles of $p_1$ and $p_2$ reversed.
\end{rmk}

\begin{prop} 
$\livf(\mathcal{G})$ is a $C^\infty(\mathcal{M})$-module.
\end{prop}
\begin{proof}
Let $X$ be a left-invariant vector field.  Then for any $f \in C^\infty(\mathcal{M})$, $s^*f \cdot X$ is also left-invariant; this follows from the fact that $s \circ m = s \circ p_2$.  Let $\widetilde{s^*f \cdot X} = (s \circ m)^*f \cdot \widetilde{X}$, which satisfies the conditions of Definition \ref{dfn-left}.
\end{proof}

\begin{dfn} The vector bundle $\mathcal{A} = \mathcal{A}(\mathcal{G}) \to \mathcal{M}$ is $e^*(T_t \mathcal{G})$, where $T_t \mathcal{G}$ is the subbundle of $T\mathcal{G}$ consisting of vectors tangent to the $t$-fibres.
\end{dfn}

\begin{rmk}It will be more practical to consider the module $\Gamma(\mathcal{A})$.  The sections of $\mathcal{A}$ are linear maps $\xi: C^\infty(\mathcal{G}) \to C^\infty(\mathcal{M})$ satisfying the ``$e$-derivation'' property 
\[\xi(fg) = \xi(f) \cdot e^*g + (-1)^{|\xi||f|}e^*f \cdot \xi(g)\] 
and the ``$t$-tangent'' property 
\[\xi \circ t^* = 0.\]
These are $t$-tangent vector fields along $\mathcal{M}$ in $\mathcal{G}$, in the sense of Remark \ref{rmk:along}.
\end{rmk}

We will next give some results that will be used to prove, in Theorem \ref{thm:lefte}, that $\livf(\mathcal{G})$ and $\Gamma(\mathcal{A})$ are isomorphic as $C^\infty(\mathcal{M})$-modules.

\begin{prop}There is a module homomorphism from $\livf(\mathcal{G})$ to $\Gamma(\mathcal{A})$.
\end{prop}
\begin{proof}
For any left-invariant vector field $X$, $\xi_X \defequal e^* \circ X$ is a section of $\mathcal{A}$.  Since $s \circ e = id$, the map $X \mapsto \xi_X$ respects the module structures.
\end{proof}

Let $\xi \in \Gamma(\mathcal{A})$.  Since $\xi \circ t^* = 0$, there exists by Lemma \ref{lemma:lift} a unique vector field $\xi_2$ along $\mathcal{G}$ in $\mathcal{G}^{(2)}$, where $\mathcal{G}$ is embedded via the degeneracy map $\Delta^1_1$, satisfying
\begin{align}
\xi_2 \circ p_1^* &= 0, \label{eqn:lift1}\\
\xi_2 \circ p_2^* &= s^* \circ \xi.
\end{align}

\begin{lemma}\label{lemma:erel}
$\xi$ and $\xi_2$ are related in the following way:
\[e^* \circ \xi_2 = \xi \circ \Delta_1^{1*}.\]
\end{lemma}

\begin{proof}
It is immediate from the definition of $\xi_2$ that
\begin{align}
e^* \circ \xi_2 \circ p_1^* &= 0, \\
e^* \circ \xi_2 \circ p_2^* &= e^* \circ s^* \circ \xi = \xi.
\end{align}
From the identities $p_1 \circ \Delta_1^1 = e \circ t$ and $p_2 \circ \Delta_1^1 = id$ it follows that
\begin{align}
\xi \circ \Delta_1^{1*} \circ p_1^* &= \xi \circ t^* \circ e^* = 0, \\
\xi \circ \Delta_1^{1*} \circ p_2^* &= \xi.
\end{align}
Then $e^* \circ \xi_2 = \xi \circ \Delta_1^{1*}$ by Lemma \ref{lemma:equal}.
\end{proof}

\begin{lemma}\label{lemma:leftlong}
The operator $X_{\xi} \defequal \xi_2 \circ m^*$ is a left-invariant vector field on $\mathcal{G}$.
\end{lemma}

\begin{proof}
To show that $X_{\xi}$ is left-invariant, a vector field $\widetilde{X}_\xi$ on $\mathcal{G}^{(2)}$ will be constructed such that the conditions of Definition \ref{dfn-left} are satisfied.

Consider $\mathcal{G}^{(3)}$ as the fibre product over $\mathcal{G}$ of two copies of $\mathcal{G}^{(2)}$.  There is the following commutative diagram:
\[\xymatrix{
\mathcal{G}^{(3)} \ar[r]^{\sigma_3^3} \ar[d]_{\sigma_0^3} & \mathcal{G}^{(2)} \ar[d]^{\sigma_0^2} \\
\mathcal{G}^{(2)} \ar[r]^{\sigma_2^2} & \mathcal{G}.
}\]
As a result of (\ref{eqn:lift1}), there exists by Lemma \ref{lemma:lift} the vector field $\xi_3$ along $\mathcal{G}^{(2)}$ in $\mathcal{G}^{(3)}$, where $\mathcal{G}^{(2)}$ is embedded via $\Delta^2_2$, by the properties 
\begin{align}
\xi_3 \circ \sigma_3^{3*} &= 0 \label{eqn:lift2a}\\
\xi_3 \circ \sigma_0^{3*} &= \sigma_0^{2*} \circ \xi_2. \label{eqn:lift2b}
\end{align}
From (\ref{eqn:lift2a}), (\ref{eqn:lift2b}) and the identities (\ref{eqn:twoface}),
\begin{equation} 
\begin{split}\label{eqn:lift2_1}
\xi_3 \circ \sigma_1^{3*} \circ \sigma_0^{2*} &= \xi_3 \circ \sigma_0^{3*} \circ \sigma_0^{2*} = \sigma_0^{2*} \circ \xi_2 \circ \sigma_0^{2*} \\
&= \sigma_0^{2*} \circ \sigma_0^{1*} \circ \xi = \sigma_1^{2*} \circ \sigma_0^{1*} \circ \xi \\
&= \sigma_1^{2*} \circ \xi_2 \circ \sigma_0^{2*},
\end{split}
\end{equation}
and
\begin{equation}\label{eqn:lift2_2}
\xi_3 \circ \sigma_1^{3*} \circ \sigma_2^{2*} = \xi_3 \circ \sigma_3^{3*} \circ \sigma_1^{2*} = 0.
\end{equation}
By Lemma \ref{lemma:equal}, it follows from (\ref{eqn:lift2_1}) and (\ref{eqn:lift2_2}) that
\begin{equation}\label{eqn:lift2c}
\xi_3 \circ \sigma_1^{3*} = \sigma_1^{2*} \circ \xi_2.
\end{equation}

The operator $\widetilde{X}_\xi \defequal \xi_3 \circ \sigma_2^{3*}$ is a vector field on $\mathcal{G}^{(2)}$.  By (\ref{eqn:lift2a}), (\ref{eqn:lift2b}), (\ref{eqn:lift2c}), and the identities (\ref{eqn:twoface}), the left-invariance conditions
\begin{align}
\widetilde{X}_\xi \circ \sigma_2^{2*} &= 0, \\
\widetilde{X}_\xi \circ \sigma_1^{2*} &= \sigma_1^{2*} \circ X_\xi, \\
\widetilde{X}_\xi \circ \sigma_0^{2*} &= \sigma_0^{2*} \circ X_\xi,
\end{align}
are satisfied.
\end{proof}

\begin{thm}\label{thm:lefte}
$\livf(\mathcal{G})$ is naturally isomorphic as a $C^\infty(\mathcal{M})$-module to $\Gamma(\mathcal{A})$.
\end{thm}

\begin{proof}
It remains to check that the maps $X \mapsto \xi_X$ and $\xi \mapsto X_\xi$ are inverses of each other.  Let $X$ be a left-invariant vector field and $\xi_X = \Delta_0^{0*} \circ X$ its image in $\Gamma(\mathcal{A})$.  Since
\begin{equation*}
\Delta_1^{1*} \circ \widetilde{X} \circ \sigma_2^{2*} = 0
\end{equation*}
and
\begin{equation*}
\Delta_1^{1*} \circ \widetilde{X} \circ \sigma_0^{2*} = \Delta_1^{1*} \circ \sigma_0^{2*} \circ X = \sigma_0^{1*} \circ \Delta_0^{0*} \circ X = \sigma_0^{1*} \circ \xi_X,
\end{equation*}
the lift $(\xi_X)_2$ is equal to $\Delta_1^{1*} \circ \widetilde{X}$ by Lemma \ref{lemma:equal}.  It follows that
\begin{equation*}
\begin{split}
X_{\xi_X} &= (\xi_X)_2 \circ \sigma_1^{2*} \\
&= \Delta_1^{1*} \circ \widetilde{X} \circ \sigma_1^{2*} \\
&= \Delta_1^{1*} \circ \sigma_1^{2*} \circ X \\
&= X.
\end{split}
\end{equation*}

In the other direction, let $\xi \in \Gamma(\mathcal{A})$ and form the left-invariant vector field $X_\xi = \xi_2 \circ \sigma_1^{2*}$.  From Lemma \ref{lemma:erel} it follows that
\begin{equation*}
\xi_{X_\xi} = \Delta_0^{0*} \circ \xi_2 \circ \sigma_1^{2*} = \xi \circ \Delta_1^{1*} \circ \sigma_1^{2*} = \xi.\qedhere
\end{equation*}
\end{proof}

\begin{thm}\label{thm:span} If $\mathcal{G}$ is a groupoid, then the pullback bundle $s^* \mathcal{A}$ is isomorphic to $T_t \mathcal{G}$.
\end{thm}

\begin{proof}
Since $\mathcal{A} = e^* (T_t \mathcal{G})$, a section of $s^* \mathcal{A}$ is a $t$-tangent vector field ``along $e \circ s(\mathcal{G})$'' in $\mathcal{G}$.  A typical element of $\Gamma (s^* \mathcal{A})$ is then a sum of operators on $C^\infty(\mathcal{G})$ of the form $f \cdot s^*\circ \xi$, where $f \in C^\infty(\mathcal{G})$ and $\xi \in \Gamma(\mathcal{A})$.  There is a natural map $\Gamma (s^* \mathcal{A}) \to \vect_t(\mathcal{G}) \defequal \Gamma(T_t \mathcal{G})$ defined by
\begin{equation} \label{eqn:span1} f \cdot s^*\circ \xi \mapsto f X_\xi, \end{equation}
where $X_\xi$ is the left-invariant vector field associated to $\xi$, as defined in Lemma \ref{lemma:leftlong}.  The map (\ref{eqn:span1}) is clearly a $C^\infty(\mathcal{G})$-module homomorphism and thus describes a bundle map $s^* \mathcal{A} \to T_t \mathcal{G}$.

Our next step will be to define an inverse map $T_t \mathcal{G} \to s^*\mathcal{A}$.  Let $Y \in \vect_t (\mathcal{G})$.  Since $Y$ is $t$-tangent, there exists by Lemma \ref{lemma:lift} a unique vector field $Y^{(2)}$ on $\mathcal{G}^{(2)}$ such that
\begin{align*}
Y^{(2)} \circ p_1^* &= 0, \\
Y^{(2)} \circ p_2^* &= p_2^* \circ Y.
\end{align*}
It is immediate from the definition that, for any $f \in C^\infty(\mathcal{G})$, 
\begin{equation}\label{eqn:fytwo}
(fY)^{(2)} = (p_2^*f) Y^{(2)}.
\end{equation}

Now let $\widehat{Y}$ be the operator on $C^\infty(\mathcal{G})$ defined by
\begin{equation*}
\widehat{Y} = i_L^* \circ Y^{(2)} \circ m^*.
\end{equation*}
This operator is not a derivation, so $\widehat{Y}$ is not a vector field on $\mathcal{G}$.  Since $m \circ i_L = e \circ s$, we can see that $\widehat{Y}$ is in fact a vector field along $e \circ s(\mathcal{G})$ in $\mathcal{G}$.  Furthermore, from the identity $t \circ m = t \circ p_1$, it follows that $\widehat{Y} \circ t^* = 0$.  Thus, $\widehat{Y}$ may be viewed as a section of $(e \circ s)^* T_t G = s^* (\mathcal{A})$.

As a consequence of (\ref{eqn:fytwo}) and the identity $p_2 \circ i_L = id$, the map $Y \mapsto \widehat{Y}$ is a $C^\infty(\mathcal{G})$-module homomorphism and thus describes a bundle map $T_t \mathcal{G} \to s^* \mathcal{A}$.

Let us now compose the map (\ref{eqn:span1}) with the map $Y \mapsto \widehat{Y}$.  In one direction, the composition takes an element of the form $f \cdot s^*\circ \xi$, where $f \in C^\infty(\mathcal{G})$ and $\xi \in \Gamma(\mathcal{A})$, to $\widehat{f X_\xi}$.  Using (\ref{eqn:fytwo}), the identities $p_2 \circ i_L = id$ and $m \circ i_L = e \circ s$, and the fact that $X_\xi$ is left-invariant, we compute
\begin{equation*}
\begin{split}
\widehat{f X_\xi} &= i_L^* \circ (f X_\xi)^{(2)} \circ m^* \\
&= i_L^* \circ \left[ (p_2^*f) X_\xi^{(2)} \right] \circ m^* \\
&= f \cdot i_L^* \circ X_\xi^{(2)} \circ m^* \\
&= f \cdot i_L^* m^* \circ X_\xi \\
&= f \cdot s^* e^* X_\xi \\
&= f \cdot s^* \circ \xi,
\end{split}
\end{equation*}
thereby proving that the composition of maps $s^* \mathcal{A} \to T_t \mathcal{G} \to s^* \mathcal{A}$
is equal to the identity map.  Since the bundles $s^* \mathcal{A}$ and $T_t \mathcal{G}$ are of the same rank, we conclude that the maps we have given are isomorphisms of vector bundles.
\end{proof}

\begin{cor}\label{cor:span}
The left-invariant vector fields span the space of $t$-tangent vector fields on $\mathcal{G}$.
\end{cor}

\begin{rmk} The proof of Theorem \ref{thm:span} uses the inverse map and therefore the conclusion does not hold for Lie categories.  As an example, consider the monoid $(\reals, \times)$ of real numbers under multiplication.  The left-invariant vector fields on $(\reals, \times)$ are of the form $\lambda x \pdiff{}{x}$, where $\lambda$ is a constant.  As all the left-invariant vector fields vanish at the noninvertible point $x=0$, they do not span the space of vector fields on $\reals$.
\end{rmk}

\begin{prop}
\begin{enumerate}
\item The Lie bracket of two left-invariant vector fields is again left-invariant.
\item If $X$ is a left-invariant vector field, then there is a vector field $\rho(X)$ on $\mathcal{M}$ that is $s$-related to $X$.
\end{enumerate}
\end{prop}

\begin{proof}
Let $X$ and $Y$ be left-invariant vector fields with $\widetilde{X}$ and $\widetilde{Y}$ the corresponding vector fields on $\mathcal{G}^{(2)}$.  It is straightforward to see that $[\widetilde{X}, \widetilde{Y}]$ satisfies the conditions of Definition \ref{dfn-left}.

For the second part of the proposition, let $\rho(X) = e^* \circ X \circ s^* = \xi_X \circ s^*$.  Since $X = (\xi_X)_2 \circ m^*$,
\begin{equation}
X \circ s^* = (\xi_X)_2 \circ m^* \circ s^* = (\xi_X)_2 \circ p_2^* \circ s^* = s^* \circ \xi_X \circ s^*,
\end{equation}
and it follows that $X$ and $\rho(X)$ are $s$-related.
\end{proof}

\begin{thm} The bundle $\mathcal{A} \to \mathcal{M}$ has the structure of a Lie algebroid.
\end{thm}
\begin{proof}
The map $X \mapsto \rho(X)$ is a module homomorphism, and hence defines an anchor map $\rho: \mathcal{A} \to T\mathcal{M}$.  Using the fact that a left-invariant vector field $X$ is $s$-related to $\rho(X)$, the anchor identity may be verified as follows:
\begin{equation}
\begin{split}
\rho \left( [X,Y] \right) &= e^* \circ [X,Y] \circ s^* \\ 
&= e^* s^* \circ [\rho(X), \rho(Y) ] \\
&= [\rho(X), \rho(Y)].
\end{split}
\end{equation}

Furthermore, for any left-invariant vector fields $X$ and $Y$ and $f \in C^\infty(\mathcal{M})$,
\begin{align}\label{eqn-leibniz}
[X, (s^*f)\cdot Y] &= X(s^*f) \cdot Y + (-1)^{|X||f|}s^*f \cdot [X,Y] \\
&= s^*\rho(X)(f) \cdot Y + (-1)^{|X||f|}s^*f \cdot [X,Y],
\end{align}
so the Leibniz rule is satisfied.
\end{proof}

\begin{prop}
If $X$ is a left-invariant vector field and $Y$ is a right-invariant vector field, then $[X,Y] = 0$.
\end{prop}
\begin{proof}
Let $\widetilde{X}$ and $\widetilde{Y}$ be the respective natural lifts of $X$ and $Y$ to $\mathcal{G}^{(2)}$.  It follows directly from the properties of Definition \ref{dfn-left} (and the corresponding properties of right-invariant vector fields) that $[\widetilde{X}, \widetilde{Y}] = 0$.  Then $m^* \circ [X,Y] = [\widetilde{X}, \widetilde{Y}] \circ m^* = 0$.  Since $m^*$ is injective, $[X,Y] = 0$.
\end{proof}

\begin{prop}\label{prop:lefttoright}
Every left-invariant vector field $X$ has a naturally associated right-invariant vector field, which is $t$-related to $-\rho(X)$.
\end{prop}
\begin{proof}
Let $X$ be a left-invariant vector field.  Then $e^* \circ X - \rho(X) \circ e^*$ is an $s$-tangent vector field along $\mathcal{M}$ in $\mathcal{G}$.  Thus there is a unique right-invariant vector field $X^R$ such that 
\begin{equation} \label{eqn:ranchor}
e^* \circ X^R = e^* \circ X - \rho(X) \circ e^*.
\end{equation}

By an argument similar to that for left-invariant vector fields, there exists a vector field on $\mathcal{M}$ that is $t$-related to $X^R$, and it is clear from (\ref{eqn:ranchor}) that this vector field must be $-\rho(X)$.
\end{proof}

\subsection{Multiplicative vector fields}\label{sec:mult}

\begin{dfn}
A vector field $\psi$ on $\mathcal{G}$ is \emph{simplicial} if there exists a vector field $\phi$ on $\mathcal{M}$ that is $s$- and $t$-related to $\psi$.  In this case, $\phi$ is called the \emph{base field} of $\psi$.
\end{dfn}

\begin{rmk}
By Lemma \ref{lemma:lift}, an equivalent condition is the existence of a vector field $\psi^{(2)}$ on $\mathcal{G}^{(2)}$ that is $p_1$- and $p_2$-related to $\psi$.  A simplicial vector field induces a vector field $\psi^{(q)}$ on $\mathcal{G}^{(q)}$ for each $q$.
\end{rmk}

\begin{prop}\label{prop:decomp}
A simplicial vector field $\psi$ can be uniquely decomposed as the sum of a left-(or right-)invariant vector field and a vector field that is $e$-related to its base field $\phi$.
\end{prop}
\begin{proof}
Let $\psi$ be a simplicial vector field.  The operator $e^* \psi - \phi e^*$ is a map from $C^\infty(\mathcal{G})$ to $C^\infty(\mathcal{M})$ that defines a vector field along $\mathcal{M} \subset \mathcal{G}$.  Since $e^* \psi - \phi e^*$ is both $t$- and $s$-tangent, there exist unique left- and right-invariant vector fields $\ell_\psi$ and $r_\psi$, respectively, such that $e^* \ell_\psi = e^* r_\psi = e^* \psi - \phi e^*$.

It is immediate that the vector field $\psi - \ell_\psi$ is $e$-related to $\phi$.  Thus there exists the decomposition $\psi =  \ell_\psi + (\psi - \ell_\psi)$.  Similarly, $\psi = r_\psi + (\psi - r_\psi)$ is a decomposition of $\psi$ into a right-invariant vector field and a vector field that is $e$-related to $\phi$.
\end{proof}

\begin{dfn}\label{dfn:multiplicative}
A vector field $\psi$ on $\mathcal{G}$ is \emph{multiplicative} if there exists a vector field $\psi^{(2)}$ on $\mathcal{G}^{(2)}$ that is $p_1$-, $p_2$-, and $m$-related to $\psi$.
\end{dfn}

\begin{prop}\label{prop:multe}
If $\psi$ is a multiplicative vector field with base field $\phi$, then $\psi$ is $e$-related to $\phi$, and $\psi$ is $i$-related to itself.
\end{prop}

\begin{proof}
Using the decompositions of Proposition \ref{prop:decomp}, $\psi^{(2)}$ may be decomposed  as $(\psi - r_\psi) \times (\psi - \ell_\psi) + \widetilde{\ell_\psi} + \widetilde{r_\psi}$, where tilde denotes the natural lifting of left- and right-invariant vector fields to $\mathcal{G}^{(2)}$ (see Definition \ref{dfn-left}.).  The first term $(\psi - r_\psi) \times (\psi - \ell_\psi)$ is $(\Delta_0^1 \circ e)$-related to $\phi$.  Then
\begin{equation}\label{eqn:tm1}
\begin{split}
e^*\Delta_0^{1*}\psi^{(2)}m^* &= \phi e^* \Delta_0^{1*}m^* + e^*\Delta_0^{1*}m^*(\ell_\psi + r_\psi) \\
&=\phi e^* + e^*(\ell_\psi + r_\psi) \\
&=2e^* \psi - \phi e^*.
\end{split}
\end{equation}
On the other hand, if $\psi^{(2)}$ is $m$-related to $\psi$, then
\begin{equation}\label{eqn:tm2}
\begin{split}
e^*\Delta_0^{1*}\psi^{(2)}m^* &= e^*\Delta_0^{1*}m^* \psi \\
&=e^* \psi.
\end{split}
\end{equation}
It follows from (\ref{eqn:tm1}) and (\ref{eqn:tm2}) that $\psi$ and $\phi$ are $e$-related.

It remains to show that $\psi$ is $i$-related to itself.  By $e$-relatedness,
\begin{equation}
t^* \circ e^* \circ \psi = \psi \circ t^* \circ e^*,
\end{equation}
and using the inverse condition $m \circ i_R = e \circ t$, we have 
\begin{equation}
\psi \circ i_R^* \circ m^* = i_R^* \circ \psi^{(2)} \circ m^*.
\end{equation}
The relation $p_1 \circ i_R = id$ trivially holds, so we also have
\begin{equation}
\psi \circ i_R^* \circ p_1^* = i_R^* \circ \psi^{(2)} \circ p_1^*.
\end{equation}
Since $p_2 = m \circ (i \times id) \circ (p_1 \times m) \circ \Delta$, it follows that 
\begin{equation}
\psi \circ i_R^* \circ p_2^* = i_R^* \circ \psi^{(2)} \circ p_2^*,
\end{equation}
implying that $\psi$ and $\psi^{(2)}$ are $i_R$-related.  Since $i = p_2 \circ i_R$, we conclude that $\psi$ is $i$-related to itself.
\end{proof}

\begin{rmk}Multiplicative vector fields were defined by Mackenzie and Xu \cite{mac-xu} (in the ordinary case) to be vector fields that are groupoid homomorphisms $\mathcal{G} \to T\mathcal{G}$.  Proposition \ref{prop:multe}, together with Definition \ref{dfn:multiplicative}, can be interpreted as the analogous statement for supergroupoids: a multiplicative vector field of degree $j$ is a groupoid homomorphism $\mathcal{G} \to [-j]T\mathcal{G}$. 
\end{rmk}

\begin{prop} \label{prop:assmult} Let $X$ be a left-invariant vector field, and let $X^R$ be the associated right-invariant vector field (see Proposition \ref{prop:lefttoright}).  Then $\psi_X \defequal X - X^R$ is a multiplicative vector field.
\end{prop}

\begin{proof} 
It is immediate from the definition of $X^R$ that $\psi_X$ is $e$-related to $\rho(X)$.  It follows that $\psi_X^{(2)}$ is $\Delta_i^1$-related to $\psi_X$ for $i = 0,1$.  As a result,
\begin{equation}\label{am1} 
\Delta_0^{1*} \psi_X^{(2)} m^* = \psi_X \Delta_0^{1*} m^* = \psi_X = \Delta_0^{1*} m^* \psi_X,
\end{equation}
where the identity $m \circ \Delta_0^1 = id$ has been used.

On the other hand, using the definition of $\psi_X$, we may write
\begin{equation}\begin{split}\label{eqn:psi2}
\psi_X^{(2)} &= \left(X- X^R\right) \times \left(X- X^R\right)\\
 &= 0 \times X - X^R \times 0 + X \times (-X^R)\\
 &= \widetilde{X} - \widetilde{X}^R + X \times (-X^R),\end{split}\end{equation}
where $\widetilde{X}$ and $\widetilde{X}^R$ are the respective lifts of left- and right-invariant vector fields to $\mathcal{G}^{(2)}$, so that

\begin{equation}\label{am2}
\Delta_0^{1*} \psi_X^{(2)} m^* = \Delta_0^{1*}\left(\widetilde{X} - \widetilde{X}^R + X \times (-X^R)\right) m^* = \Delta_0^{1*} m^* \psi_X + \Delta_0^{1*}(X \times (-X^R)) m^*.
\end{equation}

From (\ref{am1}) and (\ref{am2}), it follows that
\begin{equation}\label{am3}
\Delta_0^{1*}(X \times (-X^R)) m^* = 0.
\end{equation}

The goal will be to show, in fact, that $(X \times (-X^R)) m^* = 0$.  In simplicial notation,
\begin{equation} \begin{split}\label{am4}
(X \times (-X^R)) \sigma_1^{2*} &= \Delta_1^{2*} \sigma_1^{3*} (X \times (-X^R)) \sigma_1^{2*} \\
&= \Delta_1^{2*} (\widetilde{X} \times (-X^R)) \sigma_1^{3*} \sigma_1^{2*} \\
&=  \Delta_1^{2*} (\widetilde{X} \times (-X^R)) \sigma_2^{3*} \sigma_1^{2*}.
\end{split}\end{equation}

We claim that (\ref{am4}) vanishes since $\Delta_1^{2*} (\widetilde{X} \times (-X^R)) \sigma_2^{3*}=0$.  Indeed,
\begin{equation} \begin{split}\label{am5}
\Delta_1^{2*} (\widetilde{X} \times (-X^R)) \sigma_2^{3*} \sigma_2^{2*} 
&= \Delta_1^{2*} (\widetilde{X} \times (-X^R)) \sigma_3^{3*} \sigma_2^{2*} \\
&= \Delta_1^{2*} \sigma_3^{3*} \widetilde{X} \sigma_2^{2*},
\end{split}\end{equation}
which equals $0$ since $\widetilde{X} \in \ker \sigma_{2*}^2$, and
\begin{equation} \begin{split}\label{am6}
\Delta_1^{2*} (\widetilde{X} \times (-X^R)) \sigma_2^{3*} \sigma_0^{2*} 
&= \Delta_1^{2*} (\widetilde{X} \times (-X^R)) \sigma_0^{3*} \sigma_1^{2*} \\
&= \Delta_1^{2*} \sigma_0^{3*} (X \times (-X^R)) \sigma_1^{2*} \\
&= \sigma_0^{2*} \Delta_0^{1*} (X \times (-X^R)) \sigma_1^{2*},
\end{split}\end{equation}
which vanishes by (\ref{am3}).
It has thus been shown that
\begin{equation}\label{eqn:mtangent}
(X \times (-X^R)) m^* = 0.\end{equation}

From (\ref{eqn:psi2}) and (\ref{eqn:mtangent}), it follows that $\psi_X^{(2)}$ is $m$-related to $\psi_X$, and therefore that $\psi_X$ is multiplicative.
\end{proof}

\begin{rmk}Multiplicative vector fields, in the ordinary situation, are infinitesimal automorphisms \cite{mac-xu}.  The multiplicative vector fields of the form $\psi_X$ are the infinitesimal inner automorphisms.
\end{rmk}

The Lie bracket of two multiplicative vector fields is again multiplicative.  Furthermore, if $\psi$ is a multiplicative vector field and $X$ is a left-invariant vector field, then $[\psi, X]$ is left-invariant.  The action of $\psi$ on the left-invariant vector fields respects the algebroid structure in the following ways.

\begin{prop}\label{prop:multaction}
Let $\psi$ be a multiplicative vector field with base field $\phi$.  Then
\begin{enumerate}
\item $[\psi, \cdot ]$ is a derivation of the bracket structure;
\item $\rho([\psi, X]) = [\phi, \rho(X)].$
\end{enumerate}
\end{prop}

\begin{proof}
The first statement follows immediately from the Jacobi identity.  For the second statement, using the fact that $\psi$ and $\phi$ are $e$- and $s$-related, we have
\begin{equation}
\begin{split}
\rho([\psi, X]) &= e^* [\psi, X] s^* \\
&= e^* \left(\psi X - (-1)^{|\psi||X|}X\psi \right) s^* \\
&= \phi e^* X s^* - (-1)^{|\psi||X|} e^* X s^* \phi \\
&= [\phi, \rho(X)].\qedhere
\end{split}
\end{equation}
\end{proof}

By Proposition \ref{prop:morphderiv}, we have
\begin{cor}\label{cor:multaction}
A multiplicative vector field $\psi$ on $\mathcal{G}$ induces a morphic vector field on $\mathcal{A}(\mathcal{G})$.
\end{cor}

\subsection{Left-invariant differential forms}

The standard approach to left-invariant $1$-forms on an (ordinary) groupoid $G$ is to first consider the quotient space of $\Omega^1 (G)$ dual to the subspace $\vect_t (G)$ of $\vect(G)$ and then look for those equivalence classes satisfying a left-invariance condition.  In the simplicial approach, just as we did not need to restrict ourselves \emph{a priori} to $\vect_t (G)$ in defining left-invariant vector fields, we are able to state a left-invariance condition for differential forms that does not require us to first pass to a quotient space.  The pairing of left-invariant $1$-forms with left-invariant vector fields then has a kernel, which one must divide by in order to get the dual space $\Gamma(A^*)$.  

Let $\mathcal{G} \arrows \mathcal{M}$ be a supergroupoid.

\begin{dfn}\label{dfn:leftforms}Let $J$ be the ideal of $\Omega(\mathcal{G}^{(2)})$ generated by $p_1^* \left(\Omega^1(\mathcal{G}) \right)$.  A differential form $\alpha \in \Omega (\mathcal{G})$ is \emph{left-invariant}\footnote{One could also consider left-invariant pseudoforms in $C^\infty([-1]T\mathcal{G})$.  All the material in this section extends in a straightforward manner to this case.} if $m^* \alpha \equiv p_2^* \alpha \pmod J$.
\end{dfn}

\begin{prop} A function $f \in C^\infty(\mathcal{G})$ is left-invariant if and only if $f \in \im s^*$.
\end{prop}
\begin{proof} Since the ideal $J$ does not contain any $0$-forms, the left-invariance condition for a function $f$ is simply $m^* f = p_2^* f$.  If $f \in \im s^*$, then the equation clearly holds, since $s \circ m = s \circ p_2$.

Suppose, conversely, that $m^* f = p_2^* f$.  Composing with $\Delta_1^1$ on the left yields the equation $f = s^* e^* f$.
\end{proof}

\begin{prop}Let $I$ be the ideal of $\Omega(\mathcal{G})$ generated by $t^*\left(\Omega^1 (\mathcal{M})\right)$.  Every element of $I$ is left-invariant.
\end{prop}
\begin{proof}
It suffices to check that any element of the form $\beta \wedge t^* \omega$, where $\omega \in \Omega^1(\mathcal{M})$, is left-invariant.  Indeed, $m^* (\beta \wedge t^* \omega) = m^* \beta \wedge m^* t^* \omega = m^* \beta \wedge p_1^* t^* \omega$, which is in $J$, and $p_2^* (\beta \wedge t^* \omega) = p_2^* \beta \wedge p_2^* t^* \omega = p_2^* \beta \wedge p_1^* s^* \omega$, which is also in $J$.  Thus $m^* (\beta \wedge t^* \omega) \equiv p_2^* (\beta \wedge t^* \omega) \equiv 0 \pmod J$.
\end{proof}

\begin{rmk}The $1$-forms in $I$ are characterized by the property that for any $X \in \vect_t (\mathcal{G})$, $\langle X, \alpha \rangle = 0$.
\end{rmk}

\begin{dfn} The algebra $\Omega_{LI}(\mathcal{G})$ is defined to be the subalgebra of $\Omega(\mathcal{G})/I$ consisting of equivalence classes of left-invariant forms.
\end{dfn}

\begin{rmk} Since $\Omega(\mathcal{G}) = \bigwedge \Omega^1(\mathcal{G})$ and the ideals $I$ and $J$ are both generated by $1$-forms, it follows that $\Omega_{LI}(\mathcal{G}) = \bigwedge \Omega^1_{LI}(\mathcal{G})$.
\end{rmk}

\begin{prop}\label{prop:leftpair}Let $\alpha \in \Omega^1_{LI}(\mathcal{G})$.  Then for any $X \in \livf(\mathcal{G})$, the pairing $\langle X, \alpha \rangle$ yields a left-invariant function.
\end{prop}
\begin{proof}
Since $\widetilde{X}$ is $p_1$-tangent, it annihilates any $1$-form in $J$.  So
\begin{align*} m^*\langle X, \alpha \rangle &= \langle \widetilde{X}, m^* \alpha \rangle \\
&= \langle \widetilde{X}, p_2^* \alpha \rangle \\
&= p_2^* \langle X, \alpha \rangle.\qedhere
\end{align*}
\end{proof}

\begin{prop}\label{prop:leftdual}
There is a natural isomorphism from $\Omega^1_{LI}(\mathcal{G})$ to $\Gamma(\mathcal{A}^*)$.  
\end{prop}
\begin{proof}
It follows from Proposition \ref{prop:leftpair} that the pairing of left-invariant $1$-forms and left-invariant vector fields satisfies $\langle X, \alpha \rangle = s^*e^*\langle X, \alpha \rangle$.  Thus there is a map $\alpha \mapsto \alpha_0$ from $\Omega^1_{LI}(\mathcal{G})$ to $\Gamma(\mathcal{A}^*)$, defined by the equation
\begin{equation} \langle \xi, \alpha_0 \rangle = e^* \langle X_\xi, \alpha \rangle \end{equation}
for $X \in \Gamma(\mathcal{A})$.

To describe a map in the other direction, suppose that $\alpha_0 \in \Gamma(\mathcal{A}^*)$.  Since the left-invariant vector fields span the space of $t$-tangent vector fields, the equation 
\begin{equation}\label{eqn:leftformlift}\langle fX, \alpha \rangle = f \cdot s^* \langle \xi_X, \alpha_0 \rangle,\end{equation}
for $f \in C^\infty(\mathcal{G})$, $X \in \livf(\mathcal{G})$, determines, uniquely up to an element of $I$, a $1$-form $\alpha$ on $\mathcal{G}$.  Observe that, for any left-invariant vector field $X$, $\langle \widetilde{X}, m^*\alpha \rangle = m^*\langle X, \alpha \rangle = m^* s^* \langle \xi_X, \alpha_0 \rangle = p_2^* s^* \langle \xi_X, \alpha_0 \rangle = p_2^*\langle X, \alpha \rangle = \langle \widetilde{X}, p_2^*\alpha \rangle$.  Since the lifts of left-invariant vector fields span the space of $p_1$-tangent vector fields on $\mathcal{G}^{(2)}$, it follows that $\alpha$ is left-invariant.

It is straightforward to check the maps in both directions are inverses of each other.
\end{proof}

\begin{cor} \label{cor:liforms}
$\Omega_{LI}(\mathcal{G})$ is naturally isomorphic to $\bigwedge \Gamma(\mathcal{A}^*)$.  \end{cor}

\subsection{Cartan calculus on left-invariant forms}\label{sec:cartleft}

The Cartan calculus consists of a Lie superalgebra of operators that act as derivations on the algebra of differential forms on a manifold.  On a Lie groupoid, we would like to describe the Cartan calculus in the context of left-invariant forms.  An operator $D$ on $\Omega(\mathcal{G})$ induces an action on $\Omega_{LI}(\mathcal{G})$ if $D(I) \subseteq I$ (so that $D$ descends to the quotient $\Omega(\mathcal{G}) / I$) and, for any left-invariant form $\alpha$, $D\alpha$ is left-invariant.

\begin{prop}
The de Rham differential induces an action on $\Omega_{LI}(\mathcal{G})$.
\end{prop}
\begin{proof}
This follows directly from the fact that $d$ commutes with pullback maps. 
\end{proof}

\begin{prop}
The contraction operator $\iota_X$ induces an action on $\Omega_{LI}(\mathcal{G})$ if and only if $X$ is left-invariant.
\end{prop}
\begin{proof}
The ``if'' part is essentially the content of Proposition \ref{prop:leftpair}.  The ``only if'' part will now be shown.

It is sufficient to consider the action of $\iota_X$ on $0$- and $1$-forms.  Recall that $\iota_X$ annihilates $\Omega^0 (\mathcal{G})$ and maps $\Omega^1 (\mathcal{G})$ to $\Omega^0 (\mathcal{G})$.  Since there are no nontrivial $0$-forms in $I$, $\iota_X (I) \subseteq I$ if and only if $\iota_X t^*\omega = 0$ for any $\omega \in \Omega^1(\mathcal{M})$.  Thus $X$ is $t$-tangent or, equivalently, there exists a vector field $\widetilde{X}$ on $\mathcal{G}^{(2)}$ that is $p_1$-related to $0$ and $p_2$-related to $X$.  

Next suppose that $\alpha$ is a left-invariant $1$-form.  Since $\iota_X \alpha$ is a left-invariant function, there must exist an $f \in C^\infty(\mathcal{M})$ such that $\iota_X \alpha = s^*f$.  Then $\iota_{\widetilde{X}} m^* \alpha = \iota_{\widetilde{X}} p_2^* \alpha = p_2^* \iota_X \alpha = p_2^* s^* f = m^* s^* f = m^* \iota_X \alpha$.  By the dual statement to Corollary \ref{cor:span}, it follows that $X$ and $\widetilde{X}$ are $m$-related and therefore that $X$ is left-invariant.
\end{proof}

\begin{prop}\label{prop:lieaff}
The Lie derivative $\lie_X$ induces an action on $\Omega_{LI}(\mathcal{G})$ if and only if $X$ is the sum of a left-invariant vector field and a multiplicative vector field.
\end{prop}
\begin{proof}
Again, it is sufficient to consider the action on $0$- and $1$-forms.  Since left-invariant $0$-forms are of the form $s^*f$, the space of left-invariant $0$-forms is invariant under $\lie_X$ if and only if $X$ is $s$-related to some vector field $\phi_s$ on $\mathcal{M}$.

Next, consider the requirement that $I$ be invariant under $\lie_X$.  Restricting to $1$-forms of the form $d t^* f$, where $f \in C^\infty(\mathcal{M})$, we have $\lie_X d t^* f = \pm d \lie_X t^*f$, which is only in $I$ if $\lie_X t^*f \in \im t^*$.  This implies the existence of a vector field $\phi_t$ on $\mathcal{M}$ that is $t$-related to $X$.

The operator $e^*X - \phi_t e^*$ is a $t$-tangent vector field along $\mathcal{M}$ in $\mathcal{G}$.  Thus there exists a left-invariant vector field  $\ell_X$ such that $e^* \ell_X = e^*X - \phi_t e^*$.  Clearly, $\rho(\ell_X) = \phi_s - \phi_t$.  Let $\psi$ be defined as $X - \ell_X$.  Since $\psi$ is $s$-related to $\phi_t$, it is possible to form a vector field $\hat{X}$ on $\mathcal{G}^{(2)}$ that is $p_1$-related to $\psi$ and $p_2$-related to $X$.  Observe that the ideal $J$ is invariant under $\lie_{\hat{X}}$.

Using the final requirement that $\lie_X \alpha$ be left-invariant whenever $\alpha$ is left-invariant, we see that, mod $J$, $m^* \lie_X \alpha \equiv p_2^* \lie_X \alpha = \lie_{\hat{X}} p_2^* \alpha \equiv \lie_{\hat{X}} m^* \alpha$.  By the dual statement to Corollary \ref{cor:span}, it follows that $X$ and $\hat{X}$ are $m$-related.  Since $\psi^{(2)} = \hat{X} - \widetilde{\ell_X}$, it is clear that $\psi$ and $\psi^{(2)}$ are $m$-related.
\end{proof}

\begin{rmk} By making the obvious modifications to the proof of Proposition \ref{prop:lieaff}, it may similarly be shown that $\lie_X$ induces an action on $\Omega_{LI}(\mathcal{G})$ if and only if $X$ is the sum of a right-invariant vector field and a multiplicative vector field.  Because right-invariant vector fields annihilate left-invariant forms, the action of $\lie_X$, when $X$ is a left-invariant vector field, is the same as the action of $\lie_{\psi_X}$, where $\psi_X$ is the associated multiplicative vector field (see Proposition \ref{prop:assmult}).
\end{rmk}

We summarize the results of this section in the following theorem:

\begin{thm}
\begin{enumerate}
\item The isomorphism between $\Omega_{LI}(\mathcal{G})$ and $\bigwedge \Gamma(\mathcal{A}^*)$ is an isomorphism of differential graded algebras and is equivariant with respect to the operators $\iota_X$ and $\lie_X$ for all $X \in \Gamma(\mathcal{A})$.
\item If $\psi$ is a multiplicative vector field on $\mathcal{G}$, then the operator $\lie_\psi$ induces a morphic vector field on $\mathcal{A}$.
\end{enumerate}
\end{thm}

    \chapter{Double structures and equivariant cohomology}
    \label{chapter:double}
An important algebraic object associated to an (ordinary) Lie groupoid $G$ is the de Rham double complex \cite{lx} of $G$, whose cochains are differential forms on $G^{(q)}$.  The de Rham double complex is a model for the classifying space $BG$, which classifies homotopy classes of principal $G$-bundles \cite{buffet}.  The following examples are of particular interest:
\begin{itemize}
\item If $G$ is a group, then the de Rham double complex is equal to the Bott-Shulman complex \cite{bott} for the cohomology of $BG$.
\item In the case of a pair groupoid $G = M \times M$, the de Rham double complex extends Kock's theory of combinatorial differential forms \cite{bm} \cite{kock} to include symmetric tensors.
\item In the case of an action groupoid, the de Rham double complex is equal to the simplicial model of equivariant cohomology \cite{am:cw}.
\end{itemize}

In this chapter, we demonstrate that the de Rham double complex of a groupoid $G$ may be viewed as a ``twisted'' groupoid cohomology of the supergroupoid $[-1]TG$.  We introduce the more general concept of a \emph{$Q$-groupoid}, which is a groupoid equipped with a homological multiplicative vector field.  Every $Q$-groupoid has an associated double complex.  

The Lie functor may be applied to a $Q$-groupoid to produce a \emph{$Q$-algebroid}, or an algebroid equipped with a homological morphic vector field.  In particular, the $Q$-groupoid $[-1]TG$ induces the $Q$-algebroid $[-1]TA$ (see \S \ref{sec:oneta}).  The algebroid cohomology of $[-1]TA$ may also be ``twisted'' to produce a double complex.  This double complex is in fact a model for $EG$, and it is necessary to restrict to a basic subcomplex in order to obtain interesting cohomology.  We describe the basic subcomplex, the definition of which requires a choice of a connection on $A$.  In the case of an action algebroid $A = M \times \mathfrak{g}$, there is a natural choice of connection, and the double complex and resulting basic subcomplex may be identified with those of the BRST-Cartan model of equivariant cohomology \cite{kalkman}.  Via a generalization of the Mathai-Quillen-Kalkman isomorphism \cite{kalkman} \cite{mq}, we also obtain a Weil model for algebroids.  

The generalizations described in this chapter of the BRST-Cartan and Weil models are new.  These are ``infinitesimal'' models, in the sense that they are defined without reference to an overlying groupoid.  Since in general an algebroid may not arise from a Lie groupoid, and even if one does it may be difficult to describe the groupoid explicitly, such infinitesimal models have both a theoretical and a practical value.

\section{$\mathcal{LA}$-groupoids and $Q$-groupoids}

One of the double structures described by Mackenzie \cite{mac:dbl} is that of an $\mathcal{LA}$-groupoid, or a groupoid in the category of algebroids.  In this section, it is shown that the $[-1]$ functor of Theorem \ref{thm:minusone}, when applied to an $\mathcal{LA}$-groupoid, yields a $Q$-groupoid.  The corresponding homological vector field, when introduced into the groupoid cochain complex, produces a double complex.

A simple class of $\mathcal{LA}$-groupoids consists of those of the form $TG \arrows TM$, where $G \arrows M$ is a groupoid.  In \S\ref{sub:onetg} we describe $Q$-groupoids of the form $[-1]TG \arrows [-1]TM$ and identify the associated double complex with the de Rham double complex of the groupoid $G$.  In \S\ref{sub:lieonetg} we apply the Lie functor to the groupoid $[-1]TG$ and obtain the Lie algebroid $[-1]TA$ that was described in \S\ref{sec:oneta}.

\subsection{Groupoid cohomology}
\label{section:cohomology}
The $q$-cochains for groupoid cohomology are smooth functions on $G^{(q)}$.  If $\mathcal{G}$ is a supergroupoid, then the algebra of functions $C^\infty(\mathcal{G}^{(q)})$ is itself a graded algebra, and the space of cochains therefore has a double grading.  Denote by $C^{p,q}=C^{p,q}(\mathcal{G})$ the space of degree $p$ functions on $\mathcal{G}^{(q)}$.  If $f \in C^{p,q}$ and $g \in C^{p^\prime, q^\prime}$, then the product function $(-1)^{q p'} f \times g \in C^\infty(\mathcal{G}^{q} \times \mathcal{G}^{q'})$ restricts to a degree $p + p^\prime$ function, denoted $f*g$, on $\mathcal{G}^{(q + q^\prime)}$.  The product $*$ thus gives the space of cochains a ring structure that respects the double grading.

\begin{dfn}The coboundary operator $\delta^q : C^{\bullet,q-1} \to C^{\bullet,q}$ is
\begin{equation}
\delta^q = \sum_{i=0}^q (-1)^i \sigma_i^{q*}.
\end{equation}
\end{dfn}

\begin{rmk}
The property $\delta^2 = 0$ is a consequence of (\ref{eqn:twoface}).  It also follows from the definitions of the face maps that $\delta$ is a degree $1$ derivation with respect to the $q$-grading.  
\end{rmk}

\subsection{Double complexes}

Recall that a simplicial vector field $\psi$ naturally lifts to vector fields $\psi^{(q)}$ on $\mathcal{G}^{(q)}$ for each $q$.  A simplicial vector field thus acts via the derivations $\psi^{(q)}$ on the groupoid cochains.

\begin{prop} The action of a multiplicative vector field on the groupoid cochains commutes with the coboundary operator.  
\end{prop}
\begin{proof}
Since the face maps are all built out of the three maps $p_1$, $p_2$, and $m$, it is clear that if $\psi$ is a multiplicative vector field then $\psi^{(q)}$ is $\sigma_i^q$-related to $\psi^{(q-1)}$ for all $i$ and $q$.  
\end{proof}

\begin{rmk}\label{rmk:multface}
It is also true, as a consequence of Proposition \ref{prop:multe}, that $\psi^{(q)}$ is $\Delta_i^q$-related to $\psi^{(q+1)}$ for all $i$ and $q$.
\end{rmk}

Recall that a homological vector field $\psi$ is a degree $1$ vector field satisfying $\psi^2 = 0$.

\begin{dfn}
A supergroupoid equipped with a multiplicative homological vector field is called a $Q$-groupoid.
\end{dfn}

\begin{rmk} For a $Q$-groupoid, the groupoid cohomology complex becomes a double complex $\left( C^{p,q}, \delta, \psi \right)$.  The total differential $D = \psi - (-1)^p \delta$ clearly satisfies the equation $D^2 = 0$.  Furthermore, $D$ is a derivation with respect to the total grading.
\end{rmk}

\begin{ex}\label{ex:onetg}
Consider the odd tangent prolongation groupoid $[-1]TG$.  The de Rham differential $d$ is a homological vector field on $[-1]TG$.  Since $d$ commutes with pullback maps, it is clear that $d$ is multiplicative.  This example will be discussed in more detail below.
\end{ex}

\begin{ex}\label{ex:onetstar}
If the groupoid $G$ has a Poisson bivector $\pi$, then the Poisson cohomology operator $d_\pi$ is a homological vector field on $[-1]T^*G$.  The vector field $d_\pi$ is multiplicative if and only if $G$ is a Poisson groupoid \cite{mac:dbl} \cite{mac-xu}.
\end{ex}

\begin{rmk} Both of the previous examples have the common property of being $\mathcal{VB}$-groupoids.  In fact, there are more specific structures involved, which Mackenzie has called \emph{$\mathcal{LA}$-groupoids} \cite{mac:dbl}.  
\end{rmk}

\begin{dfn}
An \emph{$\mathcal{LA}$-groupoid} is a groupoid in the category of algebroids, i.e.\ a pair of groupoids, $\Omega \arrows A$ and $G \arrows M$, such that $\Omega$ is an algebroid over $G$, $A$ is an algebroid over $M$, and all of the structure maps for the groupoid structure of $\Omega$ are algebroid morphisms over the corresponding structure maps for $G$.
\end{dfn}

\begin{thm}\label{thm:lagpd}
Let $(\Omega \arrows A, G \arrows M)$ be an $\mathcal{LA}$-groupoid. Then $[-1]_G\Omega \arrows [-1]_M A$ is a $Q$-groupoid.
\end{thm}
\begin{proof}The multiplicative vector field on $[-1]\Omega$ is the algebroid differential $d_{\Omega}$.
\end{proof}

\begin{ex} The $Q$-groupoids of Examples \ref{ex:onetg} and \ref{ex:onetstar} arise from the $\mathcal{LA}$-groupoids 
\begin{equation}
\la{TG}{G}{TM}{M}
\end{equation}
and
\begin{equation}
\la{T^*G}{G}{A^*}{M,}
\end{equation}
respectively.
\end{ex}

\subsection{$[-1]TG$}\label{sub:onetg}

Let $G \arrows M$ be a groupoid.  Because $[-1]T$ is a covariant functor, the structure maps for the groupoid structure of $G$ naturally induce a groupoid structure on $[-1]TG$ over $[-1]TM$.  In particular, the identification 
\[[-1]T\left(G \bitimes{s}{t} G\right) = \left([-1]TG\right) \bitimes{[-1]Ts}{[-1]Tt} \left([-1]TG\right)\]
is used to identify the map $[-1]Tm$ with multiplication on $[-1]TG$.

\begin{ex}
Let $G$ be a group, i.e.\ a groupoid whose base is just a point.  Since $[-1]T(\{pt.\})$ is again a point, $[-1]TG$ is a supergroup that might be called the \emph{odd tangent group} of $G$.
\end{ex}

\begin{ex}
Let $\Gamma$ be a group that acts on a manifold $M$ from the right.  Then if the $[-1]T$ functor is applied to the action groupoid $M \times \Gamma \arrows M$, the result may be identified with $[-1]TM \times [-1]TG \arrows [-1]TG$, which is the action groupoid for an induced action of $[-1]TG$ on $[-1]TM$.
\end{ex}

For each $q$, there is a natural identification of $\left([-1]TG\right)^{(q)}$ with $[-1]T\left(G^{(q)}\right)$.  Thus, when $G$ is an ordinary groupoid, the space of groupoid cochains on $[-1]TG$ is
\begin{equation}C^{p,q}([-1]TG) = \Omega^p(G^{(q)}).\end{equation}
The double complex in this case is in fact identical to the de Rham complex of the simplicial manifold $NG$ (see \S\ref{sec:drc} and Definition \ref{dfn:nerve}), whose total cohomology is equal to the cohomology of the geometric realization $|NG|$.  Segal \cite{segal} takes (in the context of topological categories) the space $|NG|$ to be the definition of the classifying space $BG$, and it has been shown by Buffet and Lor \cite{buffet} that $BG$ classifies the homotopy classes of $G$-bundles.

In the case where $G$ is a group, $|NG|$ is the classying space $BG$ in the usual topological sense.  In the case where $G = \Gamma \times M \arrows M$ is the action groupoid for a group action, then $|NG|$ is the homotopy quotient $(M \times EG)/G$ (see \S\ref{sec:equiv}).  The following propositions are immediate:

\begin{prop} If $G$ is a Lie group, then the double complex of $[-1]TG$ computes the cohomology of the classifying space $BG$.
\end{prop}

\begin{prop} \label{prop:equiv}If $G = M \times \Gamma \arrows M$ is the action groupoid for the action of a Lie group on a manifold $M$, then the double complex of $[-1]TG$ computes the equivariant cohomology $H_\Gamma^\bullet(M)$.
\end{prop}

\begin{rmk}\label{rmk:superbg}
If $\mathcal{G}$ is a supergroup, then there is not any obvious way to construct the (super\footnote{It may be possible to define the realization functor for simplicial supermanifolds such that $B\mathcal{G}$ is a superspace, in the sense of Manin \cite{manin}.})space $B\mathcal{G}$.  However, the simplicial supermanifold $N\mathcal{G}$ and the double complex of $[-1]T\mathcal{G}$ do exist and might be considered a model for ``$B\mathcal{G}$''.  Similarly, the equivariant cohomology for the action of a supergroup $\Gamma$ on a supermanifold $\mathcal{M}$ may be defined to be the total cohomology of the double complex of $[-1]T\mathcal{G}$, where $\mathcal{G} = \mathcal{M} \times \Gamma$ is the action supergroupoid.
\end{rmk}

\subsection{Applying the Lie functor to $[-1]TG$}\label{sub:lieonetg}

\begin{prop}
\begin{enumerate}
\item If $X \in \livf(G)$, then $\iota_X$ and $\lie_X$ are left-invariant vector fields on $[-1]TG$.
\item If $\psi$ is a multiplicative vector field on $G$, then $\iota_\psi$ and $\lie_\psi$ are multiplicative vector fields on $[-1]TG$.
\end{enumerate}
\end{prop}
\begin{proof}
By Proposition \ref{prop:iota}, $\iota_{\widetilde{X}}$ is $[-1]Tp_1$-related to $0$ and $[-1]Tm$- and $[-1]Tp_2$-related to $\iota_X$.  The result for $\lie_X$ follows similarly from Corollary \ref{cor:lie}.  The proof of the second statement is also similar.
\end{proof}

\begin{rmk}
By looking at local coordinates about $M$, it is easy to see that the left-invariant vector fields of the form $\iota_X$ and $\lie_X$ span (over $C^\infty([-1]TM)$) the left-invariant vector fields on $[-1]TG$.  Thus the algebroid structure for $\mathrm{Lie}([-1]TG)$ can be described in terms of left-invariant vector fields of those forms.
\end{rmk}

\begin{prop}\label{prop:lieg}
\begin{enumerate}
\item $\rho(\iota_X) = \iota_{\rho(X)}$, and $\rho(\lie_X) = \lie_{\rho(X)}$.
\item $[\iota_X, \iota_Y] = 0$, $[\lie_X, \iota_Y] = \iota_{[X,Y]}$, and $[\lie_X, \lie_Y] = \lie_{[X,Y]}$.
\end{enumerate}
\end{prop}
\begin{proof}
Recall that the anchor is defined by the fact that $X$ and $\rho(X)$ are $s$-related.  Consequently, $\iota_X$ and $\iota_{\rho(X)}$ are $[-1]Ts$-related, implying that $\iota_{\rho(X)}$ is the anchor of $\iota_X$.  The result for $\rho(\lie_X)$ holds in exactly the same manner.

The second statement is immediate from the Cartan commutation relations (see Proposition \ref{prop:cartancomm}).
\end{proof}

\begin{thm} Let $G$ be a groupoid with $A = \mathrm{Lie}(G)$.  Then $\mathrm{Lie}([-1]TG)$ is naturally isomorphic to $[-1]TA$.  The morphic vector field induced by the multiplicative vector field $d$ on $[-1]TG$ is the de Rham differential on $[-1]T([-1]A)$.
\end{thm}

\begin{proof}Comparing Proposition \ref{prop:lieg} with \S\ref{sec:oneta}, it is clear that the isomorphism is achieved by identifying, for each $X \in \Gamma(A)$, $\iota_X$ with the vertical lift $X^V$ and $\lie_X$ with the complete lift $X^C$.  

The multiplicative vector field $d$ induces an operator $D_{\widetilde{d}}$ on $\Gamma([-1]TA)$, with base vector field equal to the de Rham differential on $[-1]TM$ and determined by the equations $D_{\widetilde{d}}(X^V) = (-1)^{|X|} X^C$ and $D_{\widetilde{d}}(X^C) = 0$.  By (\ref{eqn:morph}), the corresponding morphic vector field $\widetilde{d} \in \vect([-1]([-1]TA))$ satisfies 
\begin{align*}
\left[\iota_{X^V}, \widetilde{d}\right] = -\iota_{X^C}, && \left[\iota_{X^C}, \widetilde{d}\right] = 0.
\end{align*}
  In a local coordinate system as in \S\ref{sec:oneta}, where $\{X^C_\alpha, - X^V_\alpha\}$ are the sections dual to $\{\lambda^\alpha,\dot{\lambda}^\alpha\}$, it is clear that $\widetilde{d}$ is equal to the de Rham differential $\dot{\lambda}^\alpha \pdiff{}{\lambda^\alpha} + \dot{x}^i \pdiff{}{x^i}$.
\end{proof}

\section{$Q$-algebroids}
\label{sec:inf}
The infinitesimal objects corresponding to $Q$-groupoids are $Q$-algebroids, which are defined as superalgebroids equipped with homological morphic vector fields.  As in the case of a $Q$-groupoid, we will see that every $Q$-algebroid has an associated double complex.  At first, we will focus our attention on $Q$-algebroids of the form $[-1]TA$, where $A \to M$ is an algebroid.  
In \S\ref{sec:doubleoneta}, we describe the special case where $A = M \times \mathfrak{g}$ is an action algebroid.  In this case, the double complex associated to $[-1]TA$ is equal to the complex of the BRST model of equivariant cohomology.  Many of the properties of the BRST model have been described by Kalkman \cite{kalkman}.  In particular, he has described an algebra automorphism, extending that of Mathai and Quillen \cite{mq}, that relates the BRST model to the Weil model.  In each of these models, the cohomology of the total complex is known to be $H^\bullet(M)$, and to obtain an infinitesimal model of equivariant cohomology it is necessary to pass to a ``basic'' subcomplex \cite{gs}.  Restricting to the basic subcomplex of the BRST model yields the Cartan model \cite{cartan}.

The primary goal of this section is to generalize the BRST-Cartan and Weil models to arbitrary algebroids.  In \S\ref{sec:doubleoneta}, we describe a generalization of the Mathai-Quillen-Kalkman isomorphism and use it to show that for any $A$ the cohomology of the double complex associated to $[-1]TA$ is equal to $H^\bullet(M)$.  Then, in \S\ref{sub:basic}, we propose a definition of the basic subcomplex when $A$ is equipped with a connection.  In the case of an action algebroid $A = M \times \mathfrak{g}$ with the canonical flat connection, the basic subcomplex is equal to the complex of the Cartan model.

In \S\ref{sec:geom}, we construct from any $Q$-algebroid $A$ a superalgebroid $\widetilde{A}$ whose structure incorporates both the algebroid structure and the $Q$-structure of $A$.

\subsection{$\mathcal{VB}$-algebroids}\label{sec:vba}
This section relies heavily on the results from \S\ref{section:dvb}.

\begin{dfn}
A \emph{$\mathcal{VB}$-algebroid}\footnote{c.f. Mackenzie's notion of $\mathcal{LA}$-vector bundle \cite{mac:dbl2}.} $(D,A;E,M)$ is a double vector bundle
\begin{equation}
\dvb{D}{A}{E}{M}
\end{equation}
where the vertical sides are algebroids such that the anchor maps form a vector bundle morphism
\begin{equation}\label{eqn:vbanchor}
\dvb{D}{TE}{A}{TM}
\end{equation}
and the bracket is such that
\begin{enumerate}
\item $[\Gamma_V (D, E), \Gamma_V (D, E)] = 0$,
\item $[\Gamma_{lin} (D, E), \Gamma_V (D, E)] \subseteq \Gamma_V (D, E)$, and
\item $[\Gamma_{lin} (D, E), \Gamma_{lin} (D, E)] \subseteq \Gamma_{lin} (D, E)$.
\end{enumerate}
\end{dfn}

\begin{prop}\label{prop:vbashift}
If $(D,A;E,M)$ is a $\mathcal{VB}$-algebroid, then for any $j$, there is an induced algebroid structure on $[j]_A D \to [j]E$.
\end{prop}
\begin{proof}
The anchor is obtained by applying the $[j]$ functor to the top row of (\ref{eqn:vbanchor}) and using the identification $[j]_{TM} TE = T([j]E)$ of Proposition \ref{prop:dvb2}.  The bracket arises from the isormorphism of Proposition \ref{prop:shiftsection} by extending the induced bracket on $\Gamma_{lin} ([j]_A D, [j]E) \oplus \Gamma_V ([j]_A D,[j]E)$ by the Leibniz rule.  Since it is sufficient to verify the Jacobi identity locally on a frame of sections, it follows immediately that the algebroid axioms are satisfied.
\end{proof}

\begin{ex}We have seen in \S\ref{sec:oneta} that if $A \to M$ is an algebroid, then there is an induced algebroid structure on $[-1]TA \to [-1]TM$.  It follows directly from the construction of \S\ref{sec:oneta} that
\begin{equation}
\dvb{[-1]TA}{A}{[-1]TM}{M}
\end{equation}
is a $\mathcal{VB}$-algebroid.  Proposition \ref{prop:vbashift} implies that for all $j$, there is an associated algebroid structure on $[j-1]TA \to [j-1]TM$.  In particular, by choosing $j=1$, we can recover the algebroid structure on $TA \to TM$ described in, e.g., \cite{mac-xu}.
\end{ex}

\begin{ex}Consider an $\mathcal{LA}$-groupoid
\begin{equation}
\la{\Omega}{G}{A}{M.}
\end{equation}
We may apply the Lie functor to obtain the $\mathcal{VB}$-algebroid
\begin{equation}\label{eqn:dblie}
\dvb{\Lie(\Omega)}{\Lie(G)}{A}{M.}
\end{equation}
In fact, (\ref{eqn:dblie}) has the structure of a \emph{double Lie algebroid} \cite{mac:dblie2} \cite{mac:dbl}, in that all four sides are Lie algebroids in a compatible way.  We may apply the $[-1]$ functor to the left hand side to get an algebroid $[-1]_{\Lie(G)} \Lie(\Omega) \to [-1]A$.  Furthermore, the horizontal Lie algebroid structures of (\ref{eqn:dblie}) induce homological vector fields on $[-1]_{\Lie(G)} \Lie(\Omega)$ and $[-1]A$.  We will see in \S\ref{sec:qalg} that these vector fields make $[-1]_{\Lie(G)} \Lie(\Omega) \to [-1]A$ into a $Q$-algebroid, and that this $Q$-algebroid may be equivalently obtained by first applying the $[-1]$ functor to the $\mathcal{LA}$-groupoid and then applying the Lie functor.
\end{ex}

\subsection{$Q$-algebroids}\label{sec:qalg}

Let $G \arrows M$ be a groupoid, and let $A \to M$ be $\mathrm{Lie}(G)$.  It has been shown (see \S\ref{sec:mult} and \S\ref{sec:cartleft}) that a multiplicative vector field $\psi$ on $G$ induces a morphic vector field $\widetilde{\psi}$ on $A$.  Furthermore, if $\psi$ is a homological vector field, then so is $\widetilde{\psi}$.

\begin{dfn} A \emph{$Q$-algebroid} is an algebroid equipped with a homological morphic vector field $\Xi$. \end{dfn}

\begin{rmk} The space $\Gamma \left(\bigwedge A^*\right)$ of algebroid cochains has a double grading.  If $A$ is a $Q$-algebroid, then $d_A$ and $\Xi_1$ anticommute, so $\left(\Gamma_p \left(\bigwedge^q A^*\right), d_A, \Xi_1 \right)$ is a double complex.  The total differential is $D = d_A + \Xi_1$.
\end{rmk}

\begin{rmk} Given an $\mathcal{LA}$-groupoid $(\Omega \arrows A, G \arrows M)$, there are three different ways to obtain a $Q$-algebroid.  
\begin{enumerate}
\item Apply the $[-1]$ functor to get a $Q$-groupoid, then apply the Lie functor to get a $Q$-algebroid:
\[
\xymatrix{\Omega \ar[r] \ar@<2pt>[d] \ar@<-2pt>[d] & G \ar@<2pt>[d] \ar@<-2pt>[d] \ar@{}[dr]|{\stackrel{[-1]}{\Longrightarrow}} & ([-1]\Omega, d_\Omega) \ar@<2pt>[d] \ar@<-2pt>[d] \ar@{}[dr]|{\stackrel{\mathrm{Lie}}{\Longrightarrow}} & (\mathrm{Lie}([-1]\Omega), \widetilde{d_\Omega}) \ar[d] \\
A \ar[r] & M & ([-1]A, d_A) & ([-1]A, d_A).
}
\]

\item Apply the Lie functor to get a \emph{double Lie algebroid} \cite{mac:dbl}, then apply the $[-1]$ functor to the horizontal bundles:
\[
\xymatrix{\Omega \ar[r] \ar@<2pt>[d] \ar@<-2pt>[d] & G \ar@<2pt>[d] \ar@<-2pt>[d] \ar@{}[dr]|{\stackrel{\mathrm{Lie}}{\Longrightarrow}} & \mathrm{Lie}(\Omega) \ar[r] \ar[d] & \mathrm{Lie}(G) \ar[d] \ar@{}[dr]|{\stackrel{[-1]}{\Longrightarrow}} & ([-1]_{\mathrm{Lie}(G)}\mathrm{Lie}(\Omega), d_{\mathrm{Lie}(\Omega) \to \mathrm{Lie}(G)}) \ar[d]\\
A \ar[r] & M & A \ar[r] & M & ([-1]A, d_A).
}
\]

\item Apply the Lie functor to get a double Lie algebroid, then apply the $[-1]$ functor to the vertical bundles:
\[
\xymatrix{\Omega \ar[r] \ar@<2pt>[d] \ar@<-2pt>[d] & G \ar@<2pt>[d] \ar@<-2pt>[d] \ar@{}[dr]|{\stackrel{\mathrm{Lie}}{\Longrightarrow}} & \mathrm{Lie}(\Omega) \ar[r] \ar[d] & \mathrm{Lie}(G) \ar[d] \ar@{}[dr]|{\stackrel{[-1]}{\Longrightarrow}} & ([-1]_A \mathrm{Lie}(\Omega), d_{\mathrm{Lie}(\Omega) \to A}) \ar[d]\\
A \ar[r] & M & A \ar[r] & M & ([-1]\Lie(G), d_{\Lie(G)}).
}
\]

\end{enumerate}
We conjecture that the first two methods produce the same $Q$-algebroid and that, although the third method produces a different $Q$-algebroid, it results in the same double complex.  The following theorem proves a part of this conjecture; we hope to address the remaining parts in the future.
\end{rmk}

\begin{thm} The first two methods produce the same algebroid.\end{thm}

\begin{proof}
A thorough treatment of this theorem would be rather lengthy; for the sake of brevity, we will provide a sketch of the proof and leave several details to the reader.

Recall from Proposition \ref{prop:vertandlin} that the space of sections of $\Lie(\Omega)$ is spanned by the linear and the vertical sections.  At the same time, a section of $\Lie(\Omega)$ can be viewed as a $t$-tangent vector field along $A$ in $\Omega$ or, in other words, a section of $e^* (T_t \Omega)$.  This allows us to define linear and vertical sections of $e^* (T_t \Omega)$.  From the local forms (\ref{eqn:locallin})-(\ref{eqn:localvert}), we can see that $\xi \in \Gamma(e^* (T_t \Omega))$ is linear if and only if, viewed as a map $\xi: C^\infty(\Omega) \to C^\infty(A)$, it sends $C^\infty_{lin}(\Omega)$ to $C^\infty_{lin}(A)$.  Similarly, $\xi$ is vertical if and only if it sends $C^\infty_{lin}(\Omega)$ to $C^\infty(M)$.  The construction of the left-invariant vector field $X_\xi$ (see Lemma \ref{lemma:leftlong}) and the fact that the source, target, and multiplication maps for $\Omega$ are bundle maps imply the following:
\begin{enumerate}
\item $X_\xi$ is a linear vector field (in the sense of Definition \ref{dfn:linvf}) if and only if $\xi$ is a linear section.
\item $X_\xi$ is a vertical vector field, in the sense that $X_\xi(\alpha) \in C^\infty(G)$ for all $\alpha \in C^\infty_{lin} (\Omega)$, if and only if $\xi$ is a vertical section.
\end{enumerate}
It follows that $\livf(\Omega) = C^\infty(A) \otimes \left(\vect_{lin, LI} (\Omega) \oplus \vect_{V, LI} (\Omega) \right)$, where $\vect_{lin, LI} (\Omega)$ is the space of linear left-invariant vector fields and $\vect_{V, LI} (\Omega)$ is the space of vertical left-invariant vector fields.

Using the isomorphism of Proposition \ref{prop:liniso}, we can identify the $C^\infty(G)$-modules $\vect_{lin} (\Omega)$ and $\vect_{lin}([-1]\Omega)$, and this isomorphism identifies $\vect_{lin, LI} (\Omega)$ and $\vect_{lin, LI} ([-1]\Omega)$.  In a similar manner, we can identify $\vect_{V, LI} (\Omega)$ and $\vect_{V, LI} ([-1]\Omega)$, but with a degree shift.  Thus we have that
\begin{equation} \label{eqn:livfsplit} \livf([-1]\Omega) = C^\infty([-1]A) \otimes \left(\vect_{lin, LI} (\Omega) \oplus [-1]\vect_{V, LI} (\Omega)\right),
\end{equation}
which is, by definition, the space of sections of $\Lie([-1]\Omega)$.  On the other hand, using Proposition \ref{prop:shiftsection}, we can see that (\ref{eqn:livfsplit}) is equal to $\Gamma \left([-1]_{\Lie(G)} \Lie(\Omega)\right)$.  This proves that, as vector bundles over $[-1]A$, $\Lie([-1]\Omega)$ and $[-1]_{\Lie(G)}\Lie(\Omega)$ are isomorphic.

Next, we need to check that the algebroid structures are the same.  Recall that the anchor map is obtained by projecting left-invariant vector fields by the source map.  In the case of $\Lie([-1]\Omega)$, the anchor map is $([-1]s)_*$.  In Proposition \ref{prop:vbashift}, the anchor map for $[-1]_{\Lie(G)}\Lie(\Omega)$ arises from the map
\begin{equation}[-1](s_*): [-1]_{\Lie(G)}\Lie(\Omega) \to [-1]_{TM} TA,
\end{equation}
where we use Proposition \ref{prop:dvb2} to identify $[-1]_{TM} TA$ with $T([-1]A)$.  Thus, agreement of the two anchor maps amounts to the fact that, since $s$ is a bundle map, ``push-forward commutes with $[-1]$,'' which can be verified by using the decomposition (\ref{eqn:livfsplit}).  

The bracket for $\Lie([-1]\Omega)$ is, by definition, the Lie bracket of left-invariant vector fields.  We apply Proposition \ref{prop:hathom} and a similar result for the bracket of linear and vertical vector fields, and we conclude that the correspondences $\vect_{lin, LI} (\Omega) \cong \vect_{lin, LI} ([-1]\Omega)$ and $\vect_{V, LI} (\Omega) \cong \vect_{V, LI} ([-1]\Omega)$ respect the Lie brackets among vertical and horizontal sections.  This agrees with the definition of the bracket on $[-1]_{\Lie(G)}\Lie(\Omega)$, defined in Proposition \ref{prop:vbashift}.  Thus we have shown that $\Lie([-1]\Omega)$ and $[-1]_{\Lie(G)}\Lie(\Omega)$ are isomorphic as algebroids.
\end{proof}

\subsection{Example: The Weil algebra}
\label{ex:weil}

Let $G$ be a Lie group with Lie algebra $\mathfrak{g}$.  Then $[-1]TG$ is a supergroup, and the de Rham operator is a multiplicative vector field on $[-1]TG$.  The Lie algebra of $[-1]TG$ is $[-1]T\mathfrak{g}$, and the induced morphic vector field is the de Rham operator on $[-1]T\mathfrak{g}$.  The space of algebroid forms is $C^\infty([-1]_{\{pt.\}}([-1]T\mathfrak{g}))$.

Let $\{v_i\}$ be a basis for $\mathfrak{g}$ with dual basis $\{\theta^i\}$.  Then $\{\theta^i, \dot{\theta^i}\}$ is a basis for $([-1]T\mathfrak{g})^*$.  Note that $[-1]T\mathfrak{g} = \mathfrak{g} \oplus [-1]\mathfrak{g}$.  After applying the $[-1]$ functor to get $[-1]_{\{pt.\}}([-1]T\mathfrak{g}) = [-1]\mathfrak{g} \oplus [-2]\mathfrak{g}$, $\theta^i$ and  $\dot{\theta^i}$ are degree $1$ and $2$ coordinates, respectively.  The algebra $C^\infty([-1]\mathfrak{g} \oplus [-2]\mathfrak{g})$ is equal to $\bigwedge \mathfrak{g}^* \otimes \bigS \mathfrak{g}^*$.

The algebra $\bigwedge \mathfrak{g}^* \otimes \bigS \mathfrak{g}^*$ is known as the \emph{Weil algebra}, and denoted as $\mathcal{W}$.  The de Rham differential on $[-1]T\mathfrak{g}$ induces the differential operator $d_K = \dot{\theta}^i \pdiff{}{\theta^i}$, which is called the \emph{Koszul operator}.

The Lie algebra differential $d_{[-1]T\mathfrak{g}}$ is, in coordinates,
\begin{equation}
d_{[-1]T\mathfrak{g}} = -\frac{1}{2} c_{ij}^k \theta^i \theta^j \pdiff{}{\theta^k} - c_{ij}^k \theta^i \dot{\theta}^j \pdiff{}{\dot{\theta}^k},
\end{equation}
and the total differential $d_{[-1]T\mathfrak{g}} + d_K$ is known as the Weil differential and is denoted $d_\mathcal{W}$. 

The Weil algebra $\mathcal{W}$, equipped with the Weil differential, is known to be an acyclic complex.  It is, in fact, a model for $EG$, and in order to get a model for $BG$ it is necessary to restrict to the \emph{basic} elements of $\mathcal{W}(\mathfrak{g})$, defined as follows.

For any element $v = a^i v_i \in \mathfrak{g}$, there is a ``contraction operator'' $I_v = a^i \pdiff{}{\theta^i}$.  Since $\mathcal{W}$ represents the algebra of differential forms on $EG$, the contraction operators describe an action of $\mathfrak{g}$ on $EG$.

\begin{dfn} An element $\omega \in \mathcal{W}$ is called
\begin{enumerate}
\item \emph{horizontal} if $I_v \omega = 0$ for all $v \in \mathfrak{g}$,
\item \emph{invariant} if $L_v \omega \defequal [I_v, d_\mathcal{W}] \omega = 0$ for all $v \in \mathfrak{g}$, and
\item \emph{basic} if $\omega$ is both horizontal and invariant.
\end{enumerate}
\end{dfn}

It is clear that the horizontal elements comprise $\bigS\mathfrak{g}^*$.  Furthermore, $d_\mathcal{W}$ vanishes on the basic subcomplex $\left(\bigS\mathfrak{g}^*\right)^G$.

The following theorem is a classic result due to Cartan.

\begin{thm}[\cite{cartan}]\label{thm:cartan}
If $G$ is a compact and connected Lie group, then cohomology of $BG$ is equal to the cohomology of the basic subcomplex of $\mathcal{W}$, which is equal to $\left(\bigS\mathfrak{g}^*\right)^G$.
\end{thm}

\subsection{The double complex of $[-1]TA$}\label{sec:doubleoneta}

An important class of $Q$-algebroids is of the form $[-1]TA \to [-1]TM$, where $A \to M$ is an algebroid (see \S\ref{sec:oneta}).  The previous section dealt with the special case where $A$ is a Lie algebra (or, equivalently, where $M$ is a point).  The goal of the remainder of \S\ref{sec:inf} is essentially to generalize the definitions and results of that case.

Let us first describe the double complex associated to $[-1]TA$.  Using the identification $[-1]_{[-1]TM}[-1]_A TA = [-1]_{[-1]A}T([-1]A)$ (see \S\ref{section:dvb}), the space of cochains in the associated double complex can be identified with the space of differential forms on $[-1]A$.  It was shown in Theorem \ref{thm:lieda} that the algebroid differential of $[-1]TA$ is equal to $\lie_{d_A}$.  The total differential is the sum of the algebroid differential and the given morphic vector field, which in this case is the de Rham operator $d$.  Thus the total complex is
\begin{equation}\label{eqn:totalcom}
\left(C^\infty([-1]T([-1]A)), \lie_{d_A} + d \right).
\end{equation}

\begin{ex}[BRST model of equivariant cohomology, part I]\label{ex:brst1}

Let $M$ be a manifold and let $G$ be a Lie group with a right action on $M$.  As was stated in Proposition \ref{prop:equiv}, the double complex associated to $[-1]T(M \times G)$, where $M \times G$ is the action groupoid, computes the equivariant cohomology $H_G^\bullet (M)$.  

The associated algebroid is $[-1]T(M \times \mathfrak{g}) = [-1]TM \times [-1]T\mathfrak{g}$.  The vector field $d_M \times d_K$, where $d_M$ is the de Rham differential on $M$ and $d_K$ is the Koszul operator, is a morphic vector field, so $[-1]TM \times [-1]T\mathfrak{g}$ is a $Q$-algebroid.  The space of algebroid forms is $C^\infty\left([-1]_{[-1]TM}([-1]TM \times [-1]T\mathfrak{g})\right)$, which is equal to $\Omega(M) \otimes \mathcal{W}(\mathfrak{g})$.

The total differential $D_B$ may be written as
\begin{equation}\label{brstdiff}
D_B = d_M + d_{\mathcal{W}} + \theta^i \lie_{\rho(v_i)} - \dot{\theta}^i \iota_{\rho(v_i)},
\end{equation}
where $d_{\mathcal{W}}$ is the Weil differential and $\rho: \mathfrak{g} \to \vect(M)$ describes the infinitesimal action of $\mathfrak{g}$ on $M$.

The algebra $\Omega(M) \otimes \mathcal{W}(\mathfrak{g})$ and the differential (\ref{brstdiff}) form the \emph{BRST model} of equivariant cohomology \cite{kalkman} \cite{mq}.  In fact, this complex is a model for $M \times EG$ (hence its cohomology is equal to $H^\bullet(M)$), and one must restrict to a suitably defined basic subcomplex to obtain equivariant cohomology.  We will return to this example in \S\ref{sub:basic}, after defining the basic subcomplex in the general setting.
\end{ex}

\begin{rmk}
A more well-known model of equivariant cohomology is the Weil model, which has the same algebra $\Omega(M) \otimes \mathcal{W}(\mathfrak{g})$ as the BRST model, but the simpler differential $d_M + d_K$.  In his thesis, Kalkman \cite{kalkman} describes an extension of the Mathai-Quillen isomorphism that relates the BRST differential and the Weil model differential, thus showing that the two models are equivalent\footnote{The original isomorphism of Mathai and Quillen \cite{mq} sent the Cartan model, which may be identified with the basic subcomplex of the BRST model, to the basic subcomplex of the Weil model.
 Kalkman extended the isomorphism to the total complexes.}.
\end{rmk}

In order to compute the cohomology of the total complex (\ref{eqn:totalcom}), we will use the following generalization of the Mathai-Quillen-Kalkman isomorphism.  Since $d_A$ is a degree $1$ vector field on $[-1]A$, it follows that $\iota_{d_A}$ is a degree $0$ vector field on $[-1]T([-1]A)$.

\begin{dfn} The \emph{generalized Mathai-Quillen-Kalkman isomorphism} is $\gamma \defequal \exp (\iota_{d_A})$.
\end{dfn}

\begin{lemma} The de Rham differential $d$ is $\gamma$-related to $d + \lie_{d_A}$.
\end{lemma}
\begin{proof} Since $\iota_{d_A}$ is nilpotent, the identity $\Ad_{\exp (\iota_{d_A})} = \exp (\ad_{\iota_{d_A}})$ holds.  From the Cartan commutation relations, it is immediate that
\begin{align}
\ad_{\iota_{d_A}} (d) &\defequal [\iota_{d_A}, d] = \lie_{d_A}, \\
\ad^2_{\iota_{d_A}} (d) &= [\iota_{d_A}, \lie_{d_A}] = -\iota_{[d_A,d_A]} = 0,
\end{align}
and it follows that $\Ad_{\gamma} (d) = d + \lie_{d_A}$.
\end{proof}

\begin{rmk} The differential for the Weil model is sometimes taken to be $d_M + d_\mathcal{W}$, where $d_\mathcal{W}$ is the Weil differential (see Example \ref{ex:weil}).  The isomorphism given by Kalkman relates this differential to the BRST differential.  Thus, to compare his isomorphism with ours, it is necessary to use an automorphism of $\mathcal{W}(\mathfrak{g})$ that relates $d_K$ and $d_\mathcal{W}$.  Such an isomorphism is well-known (e.g. \cite{gs}), but we point out that it is a special case, when $A = \mathfrak{g}$, of the isomorphism $\gamma$.
\end{rmk}

\begin{cor} \label{cor:eg} The cohomology of the total complex (\ref{eqn:totalcom}) is equal to $H^\bullet(M)$.
\end{cor}

\begin{proof}
The isomorphism $\gamma$ provides an isomorphism of the total complex and the de Rham complex of $[-1]A$.  Using the Euler vector field for the vector bundle $[-1]A \to M$, one can construct a chain homotopy between $\Omega([-1]A)$ and $\Omega(M)$.  For details, see \cite{tuynman}.
\end{proof}

\begin{rmk}
If $A$ is the algebroid of a Lie groupoid $G$, then the result of Corollary \ref{cor:eg} is consistent with the assertion that the total complex is a model for $EG$.  A model for $BG$ is obtained by restricting to a basic subcomplex, which will be described in \S\ref{sub:basic}.
\end{rmk}

\subsection{Connections and prolongations}
Let $E \stackrel{\pi}{\to} M$ be a vector bundle.  Recall that $\vect_{\pi}(E)$ denotes the space of vector fields that are tangent to the $\pi$-fibres.  We refer to elements of $\vect_{\pi}(E)$ as \emph{vertical vector fields}.

\begin{dfn}
A differential form $\omega \in \Omega(E)$ is \emph{horizontal} (with respect to the projection $\pi$) if $\iota_X \omega = 0$ for all $X \in \vect_\pi(E)$.  The subalgebra of horizontal forms is denoted by $\Omega_H(E)$.
\end{dfn}

Consider the subalgebra $\pi^*(\Omega(M))$ of $\Omega(E)$.  All forms in this subalgebra are horizontal (In fact, it is by definition the subalgebra of forms that are basic for the projection $\pi$).  It is a simple exercise to see that $\pi^*(\Omega(M))$ spans the space of horizontal forms, in the sense that 
\begin{equation}\label{eqn:hor}
\Omega_H(E) = C^\infty(E) \otimes_{C^\infty(M)} \pi^*(\Omega(M)).
\end{equation}

A connection on $E$ may be described by a linear $1$-form $P \in \Omega^1 (E; T_\pi E)$, which assigns to a vector $v \in TE$ its vertical component.  Equivalently, we may view $P$ as a $C^\infty(E)$-module projection from $\vect(E)$ to $\vect_{\pi}(E)$.  

By identifying $\vect_{\pi}(E)$ with $\Gamma(\pi^*E) = C^\infty(E) \otimes_{C^\infty(M)} \Gamma(E)$, we have the dual map
\begin{equation}
P^*: C^\infty(E) \otimes_{C^\infty(M)} \Gamma(E^*) \to \Omega^1(E),
\end{equation}
which may be extended as an algebra homomorphism to
\begin{equation}\label{eqn:vert}
\bar{P}: C^\infty(E) \otimes_{C^\infty(M)} \bigwedge \Gamma(E^*) \to \Omega(E).
\end{equation}

\begin{dfn}
A differential form $\omega \in \Omega(E)$ is called \emph{vertical} if $\omega \in \im \bar{P}$.  The subalgebra of vertical forms is denoted by $\Omega_V (E)$.
\end{dfn}

The introduction of a connection thus determines a splitting 
\[\Omega(E) = \Omega_H(E) \otimes_{C^\infty(E)} \Omega_V(E),\]
and, using (\ref{eqn:hor}) and (\ref{eqn:vert}), we may make the identification
\begin{equation}\label{eqn:split1}
\Omega(E) = C^\infty(E) \otimes_{C^\infty(M)} \Omega(M) \otimes_{C^\infty(M)} \bigwedge \Gamma(E^*).
\end{equation}
In the language of supermanifolds, (\ref{eqn:split1}) may be rewritten as
\begin{equation}\label{eqn:split2}
[-1]TE = E \oplus_M [-1]TM \oplus_M [-1]E.
\end{equation}
We thus have the following result.

\begin{lemma}
Let $E \to M$ be a vector bundle.  Then a choice of connection on $E$ induces
\begin{enumerate}
\item a vector bundle structure on $[-1]TE \to M$ and
\item an injective $C^\infty(M)$-module homomorphism $a: \Gamma(E) \to \Gamma([-1]TE, M)$.
\end{enumerate}
\end{lemma}

\subsection{The basic subcomplex}\label{sub:basic}

Let $A \to M$ be an algebroid equipped with a connection.  Applying (\ref{eqn:split2}) to the bundle $[-1]A \to M$, we obtain a decomposition
\begin{equation}
[-1]T_{[-1]A}([-1]A) = [-1]A \oplus [-1]TM \oplus [-2]A.
\end{equation}
Applying the $[1]$ functor over $M$ yields
\begin{equation}\label{eqn:split3}
[1]_M \left([-1]T_{[-1]A}([-1]A)\right) = A \oplus TM \oplus [-1]A.
\end{equation}

Let $a: \Gamma(A) \hookrightarrow \Gamma([1]_M([-1]T_{[-1]A}([-1]A)))$ be the inclusion map arising from the splitting (\ref{eqn:split3}).  Then composing $a$ with contraction gives, for any $X \in \Gamma(A)$, a degree $-1$ vector field $I_X \in \vect_{-1}([-1]T([-1]A))$.

\begin{dfn} An element $\omega \in C^\infty([-1]([-1]TA))$ is called
\begin{enumerate}
\item \emph{horizontal} if $I_X \omega = 0$ for all $X \in \Gamma(A)$,
\item \emph{invariant} if $L_X \omega \defequal [I_X, D] \omega = 0$ for all $X \in \Gamma(A)$, and
\item \emph{basic} if $\omega$ is both horizontal and invariant.
\end{enumerate}
\end{dfn}

\begin{ex}[BRST model, part II]\label{ex:brst2}
Recall the BRST model of Example \ref{ex:brst1}.  The action algebroid $M \times \mathfrak{g} \to M$ has a canonical flat connection and thus a natural splitting
\begin{equation}
[-1]T(M \times [-1]\mathfrak{g}) = [-1]\mathfrak{g} \times [-1]TM \times [-2]\mathfrak{g},
\end{equation}
or, equivalently,
\begin{equation}
\Omega(M \times [-1]\mathfrak{g}) = \bigwedge \mathfrak{g}^* \otimes \Omega(M) \otimes \bigS \mathfrak{g}^*.
\end{equation}

If $\{X_i\}$ is the global frame of flat sections corresponding to the basis $\{v_i\}$ of $\mathfrak{g}$,
then the contraction operators $I_X$ are generated freely as a $C^\infty(M)$-module by those of the form
\begin{equation}
I_{X_i} = \pdiff{}{\theta^i}.
\end{equation}
Thus the horizontal elements are simply those that do not depend on any $\theta^i$, i.e.\ those that lie in $\Omega(M) \otimes S \mathfrak{g}^*$.  On the subalgebra of horizontal elements, 
\begin{equation}
L_{X_i} = \left[\pdiff{}{\theta^i}, D_B\right] = -\dot{\theta}^j c_{ij}^k \pdiff{}{\dot{\theta}^k} + \lie_{\rho(v_i)}.
\end{equation}
On the basic subalgebra $\left(\Omega(M) \otimes S \mathfrak{g}^*\right)^G$, the differential becomes the Cartan differential
\begin{equation}
d_C \defequal d_M - \dot{\theta}^i \iota_{\rho(v_i)}.
\end{equation}

We have thus recovered the well-known fact \cite{gs}\cite{kalkman}\cite{mq} that the basic subcomplex of the BRST model is equal to the Cartan model.
\end{ex}

\subsection{A geometric construction for the total differential}\label{sec:geom}

Consider $[-1]T\mathfrak{g}$, as in Example \ref{ex:weil}.  As was pointed out in \cite{gs}, there is in fact a larger superalgebra $\widetilde{\mathfrak{g}}$ which, as a vector space, is equal to $[1]\reals \oplus [-1]T\mathfrak{g}$.  If $\epsilon$ is the standard basis vector of $[1]\reals$, then the Lie bracket is defined by
\begin{align}
[v^C, w^C] &= [v,w]^C, & [v^C, w^V] &= [v,w]^V, \\
[v^V, w^V] &= 0, & [v^V, \epsilon] &= v^C \\
[v^C, \epsilon] &= 0, & [\epsilon, \epsilon] &= 0.
\end{align}

The superalgebra $\widetilde{\mathfrak{g}}$ arises in the context of (infinitesimal) Lie algebra actions in the following way.  An action $\mathfrak{g} \to \vect(M)$ naturally induces an action $\widetilde{\mathfrak{g}} \to \mathrm{Der} (\Omega(M))$, where $\epsilon$ maps to the de Rham operator.

If $\gamma$ is the coordinate dual to $\epsilon$, then the Lie algebra differential is
\begin{equation}
d_{\widetilde{\mathfrak{g}}} = d_{[-1]T\mathfrak{g}} - \gamma d_K.
\end{equation}
Note that in $C^\infty([-1]\widetilde{\mathfrak{g}})$, the coordinate $\gamma$ is of degree $0$, so there is a one-parameter family of differentials $d_{\widetilde{\mathfrak{g}}, \gamma}$.  In particular, $d_{\widetilde{\mathfrak{g}}, 0} = d_{[-1]T\mathfrak{g}}$ and $d_{\widetilde{\mathfrak{}g}, -1} = d_{\mathcal{W}}$.  In this sense, the superalgebra $\widetilde{\mathfrak{g}}$ interpolates between $[-1]T\mathfrak{g}$ and a version of $[-1]T\mathfrak{g}$ that is ``twisted'' by $d_K$.

This construction extends to the $Q$-algebroid setting.  Let $A \to M$ be a $Q$-algebroid with morphic vector field $\Xi$ with base field $\phi$.

The superalgebroid $\widetilde{A}$ is defined, as a bundle, to be $A \oplus [1]L \to M$, where $[1]L \defequal [1]\reals \times M$ is the trivial degree $1$ line bundle over $M$.  In terms of the standard basis vector $\epsilon$ of $[1]\reals$, the algebroid structure of $\widetilde{A}$ extends the algebroid structure of $A$ by the relations
\begin{align} \label{eqn:twist}
\widetilde{\rho}(\epsilon) = \phi, && [X,\epsilon] = (-1)^{|X|-1}D_\Xi X, && [\epsilon, \epsilon] = 0.
\end{align}
The fact that these relations define an algebroid structure follows from the following proposition:

\begin{prop}
The differential $d_{\widetilde{A}} = d_A + \gamma \Xi_1$ corresponds to the anchor and bracket relations (\ref{eqn:twist}).
\end{prop}

\begin{proof}
Using the fact that $\iota_\epsilon = \pdiff{}{\gamma}$, we have for $f \in C^\infty(M)$ that
\begin{equation}
\begin{split}
\iota_\epsilon d_{\widetilde{A}}f &= \pdiff{}{\gamma}(d_A f + \gamma \phi f) \\
&= \phi(f),
\end{split}
\end{equation}
which is equal to $\widetilde{\rho}(\epsilon) (f)$, and for $\omega \in \Gamma(A^*)$ that
\begin{equation}
\begin{split}
& \rho(X) \iota_\epsilon \omega - (-1)^{|X|}\phi \iota_X \omega - \iota_\epsilon \iota_X d_{\widetilde{A}}\omega \\
&= \rho(X) \pdiff{\omega}{\gamma} - (-1)^{|X|}\phi \iota_X \omega - \pdiff{}{\gamma}\left[ \iota_X (d_A + \gamma \Xi_1)\omega \right] \\
&= \rho(X) \pdiff{\omega}{\gamma} - (-1)^{|X|}\phi \iota_X \omega - \left[ \iota_X d_A \pdiff{\omega}{\gamma} + \pdiff{}{\gamma}(\gamma \iota_X  \Xi_1\omega )\right] \\
&= - (-1)^{|X|}\phi \iota_X \omega - \iota_X  \Xi_1\omega \\
&= (-1)^{|X| -1} \iota_{D_\Xi X} \omega,
\end{split}
\end{equation}
which is equal to $\iota_{[X, \epsilon]}\omega$.  Thus equations (\ref{eqn:anchor}) and (\ref{eqn:bracket}), which relate the differential to the bracket and anchor, hold.
\end{proof}

		\chapter{$Q$-groupoids and $Q$-algebroids in Poisson geometry}
		\label{chapter:examples}
\label{ch:examples}
\section{Poisson-Lie groups, Lie bialgebras, and Drinfel'd doubles}

\begin{dfn}[\cite{drinfeld:ybe}, \cite{drinfeld:quantum}] A Lie group $G$ with a Poisson structure $\pi$ is a \emph{Poisson-Lie group} if $\pi$ is a multiplicative bivector field.
\end{dfn}

Let $G$ be a Lie group with a Poisson bivector $\pi$.  Since $G$ is, in particular, a Poisson manifold, there is a Lie algebroid structure on the vector bundle $T^*G \to G$.  It was shown by Mackenzie \cite{mac:dblie1} that $G$ is a Poisson-Lie group if and only if $\mathfrak{g}^*$ possesses a Lie algebra structure such that the following is an $\mathcal{LA}$-groupoid:
\begin{equation}
\xymatrix{T^*G \ar[r] \ar@<2pt>[d] \ar@<-2pt>[d] & G \ar@<2pt>[d] \ar@<-2pt>[d] \\
\mathfrak{g}^* \ar[r] & \{pt.\}}
\end{equation}
The corresponding $Q$-groupoid is $[-1]T^*G \arrows [-1]\mathfrak{g}^*$ with homological vector field $d_\pi$, the Poisson cohomology operator.   

Consider the resulting $Q$-algebroid $[-1]T^*\mathfrak{g} = [-1]\mathfrak{g}^* \oplus \mathfrak{g} \to [-1]\mathfrak{g}^*$.  This is the action algebroid for the right coadjoint action of $\mathfrak{g}$ on $[-1]\mathfrak{g}^*$.  In coordinates $\{\theta^i, v_i\}$, where $\{v_i\}$ and $\{\theta^i\}$ are ``dual'' coordinates on $[-1]\mathfrak{g}^*$ and $[-1]\mathfrak{g}$, respectively, the algebroid differential $d_{[-1]T^*\mathfrak{g}} \in \vect\left([-1]\mathfrak{g}^* \oplus [-1]\mathfrak{g}\right)$ is\footnote{The differential $d_{[-1]T^*\mathfrak{g}}$ may be more simply viewed as the differential for the Lie algebra cohomology of $\mathfrak{g}$ with coefficients in $\bigwedge \mathfrak{g}$.}
\begin{equation}
d_{[-1]T^*\mathfrak{g}} = \theta^i c_{ij}^k v_k \pdiff{}{v_j} - \frac{1}{2} \theta^i \theta^j c_{ij}^k \pdiff{}{\theta^k},
\end{equation}
where $c_{ij}^k$ are the structure constants for $\mathfrak{g}$.  If $\gamma^{ij}_k$ are the structure constants for $\mathfrak{g}^*$, then the homological vector field is
\begin{equation}
\Xi_{[-1]T^*\mathfrak{g}} = v_i \gamma^{ij}_k \theta^k \pdiff{}{\theta^j} - \frac{1}{2}v_i v_j \gamma^{ij}_k \pdiff{}{v_k}.
\end{equation}

The property
\begin{equation}
\left[d_{[-1]T^*\mathfrak{g}}, \Xi_{[-1]T^*\mathfrak{g}}\right] = 0
\end{equation}
is equivalent to the compatibility condition for the pair $(\mathfrak{g},\mathfrak{g}^*)$ to be a Lie bialgebroid \cite{ks1} \cite{lr}.

The algebra of cochains for the double complex is $C^\infty([-1]\mathfrak{g}^* \oplus [-1]\mathfrak{g}) = \bigwedge \mathfrak{g} \otimes \bigwedge \mathfrak{g}^*$, and the total differential is in fact equal to the Chevalley-Eilenberg differential for the Drinfel'd \cite{drinfeld:ybe} double $\mathfrak{d} = \mathfrak{g}^* \oplus \mathfrak{g}$.  Thus the Drinfel'd double of the Lie bialgebra may be recovered as
\begin{equation}
\mathfrak{d} = [1]\left([-1]\mathfrak{g}^* \oplus [-1]\mathfrak{g}\right),
\end{equation}
where the total differential gives $[-1]\mathfrak{g}^* \oplus [-1]\mathfrak{g}$ the structure of a Lie antialgebra (see Remark \ref{rmk:anti}).

In this case, the horizontal condition for cochains is simply that contraction with elements of $\mathfrak{g}^*$ is zero, or, in coordinates, that $\pdiff{\omega}{v_i} = 0$ for all $i$.  The resulting invariance condition is
\begin{equation} \gamma_k^{ij} \theta^k \pdiff{\omega}{\theta^j} = 0.
\end{equation}
Therefore the basic forms are the elements of $\bigwedge \mathfrak{g}^*$ that are invariant under the $\ad_{\mathfrak{g}^*}$ action.  On the basic subcomplex, the total differential agrees with the Chevalley-Eilenberg differential for $\mathfrak{g}$.  Thus the basic cohomology is the $\mathfrak{g}^*$-invariant Lie algebra cohomology of $\mathfrak{g}$.

The above discussion may be generalized to the situation of Poisson groupoids and their infinitesimal counterparts, Lie bialgebroids.  Specifically, if $(A, A^*)$ is a Lie bialgebroid, then we may form a $Q$-algebroid $[-1]T^*A \to [-1]A^*$, whose double complex reproduces Roytenberg's \cite{royt} ``commuting Hamiltonians'' approach to bialgebroids.  The details of this case will appear in \cite{gm}.

\section{Equivariant Poisson cohomology}

Cartan's formulation of equivariant cohomology is not limited to de Rham cohomology; in fact, one can replace $\Omega(M)$ with any \emph{$G$-differential algebra}, i.e. a differential algebra equipped with a compatible\footnote{More precisely, a $G$-differential algebra is a graded commutative algebra equipped with an action of $\widetilde{\mathfrak{g}}$, as in \S\ref{sec:geom}.} action of $[-1]T\mathfrak{g}$.  Ginzburg \cite{ginzburg} used this fact to define the concept of \emph{equivariant Poisson cohomology}\footnote{Ginzburg introduced equivariant Poisson cohomology with the hope of finding cohomological obstructions to the existence of moment maps for Poisson actions; the goal was not achieved.  It would be interesting to see whether the results here provide any insight toward that problem.}.  

In the first part of this section, we show that equivariant Poisson cohomology may be formulated as the basic cohomology of the double complex of a $Q$-algebroid.  We will actually describe the more general \emph{equivariant algebroid cohomology}, which was also addressed in \cite{ginzburg}, associated to an $A$-action of a Lie algebra $\mathfrak{g}$ on an algebroid $A \to M$.  The term ``equivariant Poisson cohomology'' refers to the case $A = T^*M$, where $M$ is a Poisson manifold, but it is worth noting that the case $A=TM$ coincides with the regular notion of equivariant cohomology, as described in Examples \ref{ex:brst1} and \ref{ex:brst2}.

The infinitesimal models of equivariant cohomology are only valid for the actions of Lie groups which are connected and compact.  In \S\ref{sub:lagpd}, we propose a model for the equivariant algebroid cohomology of a Lie group action.  This model takes the form of the double complex of a $Q$-groupoid.  This new model allows us to make sense, for example, of the equivariant Poisson cohomology of a discrete group action.

\subsection{The infinitesimal model}\label{sub:poissinf}

Let $A \to M$ be an algebroid and let $\mathfrak{g}$ be a Lie algebra. 

\begin{dfn} A right \emph{$A$-action} of $\mathfrak{g}$ on $M$ is a Lie algebra homomorphism $\widetilde{a}: \mathfrak{g} \to \Gamma(A)$.
\end{dfn}

\begin{rmk} Via composition with the anchor map, an $A$-action $\widetilde{a}$ induces a Lie algebra homomorphism $a \defequal \rho \circ \widetilde{a}: \mathfrak{g} \to \vect(M)$ that describes an action of $\mathfrak{g}$ on $M$.  Ginzburg \cite{ginzburg} begins with an action map $a: \mathfrak{g} \to \vect(M)$ and defines an \emph{equivariant pre-momentum mapping} to be a lift $\widetilde{a}$ of the action on $M$ to an $A$-action.  An equivariant pre-momentum mapping, defined in this manner, is equivalent to an algebroid morphism from the action algebroid $M \times \mathfrak{g}$ to $A$.
\end{rmk}

\begin{thm}Let $\widetilde{a}: \mathfrak{g} \to \Gamma(A)$ be an $A$-action.  Then there is an induced Lie superalgebra action $\bar{\rho}: [-1]T\mathfrak{g} \to \vect([-1]A)$.  The associated action algebroid $[-1]A \times [-1]T\mathfrak{g} \to [-1]A$ is a $Q$-algebroid with morphic vector field $d_A \times d_K$.
\end{thm}

\begin{proof}
Recall that  $[-1]T\mathfrak{g}$ is naturally isomorphic to $\mathfrak{g} \oplus [-1]\mathfrak{g}$, where the first summand consists of complete lifts $v^C$ and the second summand consists of vertical lifts $v^V$ of elements $v \in \mathfrak{g}$.

The induced action is defined as
\begin{align}
\bar{\rho}(v^V) = \iota_{\widetilde{a}(v)}, && \bar{\rho}(v^C) = \lie_{\widetilde{a}(v)} \defequal [\iota_{\widetilde{a}(v)}, d_A].
\end{align}
It is clear that $\bar{\rho}$ is a Lie algebra homomorphism and thus describes an action of $[-1]T\mathfrak{g}$ on $[-1]A$.

Consider the associated action algebroid $[-1]A \times [-1]T\mathfrak{g}$.  The algebroid differential $d_{[-1]A \times [-1]T\mathfrak{g}}$ is a homological vector field on $[-1]_{[-1]A}\left([-1]A \times [-1]T\mathfrak{g}\right) = [-1]A \times [-1]\mathfrak{g} \times [-2]\mathfrak{g}$.

Let $\{v_b\}$ be a basis for $\mathfrak{g}$ with dual basis $\{\theta^b\}$.  Then $\{\theta^b, \dot{\theta}^b\}$ forms a set of linear coordinates on $[-1]T\mathfrak{g}$, and we may view $\{\theta^b\}$ and $\{\dot{\theta}^b\}$, respectively, as degree $1$ and degree $2$ coordinates on $[-1]\mathfrak{g} \times [-2]\mathfrak{g}$.  

In terms of the basis $\{v_b\}$, the algebroid differential may be written as 
\begin{equation}
d_{[-1]A \times [-1]T\mathfrak{g}} = \theta^b \lie_{\widetilde{a}(v_b)} - \dot{\theta}^b \iota_{\widetilde{a}(v_b)} + d_{[-1]T\mathfrak{g}}.
\end{equation}
The Koszul operator $d_K = \dot{\theta}^b \pdiff{}{\theta^b}$ commutes with $d_{[-1]T\mathfrak{g}}$, so it is immediate that 
\begin{equation}\label{eqn:koszulcomm}
[d_{[-1]A \times [-1]T\mathfrak{g}}, d_K] = \dot{\theta}^b \lie_{\widetilde{a}(v_b)}.
\end{equation}
Meanwhile, $\lie_{\widetilde{a}(v_b)}$ commutes with $d_A$, so
\begin{equation}\label{eqn:liedcomm}
[d_{[-1]A \times [-1]T\mathfrak{g}}, d_A] = -\dot{\theta}^b [\iota_{\widetilde{a}(v_b}, d_A] = -\dot{\theta}^b \lie_{\widetilde{a}(v_b)}.
\end{equation}

From (\ref{eqn:koszulcomm}) and (\ref{eqn:liedcomm}) it follows that $d_A \times d_K$ is a morphic vector field on $[-1]A \times [-1]T\mathfrak{g}$.  Clearly $d_A \times d_K$ is a homological vector field, so $[-1]A \times [-1]T\mathfrak{g}$ has the structure of a $Q$-algebroid.
\end{proof}

\begin{rmk}
The Mathai-Quillen isomorphism in this context is $\exp Q$, where 
\begin{equation}Q = \theta^b \iota_{\widetilde{a}(v_b)} - \frac{1}{2}\theta^a\theta^b c_{ab}^e \pdiff{}{\dot{\theta}^e}.
\end{equation}
The differential $d_A \times d_K$ is $\exp Q$-related to the total differential $D = d_{[-1]A \times [-1]T\mathfrak{g}} + d_A \times d_K$, and it follows that the cohomology of the total complex is equal to the algebroid cohomology of $A$.  The details are left to the reader.
\end{rmk}

The next step is to describe the basic subcomplex.  In terms of the canonical splitting $[-1]A \times [-1]\mathfrak{g} \times [-2]\mathfrak{g}$, there is an induced action $I_v$, for any element $v \in \mathfrak{g}$, acting by contraction in the middle component.  In local coordinates, the action is simply $I_{v_b} = \pdiff{}{\theta^b}$.  

As in the previous sections, the basic elements of $C^\infty([-1]A \times [-1]\mathfrak{g} \times [-2]\mathfrak{g})$ are those which are annihilated by $I_v$ and $L_v \defequal [I_v, D]$ for all $v \in \mathfrak{g}$.  The elements annihilated by every $I_v$ are those which do not depend on any $\theta^b$.  The basic elements thus form a subalgebra $\left(C^\infty([-1]A) \otimes S(\mathfrak{g}^*) \right)^G \subseteq C^\infty([-1]A) \otimes S(\mathfrak{g}^*)$, and on this subalgebra the total differential is
\begin{equation}\label{eqn:pcdiff}
d_C = d_A - \dot{\theta}^b \iota_{\widetilde{a}(v_b)}
\end{equation}

The basic subcomplex equipped with the differential (\ref{eqn:pcdiff}) is identical to the complex that Ginzburg \cite{ginzburg} introduced to define equivariant Poisson cohomology.

\subsection{An $\mathcal{LA}$-groupoid construction}\label{sub:lagpd}

Let $A \to M$ be a Lie algebroid and let $G$ be a Lie group.  Let $d_K$ denote the de Rham operator on $\Omega(G)$, viewed as a homological vector field on $[-1]TG$.

\begin{dfn}\label{dfn:gpaction}
An \emph{$A$-action} of $G$ is a (right) action of $TG$ on $A$ such that the action map $\widetilde{s}: A \times TG \to A$ is an algebroid morphism.
\end{dfn}

\begin{prop} Let $\widetilde{s}$ be an $A$-action of $G$.  Then there exists an action of $G$ on $M$ such that the following diagram describes an $\mathcal{LA}$-groupoid:
\begin{equation}
\xymatrix{A \times TG \ar[r] \ar@<2pt>[d] \ar@<-2pt>[d] & M \times G \ar@<2pt>[d] \ar@<-2pt>[d] \\
A \ar[r] & M.}
\end{equation}
\end{prop}
\begin{proof}
The result is immediate from Definition \ref{dfn:gpaction} and the fact that $TG \to G$ is an ``$\mathcal{LA}$-group''.
\end{proof}

\begin{rmk}
It follows that the action groupoid $[-1]A \times [-1]TG \arrows [-1]A$, with action map $[-1]\widetilde{s}$, is a $Q$-groupoid.  The multiplicative vector field on $[-1]A \times [-1]TG$ is $d_A \times d_K$.
\end{rmk}

\begin{ex}
Given an action map $s: M \times G \to M$, then $Ts: TM \times TG \to TM$ is the unique $TM$-action that lifts $s$.
\end{ex}

\begin{ex}
If $G$ is discrete, then any action $\widetilde{s}: A \times G \to A$, where $G$ acts by algebroid automorphisms, is an $A$-action.
\end{ex}

\begin{dfn}
Let $\widetilde{s}: A \times TG \to A$ be an $A$-action.  The \emph{equivariant algebroid cohomology} is the total cohomology of the double complex associated to the $Q$-groupoid $[-1]A \times [-1]TG \arrows [-1]A$.
\end{dfn}

\section{Poisson $G$-spaces}

Let $G$ be a Poisson-Lie group acting on a Poisson manifold $M$.  When the action is Poisson, there is an associated Lie algebroid structure \cite{lu} on the bundle $A = (M \times \mathfrak{g}) \oplus T^*M$, which forms an example of a ``matched pair\footnote{A matched pair of Lie algebroids is essentially the same thing as a double Lie algebroid satisfying a \emph{vacancy} condition (see \cite{mac:dbl2}).}'' of Lie algebroids \cite{mokri}.  Recently, Bursztyn and Crainic \cite{bu-cr} have extended Lu's construction to quasi-Poisson $G$-spaces.

As was pointed out by Lu \cite{lu}, when the Lie algebroid $T^*M$ is integrable, then the Lie algebroid $A$ integrates to a double Lie groupoid.  In this section, we will describe the intermediate object, which is an $\mathcal{LA}$-groupoid.  In a way, this $\mathcal{LA}$-groupoid presents the ``best of both worlds''; it incorporates the global action of $G$ on $M$, as opposed to the infinitesimal action $\mathfrak{g} \to \vect(M)$, but does not require a (possibly nonexistant) symplectic groupoid that integrates $T^*M$.

\subsection{The algebroid structure of $s^*(T^*M)$}
\label{sub:gspacealgbd}
Let $(M, \alpha)$ be a Poisson manifold, and let $(G, \pi)$ be a Poisson-Lie group with a right Poisson action $s: M \times G \rightarrow M$.  There is a natural algebroid structure on the pullback bundle $s^*(T^*M)$, as follows.

The module of sections is $\Gamma(s^*(T^*M)) = C^\infty(M \times G) \otimes_s \Omega^1(M)$.  For $\omega \in \Omega^1(M)$, the anchor satisfies
\begin{equation}\label{eqn:bcanch}
\rho(1 \otimes \omega) = \widetilde{\pi}^\sharp s^* \omega,
\end{equation}
where $\widetilde{\pi} \defequal \alpha \times \pi$ is the Poisson bivector field on $M \times G$ and $\widetilde{\pi}^\sharp$ is the associated map from $\Omega^1(M \times G)$ to $\vect(M \times G)$.  The anchor map is extended by $C^\infty(M\times G)$-linearity to all sections.

For $\omega, \eta \in \Omega^1(M)$, the bracket satisfies
\begin{equation}\label{eqn:bcbrack}
[ 1 \otimes \omega, 1 \otimes \eta] = 1 \otimes [\omega, \eta]_\alpha,
\end{equation}
where $[\cdot, \cdot]_\alpha$ is the Koszul bracket arising from the Poisson structure on $M$.  The bracket is extended by the Leibniz rule to all sections.  The anchor identity follows from the Poisson condition, which is that $\alpha$ and $\tilde{\pi}$ are $s$-related.  The Jacobi identity follows from the anchor identity and the fact that the Koszul bracket satisfies the Jacobi identity.

\begin{rmk} Even if $s$ is not a Poisson map, we may define a bracket and anchor by the equations (\ref{eqn:bcanch}) and (\ref{eqn:bcbrack}).  The Jacobi identity will not be satisfied, so we will not have a Lie algebroid structure on $s^*(T^*M)$; however, we may still apply the $[-1]$ functor and obtain a degree $1$ vector field, which will not be homological.  In what follows, we will not assume that $s$ is Poisson.
\end{rmk}

\subsection{The $Q$-manifold structure of $s^*([-1]T^*M)$}

We may apply the $[-1]$ functor to $s^*(T^*M) \to M \times G$ to get a degree $1$ vector field, which will denote by $d_{\widetilde{\pi}}$.  Since the $[-1]$ functor commutes with pullbacks, we may make the identification $[-1]s^*(T^*M) = s^*([-1]T^*M)$, whose algebra of functions is $\Gamma\left(s^*\left(\bigwedge TM\right)\right) = C^\infty(M \times G) \otimes_s \vect^\bullet (M)$.

The action of $d_{\widetilde{\pi}}$ on elements of the form $f \otimes 1$ is given by
\begin{equation}
d_{\widetilde{\pi}}(f \otimes 1) = s_* [\widetilde{\pi}, f] \in \Gamma(s^*(TM)), \\
\end{equation}
where the right hand side is the push-forward by $s$ of the (negative of the) Hamiltonian vector field of $f$.  When $s$ is a Poisson map, we have the simple equation
\begin{equation}\label{eqn:bcdiffvf}
d_{\widetilde{\pi}}(1 \otimes X) = 1 \otimes [\alpha, X] = 1 \otimes d_\alpha X.\\
\end{equation}
Even when $s$ is not a Poisson map, $d_{\widetilde{\pi}}$ satisfies the property that for any $Z \in \vect(M \times G)$, 
\begin{equation}\label{eqn:bcdiffpush}
d_{\widetilde{\pi}}(s_* Z) = s_* [\widetilde{\pi}, Z]
\end{equation}

\subsection{The groupoid structure of $s^*(T^*M)$}
\label{sub:gspacegpd}

A point in $s^*(T^*M)$ is of the form $(x,g,\eta)$, where $\eta \in T^*_{xg}M$.  There is a groupoid $s^*(T^*M) \arrows T^*M$, defined as follows.  The source and target maps $\widetilde{s}$ and $\widetilde{t}$, respectively, are defined as
\begin{equation}\begin{split}\label{eqn:gspacest}
\widetilde{s}: (x,g, \eta) &\mapsto (xg, \eta), \\
\widetilde{t}: (x,g, \eta) &\mapsto (x, r_g^* \eta),
\end{split}
\end{equation}
where $r_g$ denotes the right action of $g \in G$.
The multiplication of two compatible elements is then defined as
\begin{equation}
(x,g,\eta) \cdot (xg, h, r_{h^{-1}}^* \eta) = (x, gh, r_{h^{-1}}^* \eta).
\end{equation}

\begin{rmk} This groupoid is actually an action groupoid in disguise.  The action of $G$ on $M$ extends naturally to an action of $G$ on $T^*M$, and the groupoid $s^*(T^*M) \arrows T^*M$ is isomorphic to the action groupoid $T^*M \times G \arrows T^*M$ via the diffeomorphism $s^*(T^*M) \to T^*M \times G$, $(x,g, \eta) \mapsto (x, r_g^* \eta, g)$.  
\end{rmk}

We may identify the space $\left(s^*(T^*M)\right)^{(2)}$ of compatible pairs with the pullback bundle $m^* s^* (T^*M)$, consisting of elements of the form $(x,g,h,\eta)$. where $\eta \in T^*_{xgh} M$.  Then the three face maps to $s^*(T^*M)$ are
\begin{equation}\begin{split}\label{eqn:gspacem}
\widetilde{p_1}(x,g,h,\eta) &= (x,g,r_h^* \eta), \\
\widetilde{m}(x,g,h,\eta) &= (x, gh, \eta), \\
\widetilde{p_2}(x,g,h,\eta) &= (xg, h, \eta). \\
\end{split}\end{equation}

\subsection{The $\mathcal{LA}$-groupoid}
Putting together the structures of \S\ref{sub:gspacealgbd} and \S\ref{sub:gspacegpd}, we may form the square
\begin{equation}\label{eqn:gspacela}
\la{s^*(T^*M)}{M \times G}{T^*M}{M,}
\end{equation}
where the vertical sides are groupoids and the horizontal sides (with the exception of the top side, if $s$ is not Poisson) are algebroids.  Since the structure maps (\ref{eqn:gspacest}) and (\ref{eqn:gspacem}) are linear over the corresponding structure maps for the action groupoid $M \times G \arrows M$, we may apply the $[-1]$ functor to the left side of (\ref{eqn:gspacela}) to get the groupoid $s^*([-1]T^*M) \arrows [-1]T^*M$ equipped with the vector field $d_{\widetilde{\pi}}$.

The structure maps for $s^*([-1]T^*M) \arrows [-1]T^*M$, which by abuse of notation we will also refer to as $\widetilde{s}$, $\widetilde{t}$, and $\widetilde{m}$, are described as follows.  For $X \in \vect{M}$, 
\begin{equation}\begin{split}
\widetilde{s}^* X = 1 \otimes X,\\
\widetilde{t}^* X = s_* \widetilde{X},
\end{split}\end{equation}
where $\widetilde{X}$ is the horizontal lift of $X$ to a vector field on $M \times G$.  For a ``vector field along $M \times G$ in $M$'' $\chi \in \Gamma(s^*(TM))$, the operator $m^* \circ \chi$ is a section of $m^*s^*(TM)$, which may be viewed as a linear function on $(s^*([-1]TM))^{(2)}$.  The multiplication map is then simply described by
\begin{equation}
\widetilde{m}^* \chi = m^* \circ \chi.
\end{equation}

\begin{lemma}\label{lemma:bcrels}
$d_{\widetilde{\pi}}$ is $\widetilde{s}$-related to $d_\alpha$ if and only if $s$ is a Poisson map.
\end{lemma}
\begin{proof}
Recall that $s$ is a Poisson map if and only if, for all $f \in C^\infty(M)$, $[\alpha, f]$ is $s$-related to $[\widetilde{\pi}, s^*f]$.  We then have that $d_{\widetilde{\pi}} \widetilde{s}^* f = s_* [\widetilde{\pi}, s^* f]$, which is equal to $1 \otimes [\alpha, f] = \widetilde{s}^* d_\alpha f$ if and only if $s$ is a Poisson map.

The ``only if'' part of the statement already follows, but for the ``if'' part we need to check that $d_{\widetilde{\pi}} \widetilde{s}^* X = \widetilde{s}^* d_\alpha X$ for $X \in \vect(M)$.  Under the hypothesis that $s$ is Poisson, we may use (\ref{eqn:bcdiffvf}) and see that $d_{\widetilde{\pi}} \widetilde{s}^* X = d_{\widetilde{\pi}} (1 \otimes X) = 1 \otimes d_\alpha X = \widetilde{s}^* d_\alpha X$.
\end{proof}

\begin{thm}
The square (\ref{eqn:gspacela}) is an $\mathcal{LA}$-groupoid if and only if $s: M \times G \rightarrow M$ is a Poisson map.  Equivalently, $d_{\widetilde{\pi}}$ is multiplicative if and only if $s$ is a Poisson map.
\end{thm}

\begin{proof}
We will prove the latter formulation of the theorem.  Since the candidate base field for $d_{\widetilde{\pi}}$ is the Poisson cohomology operator $d_\alpha \in \vect([-1]T^*M)$, we have that $d_{\widetilde{\pi}}$ is multiplicative if and only if all of the following hold:
\begin{enumerate}
\item $d_{\widetilde{\pi}}$ is $\widetilde{s}$-related to $d_\alpha$,
\item $d_{\widetilde{\pi}}$ is $\widetilde{t}$-related to $d_\alpha$,
\item $d_{\widetilde{\pi}}$ is $\widetilde{m}$-related to $d_{\widetilde{\pi}}^{(2)}$.
\end{enumerate}
We have already seen in Lemma \ref{lemma:bcrels} that the first statement holds if and only if $s$ is a Poisson map.  To complete the proof of the theorem, it remains to show, if $s$ is Poisson, that the second and third statements hold.

(\emph{$\widetilde{t}$-relatedness}).
Let $f \in C^\infty(M)$.  Then, using the identity $\widetilde{[\alpha,f]} = [\widetilde{\pi}, t^*f]$ (this is essentially the statement that $t$ is a Poisson map), we compute that $\widetilde{t}^* d_\alpha f = s_* \widetilde{[\alpha,f]} = s_* [\widetilde{\pi}, t^* f] = d_{\widetilde{\pi}} \widetilde{t}^* f$.  

Let $X \in \vect(M)$.  Using (\ref{eqn:bcdiffpush}), we similarly compute that $d_{\widetilde{\pi}} \widetilde{t}^* X = d_{\widetilde{\pi}}(s_* \widetilde{X}) = s_* [\widetilde{\pi}, \widetilde{X}] = s_* \widetilde{[\alpha, X]} = \widetilde{t}^* d_\alpha X$.  Thus $d_{\widetilde{\pi}}$ is $\widetilde{t}$-related to $d_\alpha$.

(\emph{$\widetilde{m}$-relatedness}).
Since $d_{\widetilde{\pi}}$ is $\widetilde{s}$- and $\widetilde{t}$-related to $d_\alpha$, there exists a unique lift $d_{\widetilde{\pi}}^{(2)}$ to $(s^*([-1]T^*M))^{(2)} = m^*s^*([-1]T^*M)$, satisfying the equations
\begin{equation}\begin{split}
d_{\widetilde{\pi}}^{(2)} (g \otimes 1) &= s_*m_*[\widetilde{\pi}^{(2)}, g],\\
d_{\widetilde{\pi}}^{(2)} (1 \otimes X) &= 1 \otimes d_\alpha X,
\end{split}\end{equation}
where $g \in C^\infty(M \times G^2)$, $X \in \vect(M)$, and $\widetilde{\pi}^{(2)} \defequal \alpha \times \pi^2$.  Thus, for any $f \in C^\infty(M \times G)$, we have $d_{\widetilde{\pi}}^{(2)} \widetilde{m}^* (f \otimes 1) = s_*m_*[\widetilde{\pi}^{(2)}, m^*f]$.  Since $\pi$ is a multiplicative vector field on $G$, this is equal to $m^* \circ (s_*[\widetilde{\pi},f]) = \widetilde{m}^* d_{\widetilde{\pi}} (f \otimes 1)$.

Now suppose $X \in \vect(M)$.  Then $d_{\widetilde{\pi}}^{(2)} \widetilde{m}^* (1 \otimes X) = d_{\widetilde{\pi}}^{(2)} (1 \otimes X) = 1 \otimes d_\alpha X$, whereas $\widetilde{m}^* d_{\widetilde{\pi}}(1 \otimes X) = \widetilde{m}^* (1 \otimes d_\alpha X) = 1 \otimes d_\alpha X$.  Thus $d_{\widetilde{\pi}}$ is $\widetilde{m}$-related to $d_{\widetilde{\pi}}^{(2)}$.
\end{proof}

This result suggests several further questions regarding the $\mathcal{LA}$-groupoid (\ref{eqn:gspacela}).  For example, Lu \cite{lu} has described a basic subcomplex of the complex $\left(\bigwedge \Gamma(A^*), d_A \right)$, where $A = (M \times \mathfrak{g}) \oplus T^*M$, whose cohomology is equal to the tensor product of the $G$-invariant de Rham cohomology and the $G$-invariant Poisson cohomology of $M$.  The double complex arising from the $\mathcal{LA}$-groupoid (\ref{eqn:gspacela}) should compute the same data.

In light of the work of Bursztyn and Crainic \cite{bu-cr}, it would also be interesting to see how this construction generalizes to quasi-Poisson actions.  It could result in the notion of a ``quasi-double''.
		
		\chapter{The van Est map}
		\label{chapter:vanest}
\label{ch:vanest}

The map of van Est \cite{vanest} relates the smooth group cohomology of a Lie group $G$ to the Lie algebra cohomology of its Lie algebra $\mathfrak{g}$.  This map was extended to Lie groupoids by Weinstein and Xu \cite{wx}, and later by Crainic \cite{crainic:vanest}.  In this chapter, we describe an extension, which appears in \cite{mw}, of the van Est map to $Q$-groupoids.  We define the van Est map in a way that is suitable for supergroupoids, and show that the map is equivariant with respect to the action of multiplicative vector fields.  It follows that the van Est map gives a homomorphism of double complexes from the double complex of a $Q$-groupoid to that of its $Q$-algebroid.

As we have seen in the examples of the Weil algebra (Example \ref{ex:weil}) and the BRST model (Example \ref{ex:brst1}), the double complexes of a $Q$-groupoid and its $Q$-algebroid do not generally have the same cohomology.  It was this fact that motivated the definition of the basic subcomplex, which is not canonical.  Since the van Est map \emph{is} canonical, we would not expect it to be an isomorphism onto any basic subcomplex.  However, there may exist homotopy equivalences sending the image of the van Est map to the basic subcomplexes.  It also would be interesting to see if it is possible to describe a Chern-Weil map\footnote{Fernandes \cite{fernandes} shows that Chern-Weil theory extends in a straightforward manner to Lie algebroids, but the resulting characteristic classes are just the pullbacks by the anchor map of the usual characteristic classes for the tangent bundle.  More interesting invariants arise as ``exotic'' characteristic classes \cite{crainic:vanest} \cite{cf} \cite{fernandes}, the first of which is the \emph{modular class} \cite{elw}.} that integrates forms in the $Q$-algebroid complex and to compare such a map with the van Est map.  We hope to deal with these issues in future work.

For simplicity, the formulas in this chapter are given for the case where $G$ is an ordinary groupoid.  Nonetheless, all of the results go through for supergroupoids, with some appropriate sign modifications.  

\section{Preliminaries}

Let $G \arrows M$ be a groupoid.  Recall from \S\ref{section:cohomology} the groupoid cohomology complex of $G$.  We will be dealing with the homotopy equivalent subcomplex $C^\infty_\nu (G^{(q)})$ of \emph{normalized} cochains, consisting of the functions $f \in C^\infty(G^{(q)})$ such that $(\Delta^{q-1}_i)^* f = 0$ for all $i$.

Let $X$ be a left-invariant vector field on $G$.  For each $q$, there is a natural lift of $X$ to a vector field $X_q^q$ on $G^{(q)}$, acting as $X$ on the last component of $G^{(q)}$.  The following identities follow immediately from the definition of $X_q^q$:
\begin{align}\label{eqn:facex}
X_q^q (\sigma_i^q)^* &= \left\{ \begin{aligned} 
&(\sigma_i^q)^* X_{q-1}^{q-1}, &i& < q, \\
&0, &i& = q,
\end{aligned}\right. \\
X_q^q (\Delta_i^q)^* &= (\Delta_i^q)^* X_{q+1}^{q+1}, \begin{aligned} &&i& < q. \label{eqn:degenx} \end{aligned}
\end{align}

Consider the operator $X^q \defequal (\Delta_{q-1}^{q-1})^* X_q^q$ from $C^\infty(G^{(q)})$ to $C^\infty(G^{(q-1)})$.  The above equations and the simplicial relations (\ref{eqn:twoface}), (\ref{eqn:twodegen}), and (\ref{eqn:facedeg}) imply that
\begin{align}
X^q (\sigma_i^q)^* &= \left\{ \begin{aligned} 
&(\sigma_i^{q-1})^* X^{q-1}, &i& < q-1, \\
&X^{q-1}_{q-1}, &i& = q-1 \\
&0, &i& = q,
\end{aligned}\right. \label{eqn:facex2}\\
X^q (\Delta_i^q)^* &= (\Delta_i^{q-1})^* X^{q+1}, \begin{aligned} &&i& < q. \label{eqn:degenx2} \end{aligned}
\end{align}

\begin{lemma}
$X^q$ preserves the complex of normalized cochains.
\end{lemma}
\begin{proof}
Let $f \in C^\infty_\nu(G^{(q)})$.  Using (\ref{eqn:degenx2}), for any $i$, 
\[(\Delta_i^{q-2})^* X^q f = X^{q-1} (\Delta_i^{q-1})^*  f,\]
which vanishes because $f$ is normalized.
\end{proof}

Let $\theta \in C^\infty(M)$.  Recall that $\theta$ acts on left-invariant vector fields via $s^*$; specifically, $(\theta X)(f) = s^* \theta \cdot X(f)$.  Let $p_i^q: G^{(q)} \to G$ be the map defined by projection onto the $i$th component, and define $\theta_q^q \defequal (p_q^q)^* s^* \theta$ to be the ``lift'' of $\theta$ to the last component of $G^{(q)}$.

\begin{lemma}\label{lemma:theta}
\begin{enumerate}
\item $(\theta X)_q^q = \theta_q^q \cdot X_q^q$.
\item $(\Delta_q^q)^* (\theta_{q+1}^{q+1}) = \theta_q^q$.
\end{enumerate}
\end{lemma}

\begin{proof}The lift $X_q^q$ is completely determined by the properties that it is $p_i^q$-related to $0$ for $i < q$ and $p_q^q$-related to $X$.  Using the latter property, we compute that
\begin{equation*}\begin{split}
(\theta X)_q^q (p_q^q)^* &= (p_q^q)^* (s^*\theta \cdot X) \\
&= \theta_q^q \cdot (p_q^q)^* X \\
&= \theta_q^q \cdot X_q^q (p_q^q)^*,
\end{split}\end{equation*}
so $(\theta X)_q^q$ and $\theta_q^q \cdot X_q^q$ agree on functions in the image of $(p_q^q)^*$.  Since they also agree on functions in the image of $(p_i^q)^*$ for $i < q$ (indeed, they both vanish on such functions), we conclude by Lemma \ref{lemma:equal} that they are equal.

The second statement follows from the identity $s \circ p_{q+1}^{q+1} \circ \Delta_q^q = s \circ p_q^q$, which may be easily verified.
\end{proof}

\begin{lemma}\label{lemma:lin}
If $f \in C^\infty_\nu(G^{(q)})$, then for any $X, Y \in \livf(G)$ and $\theta \in C^\infty(M)$,
\begin{equation}
X^{q-1} (\theta Y)^q f = (\theta X)^{q-1} Y^q f.
\end{equation}
\end{lemma}
\begin{proof}Using Lemma \ref{lemma:theta},
\begin{equation}
\begin{split}
X^{q-1}(\theta Y)^q f &= X^{q-1} (\Delta_{q-1}^{q-1})^*\left[\theta_q^q \cdot Y_q^q (f) \right] \\
&= X^{q-1} \left( \theta_{q-1}^{q-1} \cdot Y^q (f) \right) \\
&= (\Delta_{q-2}^{q-2})^* \left[X_{q-1}^{q-1} (\theta_{q-1}^{q-1}) \cdot Y^q f + \theta_{q-1}^{q-1} \cdot X_{q-1}^{q-1}(Y^q f) \right].
\end{split}
\end{equation}
Since $Y^q f$ is normalized, the first term vanishes, and the remaining term may be written as $(\theta X)^{q-1} Y^q f$.
\end{proof}

\begin{lemma}\label{lemma:ring}
Let $f \in C^\infty(G^{(q)})$ and $g \in C^\infty(G^{(q')})$, where $q' > 0$.  Then
\[X^{q+q'}(f*g) = f*(X^{q'}g).\]
If $f$ is normalized, then for $\theta \in C^\infty(M)$,
\[X^q(f*\theta) = (X^q f)*\theta.\]
\end{lemma}

\begin{proof} First we observe that the product of cochains may be described in terms of the simplicial structure as
\begin{equation} \label{eqn:starface}
f * g = \left(\left(\sigma_{q+q'}^{q+q'}\right)^* \cdots \left(\sigma_{q+1}^{q+1}\right)^* f \right) \cdot \left(\left(\sigma_0^{q+q'}\right)^* \cdots \left(\sigma_0^{q'+1}\right)^* g\right).
\end{equation}
Using (\ref{eqn:facex}) and (\ref{eqn:facedeg}), we have that
\begin{equation*}\begin{split}
X^{q+q'} (f * g) &= \left(\Delta_{q+q'-1}^{q+q'-1}\right)^* X_{q+q'}^{q+q'} (f * g) \\
&= \left(\Delta_{q+q'-1}^{q+q'-1}\right)^* \left(\left(\sigma_{q+q'}^{q+q'}\right)^* \cdots \left(\sigma_{q+1}^{q+1}\right)^* f \right) \cdot \left(\left(\sigma_0^{q+q'}\right)^* \cdots \left(\sigma_0^{q'+1}\right)^* X_{q'}^{q'} g\right) \\
&= \left(\left(\sigma_{q+q'-1}^{q+q'-1}\right)^* \cdots \left(\sigma_{q+1}^{q+1}\right)^* f \right) \cdot \left(\left(\sigma_0^{q+q'-1}\right)^* \cdots \left(\sigma_0^{q'}\right)^* \left(\Delta_{q'-1}^{q'-1} \right)^* X_{q'}^{q'} g\right) \\
&= f* \left(X^{q'}g\right).
\end{split}\end{equation*}

Similarly, (\ref{eqn:starface}) implies that $f * \theta = f \cdot \theta_q^q$.  Then
\begin{equation*}
X^q(f*\theta) = \left(\Delta_{q-1}^{q-1}\right)^* \left(X_q^q f \cdot \theta_q^q + f \cdot X_q^q \theta_q^q \right).
\end{equation*}
If $f$ is normalized, then the second term vanishes, leaving $X^q f \cdot \theta_{q-1}^{q-1} = (X^q f)*\theta$.
\end{proof}

\section{Definition of the van Est map}

Let $G \arrows M$ be a groupoid with algebroid $A \to M$.

\begin{dfn} Let $f \in C^\infty(G^{(q)})$.  Then $Vf \in \bigwedge^q \Gamma(A^*)$ is defined by the equation
\begin{equation}
\iota_{X_q} \cdots \iota_{X_1} Vf = \sum_{\gamma \in S_q} (\sgn \gamma) (X_{\gamma(1)})^1 \cdots  (X_{\gamma(q)})^q f,
\end{equation}
where the sum is taken over the symmetric group $S_q$.
\end{dfn}

\begin{rmk} $Vf$ is constructed so as to be antisymmetric in the arguments $X_i$.  It follows from Lemma \ref{lemma:lin} that $Vf$ is $C^\infty(M)$-linear in each argument.  
\end{rmk}

\begin{prop} Let $f \in C^\infty_\nu(G^{(q)})$ and $g \in C^\infty_\nu(G^{(q')})$.  Then $V(f*g) = Vf \wedge Vg$.
\end{prop}
\begin{proof}Recall that, if $\alpha \in \bigwedge^q \Gamma(A^*)$ and $\beta \in \bigwedge^{q'} \Gamma(A^*)$, then 
\begin{equation} \label{eqn:wedge}
\iota_{X_{q' + q}} \cdots \iota_{X_1} \alpha \wedge \beta = \sum_{\gamma \in S_{(q',q)}} (\sgn \gamma) \left(\iota_{X_{\gamma(q'+q)}}\cdots \iota_{X_{\gamma(q'+1)}} \alpha \right)\left(\iota_{X_{\gamma(q')}}\cdots \iota_{X_{\gamma(1)}} \beta \right),\end{equation}
where the sum is taken over the $(q',q)$-shuffles.

From Lemma \ref{lemma:ring}, it follows that
\begin{multline}
\iota_{X_{q+q'}} \cdots \iota_{X_1} V(f*g) = \sum_{\gamma\in S_{q+q'}} (\sgn \gamma) (X_{\gamma(1)})^1 \cdots (X_{\gamma(q+q')})^{q+q'} (f*g) \\
= \sum_{\gamma\in S_{q+q'}} (\sgn \gamma) \left[ (X_{\gamma(1)})^1 \cdots (X_{\gamma(q)})^{q\vphantom{'}} f \right] \left[ (X_{\gamma(q+1)})^1 \cdots (X_{\gamma(q+q')})^{q'} g \right].
\end{multline}
Since there is a bijection between $S_{q+q'}$ and the product of sets $S_{(q,q')} \times S_{q'} \times S_q$, the final expression may be identified with (\ref{eqn:wedge}), where $\alpha = Vf$ and $\beta = Vg$.
\end{proof}

\begin{prop}
Let $f \in C^\infty_\nu(G^{(q)})$.  Then $V(\delta f) = d_A (Vf)$.
\end{prop}
\begin{proof}
By the definition of $\delta$,
\begin{equation}
\iota_{X_{q+1}} \cdots \iota_{X_1} V (\delta f) = \iota_{X_{q+1}} \cdots \iota_{X_1} V \left( \sum_{i=0}^{q+1} (-1)^i (\sigma_i^{q+1})^* f\right). 
\end{equation}

Consider first the $i=0$ term.  Using (\ref{eqn:facex2}) to move the face map to the left, 
\begin{equation}
\begin{split}
&\sum_{\gamma \in S_{q+1}} (\sgn \gamma)(X_{\gamma(1)})^1 \cdots (X_{\gamma(q+1)})^{q+1}(\sigma_0^{q+1})^* f\\ 
=& 
\sum_{\gamma \in S_{q+1}} (\sgn \gamma)(X_{\gamma(1)})^1 (\sigma_0^1)^* (X_{\gamma(2)})^2 \cdots (X_{\gamma(q+1)})^{q+1} f \\
=& \sum_{\gamma \in S_{q+1}} (\sgn \gamma) \rho(X_{\gamma(1)}) (X_{\gamma(2)})^2 \cdots (X_{\gamma(q+1)})^{q+1} f.
\end{split}
\end{equation}
Expressing $S_{q+1}$ as $S_{(1,q)} \times S_q$ by considering for each $i$ the set of permutations that begin with $i$, the last expression is equivalent to 
\begin{equation}\label{eqn:veanch}
\sum_{i=1}^{q+1}(-1)^i \rho(X_i) \iota_{X_{q+1}} \cdots \hat{\iota_{X_i}} \cdots \iota_{X_1} Vf.
\end{equation}

The $i=q+1$ term vanishes by (\ref{eqn:facex2}), so the remaining terms are
\begin{equation}
\sum_{\gamma \in S_{q+1}} \sum_{i=1}^q (\sgn \gamma)(-1)^i (X_{\gamma(1)})^1 \cdots (X_{\gamma(q+1)})^{q+1}(\sigma_i^{q+1})^* f
\end{equation}
Again, using (\ref{eqn:facex2}), the face maps may be moved to the left so that the expression becomes
\begin{equation}\begin{split}\label{eqn:veb1}
&\sum_{\gamma \in S_{q+1}} \sum_{i=1}^q (\sgn \gamma)(-1)^i (X_{\gamma(1)})^1 \cdots (X_{\gamma(i+1)})^{i+1} (\sigma_i^{i+1})^* (X_{\gamma(i+2)})^{i+1} \cdots (X_{\gamma(q+1)})^q f \\
=& \sum_{\gamma \in S_{q+1}} \sum_{i=1}^q (\sgn \gamma)(-1)^i (X_{\gamma(1)})^1 \cdots (X_{\gamma(i)})^i (X_{\gamma(i+1)})_i^i (X_{\gamma(i+2)})^{i+1} \cdots (X_{\gamma(q+1)})^q f.
\end{split}\end{equation}

Observe in (\ref{eqn:veb1}) that the two vector fields $(X_{\gamma(i)})_i^i$ and $(X_{\gamma(i+1)})_i^i$ appear adjacent in each term.  For a fixed $i$, consider for each $\gamma \in S_{q+1}$ the ``dual'' permutation $\bar{\gamma}$ which exchanges $\gamma(i)$ and $\gamma(i+1)$.  The permutations $\gamma$ and $\bar{\gamma}$ have opposite sign, so combining the dual terms in the sum introduces a Lie bracket:

\begin{equation}
\sum_{i=1}^q \sum_{\gamma(i)<\gamma(i+1)}  (\sgn \gamma)(-1)^i (X_{\gamma(1)})^1 \cdots [X_{\gamma(i)}, X_{\gamma(i+1)}]^i (X_{\gamma(i+2)})^{i+1} \cdots (X_{\gamma(q+1)})^q f.
\end{equation}

For each $j<k$, the set of $(i,\gamma)$ such that $\gamma(i)=j$ and $\gamma(i+1)=k$ may be identified with $S_q$, so the last expression may be written as 

\begin{equation}\label{eqn:vebrack}
\sum_{j<k} (-1)^k \iota_{X_{q+1}} \cdots \iota_{[X_j,X_k]} \cdots \hat{\iota_{X_j}} \cdots \iota_{X_0} Vf.
\end{equation}

Putting (\ref{eqn:veanch}) and (\ref{eqn:vebrack}) together, we see that $V(\delta f) = d_A (Vf)$.
\end{proof}

Let $\psi$ be a multiplicative vector field on $G$.  Then, for each $q$, there exists a lift $\psi^{(q)} \in \vect(G^{(q)})$ that acts as $\psi$ on each component.  Recall that the Lie bracket of $\psi$ with any left-invariant vector field is again left-invariant.

\begin{lemma} \label{lemma:multx}Let $\psi$ be a multiplicative vector field, and let $X$ be a left-invariant vector field.  Then for any $q$,
\begin{enumerate}
\item $[\psi^{(q)}, X^q_q] = [\psi, X]_q^q$,
\item $X^q \psi^{(q)} = \psi^{(q-1)} X^q - [\psi, X]^q$.
\end{enumerate}
\end{lemma}

\begin{proof}
The first statement is clear, since $X_q^q$, thus $[\psi^{(q)}, X^q_q]$, acts as $0$ on all the components except the last.  On the last component, $[\psi^{(q)}, X^q_q]$ acts as $[\psi, X]$.

Using the first statement, it is immediate that $X_q^q \psi^{(q)} = \psi^{(q)} X^q_q - [\psi,X]_q^q$.  Since $\psi^{(q)}$ is $\Delta_{q-1}^{q-1}$-related to $\psi^{(q-1)}$ (see Remark \ref{rmk:multface}), the second statement follows.
\end{proof}

Recall that a multiplicative vector field $\psi$ on $G$ induces a morphic vector field $\widetilde{\psi}$ on $A$, and that the corresponding vector field $\widetilde{\psi}_{1}$ on $[-1]A$ commutes with $d_A$.

\begin{prop}
Let $\psi$ be a multiplicative vector field on $G$ with base vector field $\phi$ on $M$.  Then, for any $f \in C^\infty_\nu (G^{(q)})$, $V(\psi^{(q)} f) = \widetilde{\psi}_{1} (Vf)$.
\end{prop}

\begin{proof}
Using Lemma \ref{lemma:multx}, it is a direct calculation that 
\begin{multline}
\sum_{\gamma \in S_q} (X_{\gamma(1)})^1 \cdots (X_{\gamma(q)})^q \psi^{(q)} f = 
\sum_{\gamma \in S_q} \left[ \vphantom{\sum_{i=1}^q} \psi^{(0)} (X_{\gamma(1)})^1 \cdots (X_{\gamma(q)})^q  f - \right. \\
\left. - \sum_{i=1}^q (X_{\gamma(1)})^1 \cdots (X_{\gamma(i-1)})^{i-1} ([\psi, X_{\gamma(i)}])^i (X_{\gamma(i+1)})^{i+1} \cdots (X_{\gamma(q)})^q  f \right] \\
= \phi(\iota_{X_q} \cdots \iota_{X_1} Vf) - \sum_{i=1}^q \iota_{X_q} \cdots \iota_{X_{i+1}} \iota_{[\psi,X_i]} \iota_{X_{i-1}} \cdots \iota_{X_1} Vf.
\end{multline}

Comparing with (\ref{eqn:liemorph}) and Corollary \ref{cor:multaction}, we see that this result is equal to $\iota_{X_q} \cdots \iota_{X_1} \widetilde{\psi}_{1} Vf$.
\end{proof}

The results of this section may be summarized in the following
\begin{thm}\label{thm:vanesthom}
The map $V: C^\infty_\nu (G^{(q)}) \to \bigwedge^q \Gamma(A^*)$ is a homomorphism of differential algebras and is equivariant with respect to the action of multiplicative vector fields.
\end{thm}
Applying Theorem \ref{thm:vanesthom} to a $Q$-groupoid $\mathcal{G}$ with homological vector field $\psi$, we have
\begin{cor}
The van Est map is a homomorphism of double complexes 
\begin{equation}V: (C^{p,q}_\nu(\mathcal{G}), \delta, \psi) \to \left(\Gamma_p \left(\bigwedge\vphantom{\bigwedge}^q A^*\right), d_\mathcal{A}, \widetilde{\psi}_1\right).\end{equation}
\end{cor}

As was mentioned earlier, it would be interesting to know under what conditions the van Est map induces an isomorphism on cohomology of double complexes, as well as its relationship with the Chern-Weil map.  As a more immediate application, we may use the van Est map to ``integrate'' various structures from the level of algebroids to that of groupoids.  One example is the integration of (twisted) Dirac structures to (twisted) presymplectic groupoids, carried out in \cite{bcwz}.
		\appendix
		\chapter{Simplicial structures}
		\label{appendix}
This appendix provides a brief review of some results from simplicial homotopy theory.  The main objective is to describe a double complex, associated to a simplicial structure, which computes equivariant cohomology.  The material here has been collected primarily from \cite{bott}, \cite{dupont}, and \cite{segal} (see also \cite{bss}, \cite{segal:cat}, and \cite{shulman}).

All manifolds in this section are assumed to be ordinary (i.e. not super).  It is possible to define simplicial supermanifolds in the obvious way, but the geometric realization functor is no longer sensible in the super category.  However, the reader may wish to keep in mind that the de Rham complex of a simplicial supermanifold \emph{does} make sense, so it is possible to compute the ``cohomology of the geometric realization'' even if the geometric realization does not really exist (see Remark \ref{rmk:superbg}).

\section{Simplicial spaces}

\begin{dfn} \label{dfn:simplicial}A \emph{simplicial space} is a sequence $X = \{X_q\}$, $q \geq 0$ of topological spaces equipped with \emph{face maps} $\sigma_i^q: X_q \to X_{q-1}$, $i=0,\dots,q$ and \emph{degeneracy maps} $\epsilon_i^q : X_q \to X_{q+1}$ $i=0, \dots, q$ such that
\begin{align}
\sigma_i^{q-1}\sigma_j^q &= \sigma_{j-1}^{q-1}\sigma_i^q, &i < j,\label{eqn:twoface}\\
\epsilon_i^{q+1}\epsilon_j^q &= \epsilon_{j+1}^{q+1}\epsilon_i^q, &i < j, \label{eqn:twodegen}\\
\sigma_i^{q+1}\epsilon_j^q &= \left\{ \begin{aligned}
&\epsilon_{j-1}^{q-1}\sigma_i^q, & i &< j,\\
&id, & i &= j, \; i = j+1, \\
&\epsilon_j^{q-1}\sigma_{i-1}^q, & i &> j+1. \label{eqn:facedeg}\end{aligned}
\right.
\end{align}
\end{dfn}

\begin{rmk} If $\{X_q\}$ is a sequence of manifolds equipped with smooth face and degeneracy maps satisfying the conditions of Definition \ref{dfn:simplicial}, then it is called a \emph{simplicial manifold}.  \end{rmk}

\begin{ex}\label{ex:mansimp}Let $M$ be a manifold.  Then the simplicial manifold $\bar{M}$, where $\bar{M}_q = M$ and all the face and degeneracy maps are the identity map, is such that $\left|\bar{M}\right| = M$.
\end{ex}

\begin{ex}[Nerve of a manifold]\label{ex:pairnerve}
Let $M$ be a manifold.  Then the \emph{nerve} $N\bar{M}$ of $M$ is a simplicial manifold where
\begin{equation} N\bar{M}_q = M^{q+1} \defequal \underbrace{M \times \cdots \times M}_{q+1} \end{equation}
and
\begin{align} \sigma_i^q (x_0, \dots, x_q) &= (x_0, \dots, \hat{x_i}, \dots x_q), \\
\epsilon_i^q (x_0, \dots, x_q) &= (x_0, \dots, x_i, x_i, \dots, x_q).
\end{align}
\end{ex}

\begin{ex}[Nerve of a group]
Let $G$ be a (topological) group.  Then the \emph{nerve} $NG$ of $G$ is a simplicial space where
\begin{equation} NG_q = G^q \defequal \underbrace{G \times \cdots \times G}_{q} \end{equation}
and
\begin{align} 
\sigma_0^q (g_1, \dots, g_q) &= (g_2, \dots, g_q),\\
\sigma_i^q (g_1, \dots, g_q) &= (g_1, \dots, g_i g_{i+1}, \dots g_q), & 0<i<q, \\
\sigma_q^q (g_1, \dots, g_q) &= (g_1, \dots, g_{q-1}),\\
\epsilon_i^q (g_1, \dots, g_q) &= (g_1, \dots, g_i, e, g_i+1, \dots, g_q).
\end{align}
When $q=0$, $NG_0$ is a point and the degeneracy map $\epsilon_0^0$ maps the point to the identity of $G$.
\end{ex}

\begin{rmk}The above examples are special cases of the \emph{nerve of a groupoid} (see Definition \ref{dfn:nerve}).  Example \ref{ex:mansimp} is the nerve of the trivial groupoid $M \arrows M$.  Example \ref{ex:pairnerve} is the nerve of the pair groupoid $M \times M \arrows M$.
\end{rmk}

\subsection{Geometric realization}

Let $\Delta^n$ be the standard $n$-simplex.  There are maps $\hat{\sigma}_i^q : \Delta^{q-1} \to \Delta^q$ and $\hat{\epsilon}_i^q : \Delta^{q+1} \to \Delta^q$, where
\begin{align}
\hat{\sigma}_i^q (t_0, \dots, t_{q-1}) &= (t_0, \dots, t_{i-1}, 0, t_i, \dots, t_{q-1}), \\
\hat{\epsilon}_i^q (t_0, \dots, t_{q+1}) &= (t_0, \dots, t_i + t_{i+1}, \dots, t_{q+1}).
\end{align}
The maps $\hat{\sigma}_i^q$ and $\hat{\epsilon}_i^q$ satisfy identities dual to those of Definition \ref{dfn:simplicial}.

\begin{dfn} Let $X = \{X_q\}$ be a simplicial space.  The \emph{geometric realization} of $X$ is the topological space
\begin{equation} \left| X \right|= \left( \coprod_n \Delta^n \times X_q \right) / \sim, \end{equation}
with the identifications
\begin{align}
(\hat{\sigma}_i^q(t), x) &\sim (t, \sigma_i^q(x)), \\
(\hat{\epsilon}_i^q(t), x) &\sim (t, \epsilon_i^q(x)).
\end{align}
\end{dfn}

\begin{rmk}If each $X_q$ has a CW structure and the face and degeneracy maps are cellular, then $\left|X \right|$ inherits a CW structure with an $n$-cell for each nondegenerate $(n-q)$-cell on $X_q$.  If each $X_q$ has the homotopy type of a CW complex, then $\left|X \right|$ has the homotopy type of a CW complex.\end{rmk}

\begin{dfn} Let $X = \{X_q\}$ and $X' = \{X'_q\}$ be simplicial spaces.  The \emph{product simplicial space} is $X \times X' = \{X_q \times X'_q \}$, where the face and degeneracy maps are the product maps of the respective maps on $X$ and $X'$.
\end{dfn}

\begin{prop}[\cite{segal:cat}] If $X_q$ and $X'_q$ are compactly generated spaces, then $\left|X \times X'\right|$ is naturally homeomorphic to $\left|X\right| \times \left|X'\right|$.
\end{prop}

\subsection{The de Rham complex of a simplicial manifold}\label{sec:drc}

Let $X = \{X_q\}$ be a simplicial manifold, and denote by $C^{p,q}(X)$ the space $\Omega^p(X_q)$ of differential $p$-forms on $X_q$.  Then there are differential operators 
\begin{equation} d: C^{p,q}(X) \to C^{p+1, q}(X) \end{equation}
and
\begin{equation} \delta: C^{p,q}(X) \to C^{p, q+1}(X), \end{equation}
where $d$ is the de Rham differential and
\begin{equation} \delta = \sum_{i=0}^{q+1} (-1)^i \left(\sigma_i^{q+1}\right)^*. \end{equation}
Clearly, $d$ and $\delta$ commute, so the total differential $D = d - (-1)^p \delta$ satisfies $D^2 = 0$.  Denote by $H^\bullet (X)$ the cohomology of the total complex.  

\begin{prop}[\cite{bss}, \cite{dupont}] \label{prop:geoco} There is a natural isomorphism between $H^\bullet(X)$ and the singular cohomology $H^\bullet(\left|X\right|)$. \end{prop}

\begin{rmk} A statement similar to Proposition \ref{prop:geoco} is true in the case of simplicial spaces, where singular cochains are used instead of differential forms (see \cite{bss} and \cite{dupont}).  The classifying space construction which follows works for topological groups by making the appropriate modification to the cohomological arguments. \end{rmk}

\section{Equivariant cohomology}\label{sec:equiv}

\subsection{The classifying bundle}

\begin{prop} Let $M$ be a manifold.  Then $\left|N\bar{M}\right|$ is contractible. \end{prop}
\begin{proof}
In the spectral sequence of the de Rham complex of $N\bar{M}$, the first page has entries $E^{p,q}_1 = \ker \delta / \im \delta$.  Choose a basepoint $x_* \in M$ and define maps $\tau_q: M^q \to M^{q+1}$ by $\tau_q (x_1, \dots, x_q) = (x_1, \dots, x_q, x_*)$.  The following identities are satisfied:
\begin{align}
\sigma^q_i \circ \tau_q &= \tau_{q-1} \circ \sigma_i^{q-1}, & i < q, \;\; q > 1,\\
\sigma^q_q \circ \tau_q &= id.
\end{align}
It follows that, when $q > 0$, $\delta \circ \tau^* - \tau^* \circ \delta$ is equal (up to sign) to the identity map, and therefore $E^{p,q}_1 = 0$.

For the case $q=0$, observe that $\sigma_0^1 \circ \tau_1: M \to M$ is the map $x \mapsto x_*$.  Then $\tau_1^* \circ \delta = -id + \ev_{x_*}$, whose kernel only consists of constant functions.  Thus $E^{0,0}_1 = \reals$ and $E^{p,q}_1 = 0$ for all other values of $p$ and $q$.  The spectral sequence clearly terminates at this stage.  Since the (reduced) cohomology of the total complex is trivial, it follows from Whitehead's theorem\footnote{Note that for the analogous statement to hold for topological spaces, we must restrict to those that have the homotopy type of a CW complex.} that $\left|N\bar{M}\right|$ is contractible.
\end{proof}

Let $G$ be a Lie group.  There is a simplicial action of $G$ on $N\bar{G}$, where
\begin{equation} g \cdot (g_0, \dots , g_q) = (gg_0, \dots, gg_q). \end{equation}
The action is free in the sense that for all $q$ the action of $G$ on $N\bar{G}_q$ is free.  Thus there is an induced free action of $G$ on $\left|N\bar{G}\right|$.  The simplicial map $\gamma: N\bar{G} \to NG$, where
\begin{equation} \gamma(g_0, \dots , g_q) = (g_0^{-1}g_1, \dots, g_{q-1}^{-1} g_q), \end{equation}
describes an isomorphism of the quotient $N\bar{G} / G$ and $NG$.  There is an induced isomorphism $|\gamma| : \left|N\bar{G}\right| / G \to \left|NG\right|$.

Since $\left|N\bar{G}\right|$ is a contractible space on which $G$ acts freely, the quotient map $\left|N\bar{G}\right| \to \left|NG\right|$ may be interpreted as the classifying bundle $EG \to BG$.

\subsection{The simplicial manifold of a group action}

Let $M$ be a manifold, and let $G$ be a Lie group that acts (from the right) on $M$.  Recall the Borel model of equivariant cohomology, in which the equivariant cohomology $H_G^\bullet (M)$ is defined to be the (singular) cohomology of the homotopy quotient $M \times EG / G$.  

There is an associated free left action on the simplicial manifold $\bar{M} \times N\bar{G}$, where
\begin{equation} g \cdot (x, g_0, \dots , g_q) = (xg^{-1}, g g_0, \dots , g g_q). \end{equation}
The simplicial map $\gamma_M: \bar{M} \times N\bar{G} \to \bar{M} \times NG$, where
\begin{equation} \gamma_M (x, g_0, \dots , g_q) = (xg_0, g_0^{-1} g_1, \dots , g_{q-1}^{-1} g_q), \end{equation}
describes an isomorphism of the quotient $(\bar{M} \times N\bar{G}) / G$ and the nerve $N(M \times G)$ of the action groupoid (see Definition \ref{dfn:nerve}).  Thus $N(M \times G)$ is a simplicial manifold whose geometric realization is the homotopy quotient.  By Proposition \ref{prop:geoco}, we have

\begin{thm} Let $M$ be a manifold with a right $G$-action.  The equivariant cohomology $H^\bullet_G(M)$ is naturally isomorphic to $H^\bullet (N(M \times G))$.
\end{thm}
    
\bibliographystyle{habbrv}
\bibliography{bibio}
\end{document}